\documentclass[12pt]{article}
\usepackage{amssymb,bm,float,booktabs, diagbox, authblk}
\usepackage[colorlinks,citecolor=blue]{hyperref}
\usepackage{natbib,setspace,lscape,longtable}
\usepackage{epsfig,graphicx,abstract}
\usepackage{amsmath,amsthm,color}
\usepackage{extarrows,multirow,booktabs}
\usepackage{threeparttable}
\usepackage{url,epstopdf}
\RequirePackage[mathlines, displaymath]{lineno}
\usepackage{fancyhdr}
\usepackage{lastpage}
\usepackage{dsfont}

\usepackage{caption}
\usepackage{subcaption}
\usepackage{rotating}
\usepackage{graphics}
\usepackage{verbatim}

\usepackage{enumitem}
\setenumerate[1]{itemsep=0pt,partopsep=0pt,parsep=\parskip,topsep=0pt}
\setitemize[1]{itemsep=0pt,partopsep=0pt,parsep=\parskip,topsep=0pt}
\setdescription{itemsep=0pt,partopsep=0pt,parsep=\parskip,topsep=0pt}

\def\beqr{\begin{eqnarray}}
\def\eeqr{\end{eqnarray}}
\def\beqrs{\begin{eqnarray*}}
\def\eeqrs{\end{eqnarray*}}
\def\bep{\begin{prop}}
\def\eep{\end{prop}}

\def\eop{\hfill $\Box$ \\}

\def\mR{\mathbb{R}}
\newcommand{\D}{D}
\newcommand{\A}{A}
\newcommand{\Y}{Y}
\newcommand{\X}{X}
\newcommand{\K}{K}
\newcommand{\W}{W}

\newcommand{\x}{x}
\newcommand{\oo}{o}
\newcommand{\OO}{O}
\newcommand{\N}{N}
\newcommand{\1}{\mathds{1}}
\newcommand{\I}{I}
\newcommand{\be}{\bm{\beta}}
\newcommand{\ttop}{^{\top}}
\newcommand{\m}{m}
\newcommand{\Z}{Z}
\newcommand{\CX}{\mathcal{X}}
\newcommand{\al}{\bm{\alpha}}
\newcommand{\calS}{\mathcal{S}}

\newcommand{\var}{\text{var}}

\newcommand{\ORCATE}{\rm{ORCATE}}
\newcommand{\PRCATE}{\rm{PRCATE}}

\newcommand{\NRCATE}{\rm{NRCATE}}

\newcommand{\indep}{\;\, \rule[0em]{.03em}{.67em} \hspace{-.25em}
\rule[0em]{.65em}{.03em} \hspace{-.25em}
\rule[0em]{.03em}{.67em}\;\,}

\newtheorem{theo}{\bf Theorem}
\newtheorem{lemm}{\bf Lemma}
\newtheorem{prop}{\bf Proposition}
\newtheorem{remark}{\bf Remark}
\newtheorem{corollary}{\bf Corollary}

\allowdisplaybreaks

\title{ Outcome regression-based estimation of conditional average treatment effect}
\author[1]{Lu Li}
\author[2]{Niwen Zhou}
\author[3]{Lixing Zhu \thanks{Corresponding author:
lzhu@hkbu.edu.hk. The first two authors are co-first authors. The research was supported by a grant from The University Grants Council of Hong Kong.}}
\affil[1]{{\small{School of Finance and statistics, East China Normal University, Shanghai 200241, China}}}\affil[2]{{\small{School of statistics, Beijing Normal University, Beijing 100875, China}}}\affil[3]{{\small{Department of Mathematics, Hong Kong Baptist University, Kowloon Tong 999077, Hong Kong, China}}}

\date{}
\begin{document}
  \maketitle
		\begin{abstract}The research is about  a systematic investigation on the following issues.   First, we  construct  different outcome regression-based  estimators for conditional average treatment effect under, respectively, true (oracle), parametric, nonparametric and semiparametric dimension reduction structure.
		Second, according to the corresponding asymptotic variance functions, we  answer the following questions  when supposing the models are correctly specified: what is the asymptotic efficiency ranking  about the four estimators in general? how is the efficiency related to the affiliation of the given covariates in the set of arguments of the regression functions?   what do the roles of bandwidth and kernel function selections play for the estimation efficiency; and  in which scenarios should the estimator under semiparametric dimension reduction regression structure  be used in practice? As a by-product, the results show that any outcome regression-based estimation should be asymptotically more  efficient than any inverse probability weighting-based estimation. All these results give a relatively complete picture of the outcome regression-based estimation such that the theoretical conclusions could provide guidance for practical use when more than one estimations can be applied to the same problem. Several simulation studies are conducted to examine the performances of these estimators in finite sample cases and a real dataset is analyzed for illustration.
	\end{abstract}

\textbf{Keywords:} Asymptotic variance; Conditional average treatment effect; Regression casual effect;  Sufficient dimension reduction.

\section{Introduction}

Causal inference has been widely applied for decades to analyse treatment effect based on observational studies, in which treatments are assigned to observations in a non-random fashion. In this paper, we consider  casual inference under the potential outcome framework \citep{Rubin:1974,Rosenbaum:1983} where the treatment is binary and the outcome variable in the hypothetical complete data set has two components $(\Y_{(1)}, \Y_{(0)})$. In which $\Y_{(1)}$ is the potential outcome if the individual receives treatment and $\Y_{(0)}$ is the corresponding potential outcome without treatment. As we can only observe one of $\Y_{(1)}$ and $\Y_{(0)}$, a commonly used method is to impute a reasonable value in the lieu of the missing one such as  linear regression imputation \cite{Healy:1956}, kernel regression imputation \cite{Cheng:1994} and ratio imputation \cite{Rao:1996}.

In this paper, we consider average treatment effect (ATE)  conditional on  some covariates to explore the heterogeneity of ATE \cite{Rosenbaum:1983,Rosenbaum:1985}.
Let $\X\in\mR^{p}$ be a set of covariates that collects individual's personal information and $\X_{1}\in\mR^{k}$ be a subvector of $\X$, $1\leq k<p$. Conditional average treatment effect (CATE, hereafter) is defined as $E(\Y_{(1)}-\Y_{(0)}|\X_{1})$.
To estimate this function, \cite{Abrevaya:2015} proposed  estimators that are based on inverse probability weighting (IPW, hereafter) method and concluded that, according to the asymptotic variance functions, the estimator with noparametrically estimated propensity score (NCATE) is asymptotically more efficient than the one with parametrically estimated propensity score (PCATE). The relevant conclusion is similar to that in \cite{Hahn:1998} and \cite{Hirano:2003} for  the IPW estimators of ATE. But, PCATE is proved to be  asymptotically equivalent to the one with true propensity score (OCATE). This is very different from the unconditional ATE. {\color{blue}Zhou and Zhu(2020)}\footnote{Zhou, N. W. $\&$ Zhu L. X.(2020). On IPW-based estimation of conditional average treatment effect. Submitted.} proposed an estimator with semiparametically estimated propensity score (SCATE) and gave some more detailed analysis on the asymptotic efficiency on NCATE and SCATE.

As well known, for ATE,  outcome regression-based estimation is already a popularly used methodology. Thus, methodologically, the research in this aspect is not new. However, for CATE, the problem becomes more complicated as it involves double conditional expectations on the full set  $\X$, or subset $\be^{\ttop}X$ of covariates, if the curse of dimensionality is concerned within dimension reduction framework, and the subset $\X_1$ where $\be$ is a projection matrix. Three relevant references are  \cite{Luo:2017}, \cite{Luo:2019} and \cite{Ma:2019}). To focus on the estimation efficiency issue, we in this paper do not give more details about how to work on dimension reduction and feature selection, while only consider the general setting supposing that a dimension reduction structure already exists. We then  consider  a systematic investigation on their  asymptotic properties to answer  the following questions when the model is correctly specified in parametric case.
\begin{enumerate}
\item [Q1.] When CATE is estimated under nonparametric, semiparametric, parametric and true (oracle) regression structure, what ranking of the asymptotic efficiency  can be for these estimators?
\item  [Q2.]Note that CATE is a function of $X_1$ and the set of arguments of the regression function, say $\tilde X$ that  is not necessary to be the full $X$, and thus $X_1$ is not necessary to be a strict subset of $\tilde X$. Then could the affiliation of $X_1$ to $\tilde X$ affect the asymptotic efficiency of different estimators? This issue is unique for CATE and  particularly important under semiparametic dimension reduction framework as the regression function would be a function of $\tilde X=\be^{\ttop}X$ where $\be$ is a $p \times r$ matrix with $r \ll p$ in high dimensional scenarios.
\item [Q3.] As all estimators involve nonparametric estimations for the involved conditional expectations, how could the bandwidth and kernel function affect the efficiency? This study is particularly necessary.
\item [Q4.] Comparing with the IPW-based estimation, what efficiency ranking  should be concluded?
\end{enumerate}
We will have a very brief discussion in Section~5 about the misspecified cases, globally or locally, that will be investigated in the near future, but not be touched in this paper.

Note that  CATE is
\begin{equation*}
\begin{aligned}
\tau(\x_{1})=E[(\Y_{(1)}-\Y_{(0)})|\X_{1}=\x_{1}]=E[E(\Y_{(1)}-\Y_{(0)}|\X)|\X_{1}=\x_{1}],
\end{aligned}\end{equation*}{where $E(\Y_{(1)}-\Y_{(0)}|\X)$ is the treatment effect heterogeneity. We are interested in, under unconfoundedness assumption, estimating $\tau(\x_{1})$ in this paper. To well answer the above four questions, we suggest / propose four estimators when assuming that $m_1(\X)-m_0(\X)=E(\Y_{(1)}-\Y_{(0)}|\X)$ is completely known function (ORCATE), parametric function (PRCATE) ($m_1(\X)=m_1(\X, \theta_1)$ and $m_0(\X)=m_0(\X, \theta_0)$), semiparametric function with dimension reduction structure (SRCATE) ($m_1(\X)=m_1(\beta_1^{\ttop}\X)$ and $m_0(\X)=m_0(\beta_0^{\ttop}\X)$), and  nonparametric function (NRCATE). The details will be in Section \ref{estimation}.  We  derive  the asymptotically linear representations and asymptotic normality of these estimators in various scenarios and, according to the asymptotic variance functions and using the estimators with true regression / propensity score as the benchmark, we obtain the following results to give a relatively complete picture for the asymptotic efficiencies of the four estimation methods. The following newly derived results show that the estimated CATEs have rather different asymptotic behaviors from the estimated ATEs. Let $A  \preceq B$ mean that method $A$ has smaller asymptotic variance function than method B, and $A  \cong B$ stand for the asymptotic equivalence of them when the asymptotic variance functions are equal. The results are summarised as follows.

\begin{description}
	\item A1. This is the answer for Q1 and Q4. In general, the ranking for the asymptotic efficiencies of the estimators is, together with the results about the IPW-based estimators respectively in \cite{Abrevaya:2015} and {\color{blue}Zhou and Zhu(2020)}\footnote{Zhou, N. W. $\&$ Zhu L. X.(2020). On IPW-based estimation of conditional average treatment effect. Submitted.}:
	\begin{align*}
	\overbrace{{\tiny\text{ORCATE}}\cong{\tiny\text{PRCATE}}\preceq{\tiny\text{SRCATE}}\preceq{\tiny\text{NRCATE}}}^{\text{regression-based CATE estimators}}=\overbrace{{\tiny\text{NCATE}}\preceq {\tiny\text{SCATE}}\preceq {\tiny\text{PCATE}}\cong{\tiny\text{OCATE}}}^{\text{IPW-based CATE estimators}}.
	\end{align*}
	\item A2. For Q2, we have the following results to show the importance of the affiliation of $X_1$ to $X$. Under semiparametric dimension reduction structure,  when $X_1 \subset \beta_1^{\ttop}X \cap \beta_0^{\ttop}X$ or $X_1$ is just contained in one of the sets $\beta_1^{\ttop}X$ or $\beta_0^{\ttop}X$, $$ORCATE \cong PRCATE \preceq SRCATE.$$ While when $X_1$ is not fully included in both $\beta_1^{\ttop}X$ and $\beta_0^{\ttop}X$, we have $$ORCATE \cong PRCATE \cong SRCATE.$$  Some more results are included in Section~\ref{estimation}.  Also some similar results about NRCATE and more detailed comparisons are described in Section~\ref{estimation}.
	\item A3. This answer is for Q3. When the CATE functions are smooth sufficiently, and the bandwidth and kernel function are delicately selected,   the following asymptotic equivalence among the regression-based estimators can be achieved:  $$ORCATE \cong PRCATE \cong SRCATE \cong NRCATE.$$
	\item A4. In high-dimensional scenarios, semiparametric-based estimation is often preferable because it can greatly overcome the curse of dimensionality and also avoid model misspecification. Some more detailed studies and comparisons for the asymptotic efficiency are contained in Section~\ref{estimation}.  The numerical studies in Section \ref{simulation} support this observation.
\end{description}

The rest of this article is organized as follows. Section \ref{estimation} introduces the CATE function and give the estimators  respectively under the true, parametric, nonparametic and semiparametric framework. The asymptotic properties of the proposed estimators are systematically investigated in this section. Section \ref{simulation} presents some simulation studies  to examine the performances of the estimators.  Section \ref{application} is devoted to the analysis for a real data example. Conclusions and some further research problems are briefly discussed in Section \ref{conclusion}. For the ease of presentation, we defer all technical proofs to the appendix.

\section{Estimations and their asymptotic properties}\label{estimation}
Let $\D$ be a dummy variable indicating treatment status with $\D=1$ if an individual receives treatment and $\D=0$ otherwise. We only observe $\D$, $\X$ and $\Y\equiv\D\cdot\Y_{(1)}+(1-\D)\cdot\Y_{(0)}$ in the real situation. The propensity score $p(\D=1|\X)$ is denoted by $p(\X)$. Let $\{\X_{i},\Y_{i},\D_{i}\}$, $i=1,\ldots,n$ be $n$ independent copies of $(\X,\Y,\D)$. To estimate $\tau(\x_{1})$, we suggest a two-step estimation procedure when both $g_1$ and $g_0$ are unknown. Four estimators are proposed in this paper when the regression casual effect under true (oracle), parametric, nonparametric and semiparametric dimension reduction structure (ORCATE, PRCATE, NRCATE and SRCATE) respectively.

To clearly state the estimation procedures, recall that  the function $\m_{t}(\X)$ is defined as
\begin{align*}
\m_{t}(\X)=E(\Y_{(t)}|\X),\quad t=0,1.
\end{align*} Under the unconfounderness assumption that is the conditional independence as
$$(\Y_{(0)},\Y_{(1)})\indep\D|\X ,$$
we then  first estimate  $\m_{1}(\X)-\m_{0}(\X)$ and then its conditional expectation $\tau (x_1)=E(\m_{1}(\X)-\m_{0}(\X)|X_1).$  But in semiparametric dimension reduction structure, this unconfounderness assumption will have a different formula that will be specified in Section~2. However, directly estimating $\tau (X_1)$ in terms of  $\Y_{(1)}-\Y_{(0)}$ is not feasible as it is never observed. It is naturally to use $\Y_{(1)}$ and $\Y_{(0)}$ to estimate  $\m_{1}(\X)$ and $\m_{0}(\X)$ separately.  Afterwards  $\tau(\x_{1})$ can be estimated  by a nonparametric method such as the N-W estimation \citep{Nadaraya:1964,Watson:1964}.

As for SRCATE and NRCATE, we will have to use high order kernel functions, we give the notation here. A function $\K_{1}$: $\mR^{k}\rightarrow \mR$ is a kernel of order $s_{1}$ if it integrates to one over $\mR^{k}$, and
\begin{equation*}
\int u_{1}^{p_{1}}\cdots u_{k}^{p_{k}}\K_{1}(u)du=0
\end{equation*} for all nonnegative integers $p_{1},\cdots,p_{k}$ such that $1\leq\sum_{i=1}^{k}p_{i}<s_{1}$, and it is nonzero when $\sum_{i=1}^{k}p_{i}=s_{1}$. Some regularity conditions are listed below.
\begin{description}
	\item (C1). (Strong ignorability)\begin{description}
		\item (a) (Unconfoundedness) $(\Y_{(0)},\Y_{(1)})\indep\D|\X$.
		\item (b) (Common support) For some very small $c>0$, the propensity score function $p (\cdot)$ satisfies that $c<p(\X)<1-c$.
	\end{description}
	\item (C2). (Distribution of $\X$) The support $\CX$ of the $p$-dimensional covariate $\X$ is a Cartesian product of compact intervals, and the density of $\X$, $f(\x)$, is bounded away from 0 on $\CX$.
	\item (C3). (Kernel functions) $\K_{1}(u)$ is a kernel of order $s_{1}$ that is symmetric around zero and  $s^{*}$ times continuously differentiable.
	\item (C4). (Distribution of $\X_{1}$) The density function of $\X_{1}$, $f(\x_{1})$, is bounded away from zero and infinity and $s_{1}\geq 2$ times continuously differentiable.
\end{description}
Part (a) of condition (C1) is a commonly used condition on the treatment effect, see e.g., \citep{Rosenbaum:1983,Abrevaya:2015,Luo:2017}. Moreover, part (a) of  condition (C1) is a quite strong but standard assumption in the casual inference literature. Part (b) of condition (C1) implies that there exists overlap between the treated and control observations. Conditions (C2) and (C4) are traditional conditions for nonparametric estimation in the literature \citep{Pagan:1999,Yin:2010}. Specially, condition (C3) is for high order kernel \citep{Abrevaya:2015}. It is noted that Gaussian kernel satisfies this assumption when $k=1$ and $s_{1}=2$. Furthermore, the  value $s^{*}$ relies on  the smoothness of the regression function. More specifically, $s^{*}\geq 2$ in parametric situation, while $s^{*}\geq s_{2}$ and $s^{*}\geq s_{4}$ in nonparametric and semiparametric situation, respectively.

In the following, we   study the four estimations in separate subsections and give some further analysis for SRCATE and NRCATE in another subsection.

\subsection{ ORCATE}
This estimator will serve as a benchmark to examine the performance of other estimators developed and investigated later.  Assume that  $\m_{1}(\X)-\m_{0}(\X)$ is completely known with no need of estimation. Then ORCATE can be written as \begin{equation}\label{true}
\begin{aligned}
\widehat{\tau}(\x_{1})=\frac{\frac{1}{nh_{1}^{k}}\sum\limits_{i=1}^{n}\K_{1}\left(\frac{\X_{1i}-\x_{1}}{h_{1}}\right)\{\m_{1}(\X_{i})-\m_{0}(\X_{i})\}}{\frac{1}{nh_{1}^{k}}\sum\limits_{i=1}^{n}\K_{1}\left(\frac{\X_{1i}-\x_{1}}{h_{1}}\right)},
\end{aligned}
\end{equation} The asymptotically linear representation and asymptotic normality are stated below.
{\theo\label{oracle} Suppose that assumptions (C1) through (C4) are satisfied. Then, when regression casual effect is given without estimation, for each point $\x_{1}$ in the support of $\X_{1}$, we have
	\begin{align*}
	&\sqrt{nh_{1}^{k}}\{\widehat{\tau}(\x_{1})-\tau(\x_{1})\}\\
	=&\frac{1}{\sqrt{nh_{1}^{k}}}\frac{1}{f(\x_{1})}\sum\limits_{i=1}^{n}\{\m_{1}(\X_{i})-\m_{0}(\X_{i})-\tau(\x_{1})\}\K_{1}\left(\frac{\X_{1i}-\x_{1}}{h_{1}}\right)+\oo_{p}(1),
	\end{align*}and then
	\begin{align*}
	\sqrt{nh_{1}^{k}}\left\{\widehat{\tau}(\x_{1})-\tau(\x_{1})\right\}\xrightarrow{d}\N\left(0,\frac{||\K_{1}||_{2}^{2}\sigma_{O}^{2}(\x_{1})}{f(\x_{1})}\right),
	\end{align*}where $||\K_{1}||_{2}=\{\int{\K_{1}(u)}^{2}du\}^{1/2}$, and
	\begin{align*}
	\sigma_{O}^{2}(\x_{1})= E[\{\m_{1}(\X)-\m_{0}(\X)-\tau(\x_{1})\}^{2}|\X_{1}=\x_{1}].
	\end{align*}}

\subsection{ PRCATE}
Suppose that both $\m_{1}(\X)$ and $\m_{0}(\X)$ have parametric structures with unknown parameters $\al_{1}$ and $\al_{0}$ respectively. That is, $\m_{t}(\X,\al_{t})$ are parametric  functions for $t=0,1$.  
Since each response can only be observed in a subpopulation, to get unbiased estimators of parameters $\al_{1}$ and $\al_{0}$, we use a similar method to that of \cite{Wang:2004}. Write, for $i=1,\ldots,n,$ \begin{align*}
\D_{i}\Y_{i}=\D_{i}\m_{1}(\X_{i},\al_{1})+\D_{i}\epsilon_{1i},\quad (1-\D_{i})\Y_{i}=(1-\D_{i})\m_{0}(\X_{i},\al_{0})+(1-\D_{i})\epsilon_{0i},
\end{align*}where $\epsilon_{ti}$, $t=0,1$, are random error terms, and independent of  $\X_{i}$, $i=1,\ldots,n$. Use weighted least squares  to estimate $\al_{t}$ for $t=0,1$, and then $\m_{1}(\X_{i})=\m_{1}(\X_{i},\al_{1})$ see \cite{Matloff:1981}. Write them as $\hat \al_t$ and $\widehat{\m}_{1}(\X)$.   PRCATE is then defined as:
\begin{equation}\label{parametric estimation}
\begin{aligned}
\widehat{\tau}(\x_{1})=\frac{\frac{1}{nh_{1}^{k}}\sum\limits_{i=1}^{n}\K_{1}\left(\frac{\X_{1i}-\x_{1}}{h_{1}}\right)\{\widehat{\m}_{1}(\X_{i})-\widehat{\m}_{0}(\X_{i})\}}{\frac{1}{nh_{1}^{k}}\sum\limits_{i=1}^{n}\K_{1}\left(\frac{\X_{1i}-\x_{1}}{h_{1}}\right)},
\end{aligned}
\end{equation}where
\begin{align*}
\widehat{\m}_{1}(\X_{i})=\m_{1}(\X,\widehat{\al}_{1}),\quad \widehat{\m}_{0}(\X_{i})=\m_{0}(\X,\widehat{\al}_{0}),\quad i=1,\ldots,n.
\end{align*} Assume the following additional condition:
\begin{description}
	\item (A1). (Bandwidths) $h_{1}\rightarrow 0$, $nh_{1}^{k}\rightarrow \infty$, $nh_{1}^{2s_{1}+k}\rightarrow 0$.
\end{description}
The following theorem states the asymptotic properties of $\widehat{\tau}(\x_{1})$.
{\theo\label{Parametric} Suppose that conditions (C1) through (C4) and (A1) are satisfied for $s_{1}=s^{*}+2$. Then, for each point $\x_{1}$ in the support of $\X_{1}$, we have
	\begin{align*}
	&\sqrt{nh_{1}^{k}}\{\widehat{\tau}(\x_{1})-\tau(\x_{1})\}\\
	=&\frac{1}{\sqrt{nh_{1}^{k}}}\frac{1}{f(\x_{1})}\sum\limits_{i=1}^{n}\{\m_{1}(\X_{i})-\m_{0}(\X_{i})-\tau(\x_{1})\}\K_{1}\left(\frac{\X_{1i}-\x_{1}}{h_{1}}\right)+\oo_{p}(1)\\
	\xrightarrow{d}&\N\left(0,\frac{||\K_{1}||_{2}^{2}\sigma_{P}^{2}(\x_{1})}{f(\x_{1})}\right),
	\end{align*}where
	\begin{align*}
	\sigma_{P}^{2}(\x_{1})=\sigma_{O}^{2}(\x_{1})= E[\{\m_{1}(\X)-\m_{0}(\X)-\tau(\x_{1})\}^{2}|\X_{1}=\x_{1}].
	\end{align*}}
{\remark\label{bias term}This theorem states the asymptotic equivalence between PRCATE and ORCATE in the sense that their asymptotic variance functions are identical. }

\subsection{ NRCATE}
If we do not have  prior information on the structures of $\m_{1}(\X)$ and $\m_{0}(\X)$ or we try to avoid model misspecification,  a nonparametric estimation is feasible. Similarly, we estimate $\m_{1}(\X)$ and $\m_{0}(\X)$ separately. Therefore, NRCATE  is written as
\begin{equation}\label{nonparametic estimation estimator}
\begin{aligned}
\widehat{\tau}(\x_{1})=\frac{\frac{1}{nh_{1}^{k}}\sum\limits_{i=1}^{n}\K_{1}\left(\frac{\X_{1i}-\x_{1}}{h_{1}}\right)\{\widehat{\m}_{1}(\X_{i})-\widehat{\m}_{0}(\X_{i})\}}{\frac{1}{nh_{1}^{k}}\sum\limits_{i=1}^{n}\K_{1}\left(\frac{\X_{1i}-\x_{1}}{h_{1}}\right)},
\end{aligned}\end{equation}where
{\small\begin{align*}
	\widehat{\m}_{1}(\X_{i})=\frac{\frac{1}{nh_{2}^{p}}\sum\limits_{j=1}^{n}\K_{2}\left(\frac{\X_{j}-\X_{i}}{h_{2}}\right)\Y_{1j}\1(\D_{j}=1)}{\frac{1}{nh_{2}^{p}}\sum\limits_{j=1}^{n}\K_{2}\left(\frac{\X_{j}-\X_{i}}{h_{2}}\right)\1(\D_{j}=1)}, \widehat{\m}_{0}(\X_{i})=\frac{\frac{1}{nh_{2}^{p}}\sum\limits_{j=1}^{n}\K_{2}\left(\frac{\X_{j}-\X_{i}}{h_{2}}\right)\Y_{0j}\1(\D_{j}=0)}{\frac{1}{nh_{2}^{p}}\sum\limits_{j=1}^{n}\K_{2}\left(\frac{\X_{j}-\X_{i}}{h_{2}}\right)\1(\D_{j}=0)}.
	\end{align*}}
To study the asymptotic properties of $\widehat{\tau}(\x_{1})$, we give some more conditions on the kernel function and bandwidths.
\begin{description}
	\item (A2). $\K_{2}(u)$ is a kernel of order $s_{2}\geq p$,  symmetric around zero and equal to zero outside $\prod_{i=1}^{p}[-1,1]$ with continuous $(s_{2}+1)$ order derivatives.
	\item (A3). $h_{2}\rightarrow 0$, $\frac{\log n}{nh_{2}^{p+s_{2}}}\rightarrow 0$.
	\item (A4). $h_{2}^{2s_{2}}h_{1}^{-2s_{2}-k}\rightarrow 0$, $nh_{1}^{k}h_{2}^{2s_{2}}\rightarrow 0$.
\end{description}
Conditions (A2),  (A3) and (A4) are used to  affiliate with the high order derivatives of $m_1$ and $m_0$  to ensure the asymptotic normality. The following theorem states the main theoretical results of NRCATE. For convenience,  define the following function:
\begin{align*}
\Psi_{1}(\X,\Y,\D):=\frac{\D\{\Y-\m_{1}(\X)\}}{p(\X)}-\frac{(1-\D)\{\Y-\m_{0}(\X)\}}{1-p(\X)}+\m_{1}(\X)-\m_{0}(\X).
\end{align*}
{\theo\label{NonParametric 1} Suppose that conditions (C1) through  (C4) and (A1) through (A4) are satisfied for $s^{*}\geq s_{2}\geq p $. Then, for each point $\x_{1}$, we have
	\begin{align*}
	&\sqrt{nh_{1}^{k}}(\widehat{\tau}(\x_{1})-\tau(\x_{1}))\\
	=&\frac{1}{\sqrt{nh_{1}^{k}}}\frac{1}{f(\x_{1})}\sum\limits_{i=1}^{n}[\Psi_{1}(\X_{i},\Y_{i},\D_{i})-\tau(\x_{1})]\K_{1}\left(\frac{\X_{1i}-\x_{1}}{h_{1}}\right)+\oo_{p}(1)\\
	\xrightarrow{d}&\N\left(0,\frac{||\K_{1}||_{2}^{2}\sigma_{N}^{2}(\x_{1})}{f(\x_{1})}\right),
	\end{align*}where \begin{align*}
	\sigma_{N}^{2}(\x_{1})&= E[\{\Psi_{1}(\X,\Y,\D)-\tau(\x_{1})\}^{2}|\X_{1}=\x_{1}]\\ &=\sigma_{P}^{2}(\x_{1})+E\left\{\frac{\var(\Y_{(1)}|\X)}{p(\X)}+\frac{\var(\Y_{(0)}|\X)}{1-p(\X)}\Bigg|\X_{1}=\x_{1}\right\}\\
	&\geq \sigma_{P}^{2}(\x_{1})=\sigma_{O}^{2}(\x_{1}),
	\end{align*}}the equality holds if and only if $\frac{\var(\Y_{(1)}|\X)}{p(\X)}=0$ and $\frac{\var(\Y_{(0)}|\X)}{1-p(\X)}=0$, which  rarely happen. Thus, the  inequality shows that $\NRCATE$ is asymptotically less efficient than $\PRCATE$ and $\ORCATE$.

\subsection{ SRCATE }  An obvious limitation of NRCATE is  its incapability of handling  models with high-dimensional  covariates $X$ in practice. Therefore, how to alleviate the curse of dimensionality is an important issue. To this end, reducing  dimensionality is a natural idea. But we restrict ourselves to the sufficient dimension reduction framework below and use existing methods to estimate the projection directions as the focus of this paper is on asymptotics of the estimations assuming the dimension reduction structure is specified in a semiparametric manner. For other dimension reduction issues, we can see the relevant references such as \cite{Luo:2017} and \cite{Ma:2019}.

We first give a very brief review on sufficient dimension reduction. For  given $\be\ttop\X$ where $\be$ is a $p\times r$ orthonormal matrix with an unknown number $r\ll p$ of columns, suppose that   the regression of a response variable $\W$ is independent of $\X$, which is written as $E(\W|\X)\indep\X|\be\ttop\X$, where $\indep$ stands for independence. It is generally known that $E(\W|\X)$ is an unspecified function of $\be\ttop\X$, which allows full freedom in the regression with $\be\ttop\X$ being the sufficiently reduced covariates (from $p$ to $r$). This structure has a dimension reduction structure with unknown parameter $\be$ and also is  very much flexible with a nonparametric nature.    To identify the projection directions $\be$, \cite{Cook:2002} defined the notion of central mean subspace that is the intersection of all subspaces spanned by any $\be$ such that the above conditional independence holds.
To be specific, without notational confusion, write
$\calS_{E(\Y_{(1)}|\X)}$ and $\calS_{E(\Y_{(0)}|\X)}$  respectively spanned by $\be_{1}\in\mR^{p\times r(1)}$ and $\be_{0}\in\mR^{p\times r(0)}$  where $r(t)<p$ for $t=0,1$ as the central mean subspaces such that
\begin{equation}\label{conditional independence}
\begin{aligned}
\m_{1}(\X)\indep\X|\be_{1}\ttop\X,\quad \m_{0}(\X)\indep\X|\be_{0}\ttop\X.
\end{aligned}\end{equation}

There are some approaches available in the literature to identify $\be_1$ and $\be_0$. For instance, \cite{Luo:2017} and \cite{Ma:2019}  discussed the relevant dimension reduction issues  and derived the properties of ATE under semiparametric structures. As the focus of this paper is on the asymptotic properties  of CATE estimations and the comparisons amongst  them, we then do not give the details about the estimation procedures of dimension reduction matrices $\be_1$ and $\be_0$, while just assume the  root-$n$ consistency of two estimators $\widehat \be_1$ and $\widehat \be_0$ we can define.

Note that under this dimension reduction structure, we have   $\m_{t}(\X)= E(\Y_{(t)}|\X)=E(\Y_{(t)}|\be_{t}\ttop\X)=\m_{t}(\be_{t}\ttop\X)$ for $t=0,1$. Define a  SRCATE  as
{\small\begin{equation}\label{semiparametic estimation estimator}
	\begin{aligned}
	\widehat{\tau}(\x_{1})=\frac{\frac{1}{nh_{1}^{k}}\sum\limits_{i=1}^{n}\K_{1}\left(\frac{\X_{1i}-\x_{1}}{h_{1}}\right)\{\widehat{\m}_{1}(\widehat{\be}_{1}\ttop\X_{i})-\widehat{\m}_{0}(\widehat{\be}_{0}\ttop\X_{i})\}}{\frac{1}{nh_{1}^{k}}\sum\limits_{i=1}^{n}\K_{1}\left(\frac{\X_{1i}-\x_{1}}{h_{1}}\right)},
	\end{aligned}\end{equation}}where
In order to derive theoretical results, give the following conditions.
\begin{description}
	\item (A5). $\K_{4}(u)$ is a kernel of order $s_{4}$, is symmetric around zero, is equal to zero outside $\prod_{i=1}^{p}[-1,1]$, and is continuously differentiable. The density function of $\be_{t}\ttop\X$, $f_{t}(\be_{t}\ttop\X)$ is $s_{4}$ times continuously differentiable for $t=0,1$. For $t=0,1$, $p(\be_{t}\ttop\X)\in(c^{*},1-c^{*})$ almost surely for some $c^{*}\in(0,0.5)$.
	\item (A6). $h_{4}\rightarrow 0$, $\frac{\log n}{nh_{4}^{\max\{r(0),r(1)\}+s_{4}}}\rightarrow 0$.
	\item (A7). $h_{4}^{2s_{4}}h_{1}^{-2s_{4}-k}\rightarrow 0$, $nh_{1}^{k}h_{4}^{2s_{4}}\rightarrow 0$.
	\item (A8). $\widehat{\be}_{1}-\be_{1}=\OO_{p}(n^{-\frac{1}{2}})$ and $\widehat{\be}_{0}-\be_{0}=\OO_{p}(n^{-\frac{1}{2}})$.
\end{description}Since the treatment effect heterogeneity under the semiparametric structure is based on $\be_{t}\ttop\X$ for $t=0,1$, Assumptions (A5) through (A7) play the same role as Assumptions (A2) through (A4). Condition (A8) often holds.

Define three functions as\begin{align}\label{phi}
\Psi_{2}(\X,\Y,\D)&=\frac{\D\{\Y-\m_{1}(\X)\}}{p(\be_{1}\ttop\X)}+\m_{1}(\X)-\m_{0}(\X),\nonumber\\
\Psi_{3}(\X,\Y,\D)&=-\frac{(1-\D)\{\Y-\m_{0}(\X)\}}{1-p(\be_{0}\ttop\X)}+\m_{1}(\X)-\m_{0}(\X),\\
\Psi_{4}(\X,\Y,\D)&=\frac{\D\{\Y-\m_{1}(\X)\}}{p(\be_{1}\ttop\X)}-\frac{(1-\D)\{\Y-\m_{0}(\X)\}}{1-p(\be_{0}\ttop\X)}+\m_{1}(\X)-\m_{0}(\X).\nonumber
\end{align}Next, for ease of explanation of our theoretical results, we introduce some  notations. Write $A$ and $B$ as two sets of elements. Without confusion, write $card (A)$ as the cardinality of the set $A$.
\begin{description}
	\item (F1) $A\subset B$ stands for $A\cap B=A$. In other words, elements of $A$ are all in $B$ and $card (B)\geq card (A)$.
	\item (F2) $A\subset^{k-q}B$ stands for $A\cap B=C$ with $card (C)=k-q$, that is, $k-q$ elements of $A$ belong to $B$. When $k=q$, it means that $A$ and $B$ do not share the same elements, i.e. $A\cap B=\emptyset$, written as $A\not\subseteq B$.
\end{description}

The following theorem states some very detailed investigation on the asymptotic efficiency of SRCATE.
{\theo\label{SemiParametric 1} Suppose that assumptions (C1) through (C4), (A1) and (A5) through (A8) are satisfied for $s^{*}\geq s_{4}\geq \max\{r(0),r(1)\}$. Then, for each point $\x_{1}$ in the support of $\X_{1}$, noting the definitions of $\Psi_i$ for $i=2,3,4$ in (\ref{phi}),\\
	(1) when $\X_{1}\subset^{k-q}\be_{1}\ttop\X$ and $\X_{1}\subset^{k-q}\be_{0}\ttop\X$ with $s_{4}(2-k/q)+k>0$ and $0<q\leq k$, the asymptotically linear representation of $\widehat{\tau}(\x_{1})$ is
	\begin{align*}
	&\sqrt{nh_{1}^{k}}\{\widehat{\tau}(\x_{1})-\tau(\x_{1})\}\\
	=&\frac{1}{\sqrt{nh_{1}^{k}}}\frac{1}{f(\x_{1})}\sum\limits_{i=1}^{n}\{\m_{1}(\X_{i})-\m_{0}(\X_{i})-\tau(\x_{1})\}\K_{1}\left(\frac{\X_{1i}-\x_{1}}{h_{1}}\right)+\oo_{p}(1),
	\end{align*}and the asymptotic distribution of $\widehat{\tau}(\x_{1})$ is
	\begin{align*}
	\sqrt{nh_{1}^{k}}(\widehat{\tau}(\x_{1})-\tau(\x_{1}))\xrightarrow{d}\N\left(0,\frac{||\K_{1}||_{2}^{2}\sigma_{S,1}^{2}(\x_{1})}{f(\x_{1})}\right);
	\end{align*}
	(2) when $\X_{1}\subset\be_{1}\ttop\X$ and $\X_{1}\subset^{k-q}\be_{0}\ttop\X$ with $s_{4}(2-k/q)+k>0$ and $0<q\leq k$, the asymptotically linear representation of $\widehat{\tau}(\x_{1})$ is
	\begin{align*}
	&\sqrt{nh_{1}^{k}}\{\widehat{\tau}(\x_{1})-\tau(\x_{1})\}\\
	=&\frac{1}{\sqrt{nh_{1}^{k}}}\frac{1}{f(\x_{1})}\sum\limits_{i=1}^{n}\{\Psi_{2}(\X_{i},\Y_{i},\D_{i})-\tau(\x_{1})\}\K_{1}\left(\frac{\X_{1i}-\x_{1}}{h_{1}}\right)+\oo_{p}(1)\\ \xrightarrow{d} &\N\left(0,\frac{||\K_{1}||_{2}^{2}\sigma_{S,2}^{2}(\x_{1})}{f(\x_{1})}\right);
	\end{align*}
	(3) when $\X_{1}\subset^{k-q}\be_{1}\ttop\X$ and $\X_{1}\subset\be_{0}\ttop\X$ with $s_{4}(2-k/q)+k>0$ and $0<q\leq k$, the asymptotically linear representation of $\widehat{\tau}(\x_{1})$ is
	\begin{align*}
	&\sqrt{nh_{1}^{k}}\{\widehat{\tau}(\x_{1})-\tau(\x_{1})\}\\
	=&\frac{1}{\sqrt{nh_{1}^{k}}}\frac{1}{f(\x_{1})}\sum\limits_{i=1}^{n}\{\Psi_{3}(\X_{i},\Y_{i},\D_{i})-\tau(\x_{1})\}\K_{1}\left(\frac{\X_{1i}-\x_{1}}{h_{1}}\right)+\oo_{p}(1)\\ \xrightarrow{d} &\N\left(0,\frac{||\K_{1}||_{2}^{2}\sigma_{S,3}^{2}(\x_{1})}{f(\x_{1})}\right);
	\end{align*}
	(4) when $\X_{1}\subset\be_{1}\ttop\X$ and $\X_{1}\subset\be_{0}\ttop\X$, the asymptotically linear representation of $\widehat{\tau}(\x_{1})$ is
	\begin{align*}
	&\sqrt{nh_{1}^{k}}\{\widehat{\tau}(\x_{1})-\tau(\x_{1})\}\\
	=&\frac{1}{\sqrt{nh_{1}^{k}}}\frac{1}{f(\x_{1})}\sum\limits_{i=1}^{n}\{\Psi_{4}(\X_{i},\Y_{i},\D_{i})-\tau(\x_{1})\}\K_{1}\left(\frac{\X_{1i}-\x_{1}}{h_{1}}\right)+\oo_{p}(1)\\ \xrightarrow{d} &\N\left(0,\frac{||\K_{1}||_{2}^{2}\sigma_{S,4}^{2}(\x_{1})}{f(\x_{1})}\right);
	\end{align*}
	where
	\begin{align}\label{sor}
	\sigma_{S,1}^{2}(\x_{1})&=\sigma_{O}^{2}(\x_{1})= E[\{\m_{1}(\X)-\m_{0}(\X)-\tau(\x_{1})\}^{2}|\X_{1}=\x_{1}],\nonumber\\
	\sigma_{S,2}^{2}(\x_{1})&= E[\{\Psi_{2}(\X,\Y,\D)-\tau(\x_{1})\}^{2}|\X_{1}=\x_{1}],\nonumber\\
	\sigma_{S,3}^{2}(\x_{1})&= E[\{\Psi_{3}(\X,\Y,\D)-\tau(\x_{1})\}^{2}|\X_{1}=\x_{1}],\\
	\sigma_{S,4}^{2}(\x_{1})&= E[\{\Psi_{4}(\X,\Y,\D)-\tau(\x_{1})\}^{2}|\X_{1}=\x_{1}]. \nonumber
	\end{align}}
{\remark\label{interpretation}These results imply that the asymptotic behaviours of $\widehat{\tau}(\x_{1})$ rely on whether $\X_{1}$ is a subset of $\be_{t}\ttop\X$ for $t=0,1$. Note that $\X_{1}\subset^{k-q}\be_{t}\ttop\X$ implies that only $k-q$ elements of $\X_{1}$ are also the $k-q$ linear combinations of $\be_{t}\ttop\X$ for $t=0,1$. In this case, write $\be_{t}\ttop\X$  as $\be_{t}\ttop\X=(\X_{1(1)},\ldots,\X_{1(k-q)},(\widetilde{\be}_{t}\ttop\X)\ttop)\ttop$ for $t=0,1$. Therefore,  when $\X_{1}\subset^{k-q}\be_{t}\ttop\X$ with $s_{4}(2-k/q)+k>0$ and $0<q\leq k$,  we should determine the intersection between $\X_{1}$ and $\be_{t}\ttop\X$, and then estimate $\be_{t}$ through estimating $\widetilde{\be}_{t}$ for $t=0,1$. It could be done by using partial sufficient dimension reduction (e.g. \citet{Feng:2013}). As this is not the focus of this paper, we then assume that $\be_{t}$ can be estimated at the rate $1/\sqrt n$ of convergence. Obviously, the assumption $s_{4}(2-k/q)+k>0$ is satisfied for $k= 1$. 
}

\begin{corollary}\label{SemiParametric 2} We have
	\begin{align*} \sigma_{S,1}^{2}(\x_{1})&=\sigma_{P}^{2}(\x_{1})=\sigma_{O}^{2}(\x_{1}),\\ \sigma_{S,2}^{2}(\x_{1})&=\sigma_{P}^{2}(\x_{1})+E\left\{\frac{\var(\Y_{(1)}|\X)}{p(\be_{1}\ttop\X)}\Bigg|\X_{1}=\x_{1}\right\}\geq \sigma_{P}^{2}(\x_{1})=\sigma_{O}^{2}(\x_{1}),\\ \sigma_{S,3}^{2}(\x_{1})&=\sigma_{P}^{2}(\x_{1})+E\left\{\frac{\var(\Y_{(0)}|\X)}{1-p(\be_{0}\ttop\X)}\Bigg|\X_{1}=\x_{1}\right\}\geq \sigma_{P}^{2}(\x_{1})=\sigma_{O}^{2}(\x_{1}),\\ \sigma_{S,4}^{2}(\x_{1})&=\sigma_{P}^{2}(\x_{1})+E\left\{\left[\frac{\var(\Y_{(1)}|\X)}{p(\be_{1}\ttop\X)}+\frac{\var(\Y_{(0)}|\X)}{1-p(\be_{0}\ttop\X)}\right]\Bigg|\X_{1}=\x_{1}\right\}\geq \sigma_{P}^{2}(\x_{1})=\sigma_{O}^{2}(\x_{1}).
	\end{align*} Assume that $\var(\Y_{(t)}|\X)$ is a measurable function with respect to $\be_{t}\ttop\X$ for $t=0,1$. Then
	\begin{align*}
	E\left\{\frac{\var(\Y_{(1)}|\X)}{p(\be_{1}\ttop\X)}\right\}\leq E\left\{\frac{\var(\Y_{(1)}|\X)}{p(\X)}\right\},\quad \text{and} \quad E\left\{\frac{\var(\Y_{(0)}|\X)}{1-p(\be_{0}\ttop\X)}\right\}\leq E\left\{\frac{\var(\Y_{(0)}|\X)}{1-p(\X)}\right\}.
	\end{align*} Then
	\begin{align} \sigma_{O}^{2}(\x_{1})=\sigma_{P}^{2}(\x_{1})&\leq\sigma_{S,2}^{2}(\x_{1})\leq\sigma_{S,4}^{2}(\x_{1})\leq\sigma_{N}^{2}(\x_{1}),\nonumber\\ \sigma_{O}^{2}(\x_{1})=\sigma_{P}^{2}(\x_{1})&\leq\sigma_{S,3}^{2}(\x_{1})\leq\sigma_{S,4}^{2}(\x_{1})\leq\sigma_{N}^{2}(\x_{1}).
	\end{align}
\end{corollary}

{\remark\label{parametric and semiparametric} The results in the above corollary are based on some elementary calculations and the application of Theorem 3 of \cite{Luo:2017}. We then omit the detailed calculations.
	Based on these facts,  SRCATE is more efficient than NRCATE in all cases, and less efficient than PRCATE and ORCATE in cases (2) to (4). In particularly, SRCATE shares the same asymptotic distribution as PRCATE and ORCATE in case (1). Furthermore, SRCATE in case (4) is less efficient than cases (2) and (3).} 
\subsection{Further studies on NRCATE and SRCATE}
Inspired by Theorem \ref{SemiParametric 1} about the importance of affiliation of $X_1$ to the set of arguments of the regression functions, we further investigate $SRCATE$ and $NRCATE$ in more general settings. The results are stated in the following.
{\corollary\label{NonParametric 2} Suppose that conditions (C1) through (C4) and (A1) through (A8) are satisfied. Assume that there is a given $\widetilde{\X}$ such that $(\Y_{(0)},\Y_{(1)})\indep\X|\widetilde{\X}$ with $\widetilde{\X}\subset\X$ and $\X_{1}\not\subset\widetilde{\X}$, when $\X_{1}\subset^{k-q}\be_{1}\ttop\X$ and $\X_{1}\subset^{k-q}\be_{0}\ttop\X$ with $s_{4}(2-k/q)+k>0$ and $0<q\leq k$, then the four outcome regression-based CATE estimators share the same asymptotic distribution.\\
	Here, $\widetilde{\sigma}_{N}^{2}(\x_{1})\equiv E[\{\m_{1}(\X)-\m_{0}(\X)-\tau(\x_{1})\}^{2}|\X_{1}=\x_{1}]=\sigma_{P}^{2}(\x_{1})=\sigma_{O}^{2}(\x_{1})$.}
{\remark\label{parametric and nonparamtric 1} 
	Much to our surprise, NRCATE can be asymptotically more efficient  in this special case to share the same asymptotic variance of PRCATE. This shows the importance of covariate affiliation to the set of arguments of the regression function.  This is a unique property for CATE as for ATE, this does not happen.}

{\corollary\label{Comparallel to general case} In Theorem \ref{NonParametric 1} and Theorem \ref{SemiParametric 1}, if commonly used constraints on the bandwidths $h_{1}$, $h_{2}$ and $h_{4}$ are replaced with\\ $\sqrt{nh_{1}^{k}}\left(h_{2}^{s}+\sqrt{\log(n)/nh_{2}^{p}}\right)=\oo(1)$ and  $\sqrt{nh_{1}^{k}}\left(h_{4}^{s}+\sqrt{\frac{\log(n)}{nh_{4}^{\max\{r(0),r(1)\}}}}\right)=\oo(1)$ for some order $s$, NRCATE and SRCATE have the same asymptotic distribution as PRCATE and ORCATE.}
{\remark\label{Convergences rate}As mentioned above, if we choose the bandwidth to satisfy the above conditions, NRCATE and SRCATE will share the same asymptotic efficiencies as PRCATE and ORCATE. It is obvious that the condition \\ $\sqrt{nh_{1}^{k}}\left(h_{2}^{s}+\sqrt{\log(n)/nh_{2}^{p}}\right)=\oo(1)$ and $\sqrt{nh_{1}^{k}}\left(h_{4}^{s}+\sqrt{\frac{\log(n)}{nh_{4}^{\max\{r(0),r(1)\}}}}\right)=\oo(1)$ are much stronger than the assumptions in Theorem \ref{NonParametric 1} and Theorem \ref{SemiParametric 1}. However, it is possible to choose such bandwidths if the regression casual effect function is sufficiently smooth such that high order kernel can be used. For details, see \cite{Li:2007} and {\color{blue}{Zhou and Zhu, 2020}}.  Therefore, we obtain that the ranking for the asymptotic efficiencies of four regression-based CATE estimators and four propensity score-based CATE estimators   under the condition that $\sqrt{nh_{1}^{k}}\left(h_{2}^{s}+\sqrt{\log(n)/nh_{2}^{p}}\right)=\oo(1)$ and $\sqrt{nh_{1}^{k}}\left(h_{4}^{s}+\sqrt{\frac{\log(n)}{nh_{4}^{\max\{r(0),r(1)\}}}}\right)=\oo(1)$,    \begin{align}\label{or_ipw}
	\overbrace{{\tiny\text{ORCATE}}={\tiny\text{PRCATE}}={\tiny\text{SRCATE}}={\tiny\text{NRCATE}}}^{\text{regression-based CATE estimators}}\leq\overbrace{{\tiny\text{NCATE}}={\tiny\text{SCATE}}={\tiny\text{PCATE}}={\tiny\text{OCATE}}}^{\text{IPW-based CATE estimators}}.
	\end{align}The equality occurs if and only if {\small\begin{align*}
		E\left\{\left[\frac{\var(\Y_{(1)}|\X)}{p(\X)}+\frac{\var(\Y_{(0)}|\X)}{1-p(\X)}+p(\X)(1-p(\X))\left(\frac{\m_{1}(\X)}{p(\X)}+\frac{\m_{0}(\X)}{1-p(\X)}\right)^{2}\right]\Bigg|\X_{1}=\x_{1}\right\}=0.
		\end{align*}}In other words, regression based estimators are always more efficient than IPW-type estimators in this general setting.
	
	On the other hand,  the above investigations are mainly for theoretical studies, and in practice, we may avoid to choose those bandwidths as they are often very difficult to properly select otherwise, the estimators would perform worse.}

\section{Simulations}\label{simulation}
To verify our theoretical results, we in this section conduct simulation studies to compare the regression-based ORCATE, PRCATE, SRCATE, NRCATE estimators with IPW-based OCATE, PCATE, SCATE, NCATE estimators \citep{Abrevaya:2015}. Set $p=\text{dim}(\X)\in\{2,4\}$ to avoid the curse of dimensionality under nonparametric estimation. Based on our experience and the theoretical results, when $p$ is large, NRCATE is very hard to implement. As well known, bandwidth selection plays an important role in the NW estimation. Hence, we first discuss this issue.

\subsection{Bandwidth and kernel function selection}\label{bandwidth selection}
Note that ORCATE and PRCATE only involve one bandwidth $h_{1}$ used in the second step of the estimation procedure. We first check how to choose  bandwidth sequences and kernel functions satisfying the conditions A1 - A7. To this end, consider
\begin{equation}\label{criterion of bandwidth}
\begin{aligned}
h_{1}&=a_{1}\cdot n^{-\frac{1}{k+2_{s_{1}}-\delta_{1}}},\quad a_{1}>0,\quad \delta_{1}>0,\\
h_{2}&=a_{2}\cdot n^{-\frac{1}{p+s_{2}+\delta_{2}}},\quad a_{2}>0,\quad\delta_{2}>0,\\
h_{4}&=a_{3}\cdot n^{-\frac{1}{\max\{r(0),r(1)\}+s_{4}+\delta_{3}}},\quad a_{3}>0,\quad\delta_{3}>0,
\end{aligned}\end{equation}where $\delta_{1}$, $\delta_{2}$ and $\delta_{3}$ can be selected as small as necessary or desired. It is clear that $h_{1}$, $h_{2}$ and $h_{4}$ satisfy conditions A1, A2, A3, A5 and A6. To satisfy condition A4, we  set the kernel orders as $s_{2}=p$ and $p+1$ for even  and odd $p$ respectively; and $s_{1}=s_{2}+2$. To satisfy condition~7, under semiparametric  dimension reduction structure,  set $s_{4}=\max\{r(0),r(1)\}$ and $=\max\{r(0),r(1)\}+1$ respectively for even  and odd $\max\{r(0),r(1)\}$. Based on the above values of $s_1$, $s_2$ and $s_4$, we  verify th first parts of conditions A4 and A7. Next, consider the second parts of these two conditions. Note that
when $s_{2}\geq p$ and $s_{4}\geq\max\{r(0),r(1)\}$, {\small\begin{align*}
	-\frac{2s_{2}}{p+s_{2}}\leq-1,\quad \frac{2s_{2}+k}{2s_{2}+4+k}<1,-\frac{2s_{4}}{\max\{r(0),r(1)\}+s_{4}}\leq-1,\quad
	\frac{2s_{4}+k}{2s_{1}+k}<1.
	\end{align*}} Then \begin{align*}
-\frac{2s_{2}}{p+s_{2}}+\frac{2s_{2}+k}{2s_{2}+4+k}<0, \quad
-\frac{2s_{4}}{\max\{r(0),r(1)\}+s_{4}}+\frac{2s_{4}+k}{2s_{1}+k}<0.
\end{align*} Therefore,  $h_{2}^{2s_{2}}h_{1}^{-2s_{2}-k}\rightarrow 0$ and $h_{4}^{2s_{4}}h_{1}^{-2s_{4}-k}\rightarrow 0$.
Invoking condition A3, $nh_{1}^{k}h_{2}^{2s_{2}}=nh_{1}^{2s_{1}+k}h_{2}^{2s_{2}}h_{1}^{-2s_{1}}\rightarrow 0$ when $h_{2}^{2s_{2}}h_{1}^{-2s_{1}}\rightarrow 0$. Since $\delta_{1}$, $\delta_{2}$ and $\delta_{3}$ can be arbitrarily small, we get, because $-s_{2}/(s_{2}+p)\leq-1/2$ and $(s_{2}+2)/(2s_{2}+4+k)<1/2$,
\begin{align*}
-\frac{s_{2}}{s_{2}+p}+\frac{s_{2}+2}{2s_{2}+4+k}<0.
\end{align*} Thus, condition A4 is satisfied.  Similarly,  together with condition A6,  condition A7 can also be satisfied, which has $nh_{1}^{k}h_{4}^{2s_{4}}\rightarrow 0$ by \begin{align*}
-\frac{s_{4}}{\max\{r(0),r(1)\}+s_{4}}+\frac{s_{4}}{2s_{4}+k}<0.
\end{align*}

\subsection{Model setting}
To  examine the finite sample performances of the CATE estimators,  consider the following three models:\\
Model 1: $Y_{(0)}=0,\quad  Y_{(1)}=\X_{1}^{2}+\X_{2}+\epsilon_{1},\quad p_{1}(\X)=\frac{\exp(\X_{1}+\X_{2})}{1+\exp(\X_{1}+\X_{2})}$.\\
Model 2: $Y_{(0)}=0,\quad Y_{(1)}=\X_{1}+\X_{2}+\X_{3}+\X_{4}+\epsilon_{2},\quad p_{2}(\X)=\frac{\exp\{0.5(\X_{1}+\X_{2}+\X_{3}+\X_{4})\}}{1+\exp\{0.5(\X_{1}+\X_{2}+\X_{3}+\X_{4})\}}$.\\
Model 3: $Y_{(0)}=0,\quad Y_{(1)}=\X_{2}+\X_{3}+\epsilon_{3},\quad p_{3}(\X)=\frac{\exp\{\X_{2}+\X_{3}\}}{1+\exp\{\X_{2}+\X_{3})\}}$.
\vskip 0.1cm
Model~1 is a model with the dimensions 2 and 0 of the central mean subspaces for the treatment and control groups; Model~2 is used to verify Theorem~\ref{SemiParametric 1}. Model 3 is set to justify the theory in Corollary \ref{NonParametric 2}. The dimensions of central mean subspaces for the treatment and control group are 1 and 0 in Models~2 and 3. For Model~1, $\X=(\X_{1},\X_{2})\ttop$ is generated by
\begin{align*}
\X_{1}\sim U(-0.5,0.5),\quad \X_{2}=(1+2\X_{1})^{2}+\zeta,
\end{align*}where $\zeta\sim U(-0.5,0.5)$, $\epsilon_{1}\sim\N(0,0.25^{2})$. For Model~2, we generate $\X=(\X_{1},\X_{2},\X_{3},\X_{4})\ttop$ by
\begin{align*}
\X_{1}&\sim U(-0.5,0.5),\quad \X_{2}=1+\X_{1}^{2}+\zeta_{1},\\ \X_{3}&=(1+\X_{1})^{2}+\zeta_{2},\quad \X_{4}=(-1+\X_{1})^{2}+\zeta_{3},
\end{align*}where $\zeta_{j}\overset{iid}\sim U(-0.5,0.5)$, $\epsilon_{2}\sim\N(0,0.25^{2})$, $j=1,2,3$. In Model 3, $\X=(\X_{1},\X_{2},\X_{3})\ttop$ are given by
\begin{align*}
\X_{1}&\sim U(-0.5,0.5),\quad \X_{2}=1+\X_{1}^{2}+\vartheta_{1},\quad \X_{3}=(1+\X_{1})*(-1+\X_{1})+\vartheta_{2},
\end{align*}where $\vartheta_{j}\overset{iid}\sim U(-0.5,0.5)$, $\epsilon_{3}\sim\N(0,0.25^{2})$, $j=1,2$.

The sample size is taken to be respectively $n=200$ and $n=500$ and the replication time is 500. Let $T(\x_{1})=\sqrt{(nh_{1})}[\widehat{\tau}(\x_{1})-\tau(\x_{1})]$, we report the estimated standard deviation (SD) of $T(\x_{1})$, the BIAS of $T(\x_{1})$ and the MSE of $T(\x_{1})$. For the bandwidth selection described in Subsection \ref{bandwidth selection}, we have the following selections.\\
a). For Model~1 as   $p=2$, equation \eqref{criterion of bandwidth} gives  $s_{1}=4$, $s_{2}=2$, and $s_{4}=2$. We then choose $h_{1}=a_{1}\cdot n^{-\frac{1}{9}}$ for $a_{1}=0.05$, $h_{2}=a_{2}\cdot n^{-\frac{1}{4}}$ for $a_{2}\in\{0.5,0.4\}$, $h_{4}=a_{3}\cdot n^{-\frac{1}{4}}$ for $a_{3}\in\{0.4,0.45,0.5\}$. Here, $a_{1}$, $a_{2}$ and $a_{3}$ are called baselines.\\ b). For Model~2, as $p=4$, $h_{1}=a_{1}\cdot n^{-\frac{1}{13}}$ for $a_{1}=0.02$, $h_{2}=a_{2}\cdot n^{-\frac{1}{8}}$ for $a_{2}\in\{0.15,0.16,0.18,0.2\}$, $h_{4}=a_{3}\cdot n^{-\frac{1}{3}}$ for $a_{3}\in\{0.1,0.13,0.15\}$. \\
c). For Model~3,  
as $p=3$, then $h_{1}=a_{1}\cdot n^{-\frac{1}{13}}$ for $a_{1}=0.02$, $h_{2}=a_{2}\cdot n^{-\frac{1}{8}}$ for $a_{2}\in\{0.18,0.2\}$, $h_{4}=a_{3}\cdot n^{-\frac{1}{3}}$ for $a_{3}=\in\{0.15,0.16,0.17\}$.

To make the simulation results more accessible, we tubulate the results in Tables~1-3 and some results in Appendix, as well as plot  the SDs of all estimators divided by the SD of NRCATE to show the relative efficiency in Figures~\ref{fig2}-\ref{fig4}. We choose a Gaussian kernel and derive higher order kernels from it.
\begin{table}[h!]
	\centering
	\caption{{The distribution of $\sqrt{nh_{1}}[\widehat{\tau}(\x_{1})-\tau(\x_{1})]$ for model 1}}\resizebox{13cm}{6cm}{
		\begin{tabular}{|c|c|cccc|cccc|cccc|cccc|}
			\hline
			&       & \multicolumn{8}{c|}{n=200}                                     &        \multicolumn{8}{|c|}{n=500} \\
			\hline
			&    $\x_{1}$   & \multicolumn{1}{c}{OR} & \multicolumn{1}{c}{PR} & \multicolumn{1}{c}{SR} & \multicolumn{1}{c}{NR} & \multicolumn{1}{c}{N} & \multicolumn{1}{c}{S} & \multicolumn{1}{c}{P} & \multicolumn{1}{c}{O} &       \multicolumn{1}{c}{OR} & \multicolumn{1}{c}{PR} & \multicolumn{1}{c}{SR} & \multicolumn{1}{c}{NR} & \multicolumn{1}{c}{N} & \multicolumn{1}{c}{S} & \multicolumn{1}{c}{P} & \multicolumn{1}{c|}{O} \\
			\hline
			panel1   &       & \multicolumn{16}{c|}{ $h_{1}=0.05n^{-1/9}$, $h_{4}=0.6n^{-1/4}$, $h_{2}=0.5n^{-1/4}$} \\
			\hline
			& -0.4  & 0.187 & 0.221 & 0.218 & 0.213 & 0.363 & 0.375 & 0.397 & 0.399 & 0.191 & 0.222 & 0.214 & 0.217 & 0.386 & 0.395 & 0.415 & 0.419 \\
			& -0.2  & 0.203 & 0.217 & 0.210 & 0.215 & 0.381 & 0.390 & 0.399 & 0.405 & 0.182 & 0.192 & 0.179 & 0.195 & 0.349 & 0.357 & 0.367 & 0.368 \\
			\multicolumn{1}{|l|}{SD} & 0     & 0.193 & 0.201 & 0.213 & 0.213 & 0.446 & 0.467 & 0.471 & 0.480 & 0.192 & 0.202 & 0.211 & 0.213 & 0.404 & 0.415 & 0.453 & 0.466 \\
			& 0.2   & 0.196 & 0.204 & 0.238 & 0.236 & 0.430 & 0.440 & 0.468 & 0.496 & 0.195 & 0.204 & 0.230 & 0.227 & 0.410 & 0.420 & 0.453 & 0.479 \\
			& 0.4   & 0.197 & 0.213 & 0.241 & 0.239 & 0.394 & 0.415 & 0.443 & 0.437 & 0.200 & 0.225 & 0.243 & 0.241 & 0.392 & 0.395 & 0.443 & 0.446 \\
			\hline
			& -0.4  & -0.001 & 0.000 & -0.046 & -0.004 & 0.012 & 0.032 & 0.017 & 0.024 & 0.006 & -0.008 & -0.123 & -0.023 & -0.025 & 0.004 & -0.007 & -0.011 \\
			& -0.2  & 0.016 & 0.013 & 0.102 & 0.067 & -0.007 & 0.008 & 0.008 & 0.015 & 0.002 & 0.000 & 0.123 & 0.057 & -0.026 & -0.010 & -0.016 & -0.014 \\
			\multicolumn{1}{|l|}{BIAS} & 0     & -0.018 & -0.022 & 0.004 & 0.003 & -0.034 & -0.017 & -0.021 & -0.010 & 0.003 & 0.006 & 0.052 & 0.034 & -0.017 & -0.002 & 0.002 & 0.015 \\
			& 0.2   & 0.001 & 0.001 & 0.003 & 0.008 & 0.010 & 0.048 & 0.003 & 0.014 & -0.016 & -0.013 & -0.006 & -0.001 & 0.003 & 0.009 & 0.013 & 0.021 \\
			& 0.4   & -0.001 & 0.006 & -0.005 & -0.006 & 0.028 & 0.063 & 0.034 & 0.024 & 0.006 & 0.002 & -0.009 & -0.008 & 0.043 & 0.043 & 0.010 & 0.000 \\
			\hline
			& -0.4  & 0.035 & 0.049 & 0.049 & 0.045 & 0.132 & 0.142 & 0.158 & 0.160 & 0.037 & 0.049 & 0.061 & 0.048 & 0.149 & 0.156 & 0.172 & 0.176 \\
			& -0.2  & 0.041 & 0.047 & 0.054 & 0.051 & 0.145 & 0.152 & 0.159 & 0.165 & 0.033 & 0.037 & 0.047 & 0.041 & 0.123 & 0.128 & 0.135 & 0.135 \\
			\multicolumn{1}{|l|}{MSE} & 0     & 0.038 & 0.041 & 0.045 & 0.046 & 0.200 & 0.219 & 0.222 & 0.230 & 0.037 & 0.041 & 0.047 & 0.047 & 0.164 & 0.172 & 0.205 & 0.217 \\
			& 0.2   & 0.038 & 0.042 & 0.057 & 0.056 & 0.185 & 0.196 & 0.219 & 0.246 & 0.038 & 0.042 & 0.053 & 0.052 & 0.168 & 0.176 & 0.206 & 0.230 \\
			& 0.4   & 0.039 & 0.046 & 0.058 & 0.057 & 0.156 & 0.177 & 0.198 & 0.191 & 0.040 & 0.051 & 0.059 & 0.058 & 0.156 & 0.158 & 0.196 & 0.199 \\
			\hline
			panel2 &       & \multicolumn{16}{c|}{ $h_{1}=0.05n^{-1/9}$, $h_{4}=0.5n^{-1/4}$, $h_{2}=0.5n^{-1/4}$} \\
			\hline
			& -0.4  & 0.197 & 0.227 & 0.221 & 0.212 & 0.355 & 0.371 & 0.386 & 0.378 & 0.184 & 0.225 & 0.223 & 0.224 & 0.382 & 0.393 & 0.404 & 0.413 \\
			& -0.2  & 0.177 & 0.191 & 0.191 & 0.196 & 0.351 & 0.352 & 0.377 & 0.380 & 0.189 & 0.201 & 0.192 & 0.206 & 0.376 & 0.391 & 0.401 & 0.404 \\
			\multicolumn{1}{|l|}{SD} & 0     & 0.185 & 0.199 & 0.206 & 0.207 & 0.445 & 0.453 & 0.471 & 0.480 & 0.186 & 0.200 & 0.200 & 0.202 & 0.412 & 0.417 & 0.454 & 0.465 \\
			& 0.2   & 0.197 & 0.201 & 0.229 & 0.225 & 0.457 & 0.463 & 0.508 & 0.542 & 0.202 & 0.209 & 0.230 & 0.228 & 0.446 & 0.459 & 0.492 & 0.514 \\
			& 0.4   & 0.208 & 0.229 & 0.254 & 0.253 & 0.388 & 0.393 & 0.417 & 0.440 & 0.195 & 0.212 & 0.236 & 0.234 & 0.379 & 0.378 & 0.417 & 0.434 \\
			\hline
			& -0.4  & 0.007 & 0.007 & -0.060 & -0.004 & -0.014 & 0.011 & 0.008 & 0.002 & -0.004 & -0.014 & -0.150 & -0.034 & -0.047 & -0.017 & -0.028 & -0.029 \\
			& -0.2  & 0.010 & 0.008 & 0.095 & 0.068 & -0.029 & -0.013 & -0.009 & -0.014 & 0.011 & 0.005 & 0.127 & 0.062 & 0.006 & 0.022 & 0.026 & 0.025 \\
			\multicolumn{1}{|l|}{BIAS} & 0     & 0.014 & 0.010 & 0.044 & 0.041 & -0.007 & 0.006 & -0.002 & 0.000 & 0.007 & 0.005 & 0.058 & 0.040 & -0.010 & -0.002 & -0.007 & 0.000 \\
			& 0.2   & 0.007 & 0.001 & 0.000 & 0.004 & -0.017 & 0.003 & -0.014 & -0.007 & -0.012 & -0.010 & -0.001 & 0.000 & -0.017 & -0.013 & -0.016 & -0.009 \\
			& 0.4   & -0.001 & -0.008 & -0.018 & -0.023 & 0.014 & 0.029 & 0.007 & -0.006 & 0.027 & 0.032 & 0.021 & 0.019 & 0.064 & 0.066 & 0.049 & 0.043 \\
			\hline
			& -0.4  & 0.039 & 0.051 & 0.052 & 0.045 & 0.126 & 0.138 & 0.149 & 0.143 & 0.034 & 0.051 & 0.072 & 0.052 & 0.148 & 0.155 & 0.164 & 0.172 \\
			& -0.2  & 0.031 & 0.037 & 0.046 & 0.043 & 0.124 & 0.124 & 0.142 & 0.145 & 0.036 & 0.040 & 0.053 & 0.046 & 0.141 & 0.153 & 0.162 & 0.164 \\
			\multicolumn{1}{|l|}{MSE} & 0     & 0.034 & 0.040 & 0.044 & 0.045 & 0.198 & 0.205 & 0.222 & 0.230 & 0.035 & 0.040 & 0.043 & 0.043 & 0.170 & 0.174 & 0.206 & 0.216 \\
			& 0.2   & 0.039 & 0.040 & 0.052 & 0.051 & 0.209 & 0.214 & 0.258 & 0.294 & 0.041 & 0.044 & 0.053 & 0.052 & 0.199 & 0.211 & 0.242 & 0.265 \\
			& 0.4   & 0.043 & 0.053 & 0.065 & 0.064 & 0.151 & 0.155 & 0.174 & 0.194 & 0.039 & 0.046 & 0.056 & 0.055 & 0.148 & 0.148 & 0.176 & 0.190 \\
			\hline
	\end{tabular}}%
	\label{tab1}%
\end{table}%
\begin{figure}[h!]
	\centering
	\includegraphics[width=10cm]{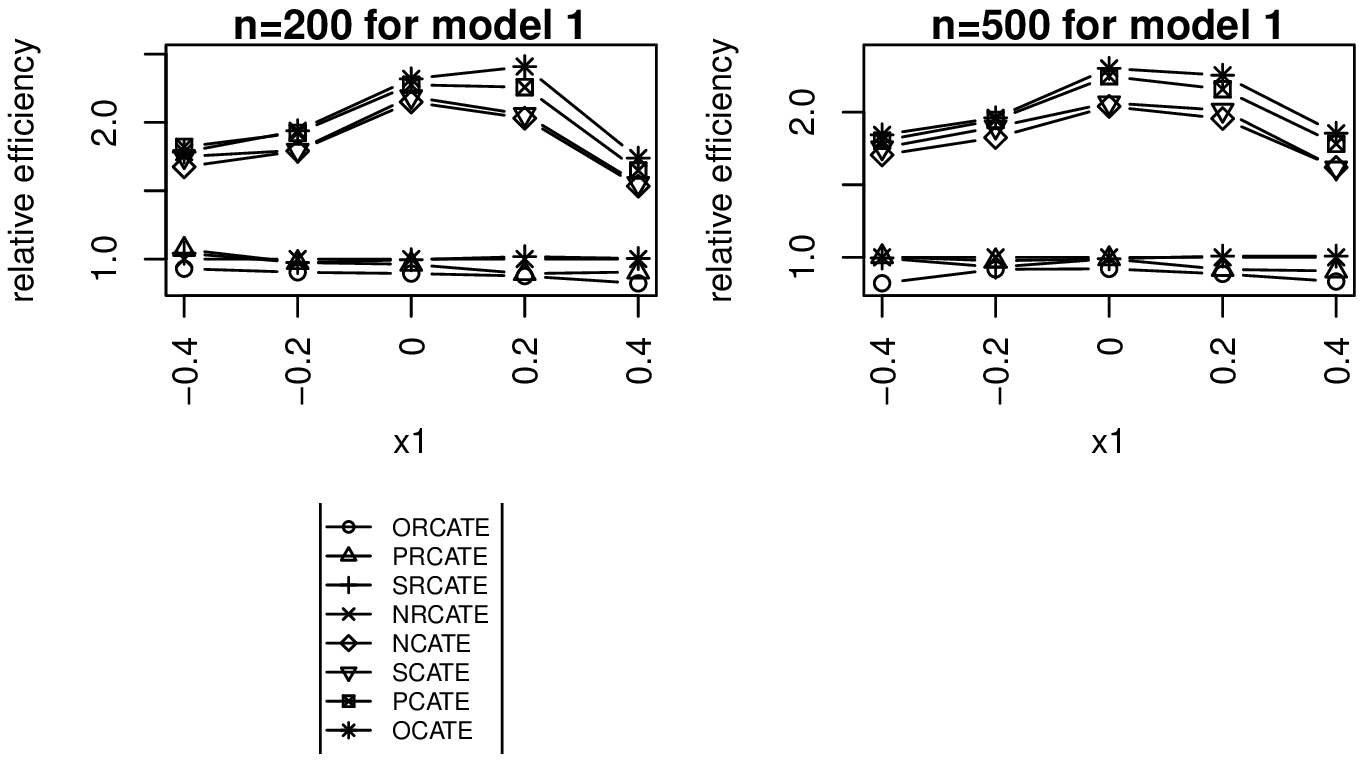}\\
	\includegraphics[width=10cm]{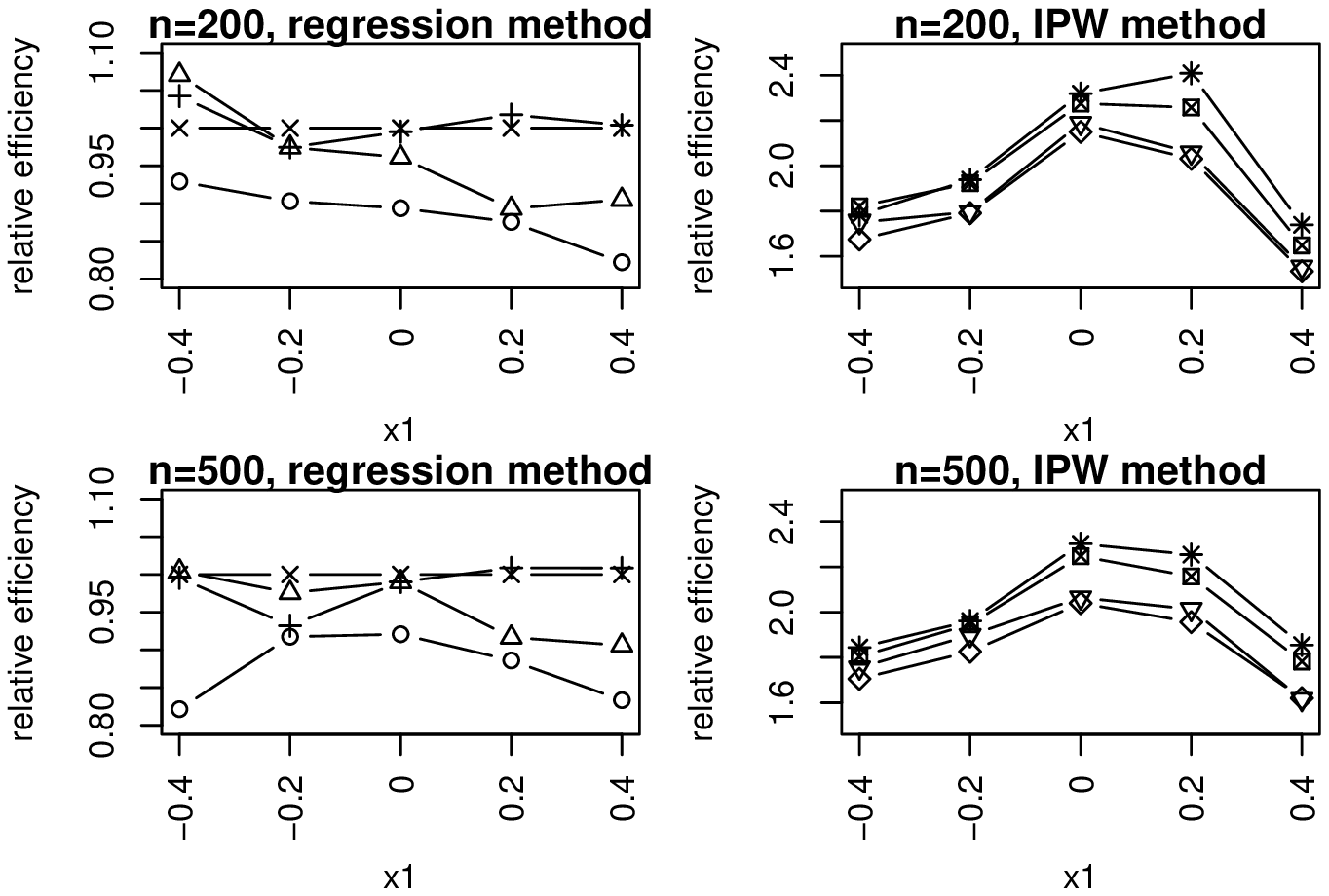}
	\caption{Relative efficiency of the  CATE estimators against NRCATE for Model 1, which are based the results in panel2 of Table~\ref{tab1}. }\label{fig2}
\end{figure}

\subsection{Simulation results}
The observations are as follows.

First, it is reasonable that larger sample size results in smaller  SD and MSE to show the estimation consistency. The dimension of $X$ also effects the estimation performance. When  $p$ increases to 4 from 2, both SD and MSE obviously increase particularly when $n=500$.

Second, the comparisons  show the significant advantage of outcome regression-based estimation over IPW-based estimation. Even though in theory, NRCATE is asymptotically equivalent to NCATE, the difference on the estimation efficiency is still very significant. All tables and figures obviously indicate this: all IPW-based estimators have much larger SD than all regression-based estimators.

Third,  as discussed before, the performances of  NRCATE and SRCATE are highly associated with the affiliation of the given covariates to the set of arguments of the outcome regression.  This finding can also be confirmed in Tables \ref{tab2} and \ref{tab3} and Figures \ref{fig3} and \ref{fig4}. In Model~2,  $\X_{1}\subset^{k-q}\be_{1}\ttop\X$ and $\X_{1}\subset^{k-q}\be_{0}\ttop\X$ with $k=1$ and $q=0$, thus in theory, SRCATE shares the same asymptotic variance as PRCATE and ORCATE and is more efficient than NRCATE. From Table~\ref{tab2} and Figure \ref{fig3}, we can see that the SDs of SRCATE are similar to those of PRCATE and ORCATE, which are smaller than that of NRCATE.  In Model~3, $\X_{1}\not\subseteq\widetilde{\X}=(\X_{2},\X_{3})\ttop$. the asymptotic efficiencies are equivalent in theory and its SDs in Table~\ref{tab3} are similar to, even slightly smaller than,  the others. In this case, all outcome regression-based estimations  have smaller SDs than all IPW-based estimations. Figure~\ref{fig4} obviously tells this.

\begin{table}[h!]
	\centering
	\caption{{The distribution of $\sqrt{nh_{1}}[\widehat{\tau}(\x_{1})-\tau(\x_{1})]$ for model 2}}\resizebox{13cm}{5cm}{
		\begin{tabular}{|c|c|cccc|cccc|cccc|cccc|}
			\hline
			&       & \multicolumn{8}{c|}{n=200}                                     &        \multicolumn{8}{|c|}{n=500} \\
			\hline
			&    $\x_{1}$   & \multicolumn{1}{c}{OR} & \multicolumn{1}{c}{PR} & \multicolumn{1}{c}{SR} & \multicolumn{1}{c}{NR} & \multicolumn{1}{c}{N} & \multicolumn{1}{c}{S} & \multicolumn{1}{c}{P} & \multicolumn{1}{c}{O} &       \multicolumn{1}{c}{OR} & \multicolumn{1}{c}{PR} & \multicolumn{1}{c}{SR} & \multicolumn{1}{c}{NR} & \multicolumn{1}{c}{N} & \multicolumn{1}{c}{S} & \multicolumn{1}{c}{P} & \multicolumn{1}{c|}{O} \\
			\hline
			panel1   &       & \multicolumn{16}{c|}{ $h_{1}=0.02n^{-1/9}$, $h_{4}=0.15n^{-1/4}$, $h_{2}=0.15n^{-1/4}$} \\
			\hline
			& -0.4  & 0.384 & 0.390 & 0.403 & 0.410 & 1.023 & 1.166 & 1.151 & 1.156 & 0.354 & 0.358 & 0.375 & 0.395 & 0.983 & 1.122 & 1.106 & 1.106 \\
			& -0.2  & 0.367 & 0.370 & 0.380 & 0.419 & 1.035 & 1.205 & 1.200 & 1.200 & 0.354 & 0.354 & 0.362 & 0.380 & 0.969 & 1.132 & 1.104 & 1.116 \\
			\multicolumn{1}{|l|}{SD} & 0     & 0.366 & 0.369 & 0.385 & 0.415 & 0.981 & 1.159 & 1.151 & 1.140 & 0.385 & 0.388 & 0.399 & 0.414 & 0.965 & 1.128 & 1.087 & 1.091 \\
			& 0.2   & 0.374 & 0.376 & 0.395 & 0.417 & 0.992 & 1.180 & 1.137 & 1.129 & 0.364 & 0.365 & 0.370 & 0.388 & 1.008 & 1.141 & 1.103 & 1.126 \\
			& 0.4   & 0.397 & 0.404 & 0.430 & 0.427 & 1.037 & 1.186 & 1.139 & 1.129 & 0.362 & 0.365 & 0.384 & 0.407 & 1.067 & 1.250 & 1.190 & 1.199 \\
			\hline
			& -0.4  & 0.014 & 0.009 & 0.056 & 0.031 & -0.692 & -0.134 & 0.048 & 0.051 & -0.017 & -0.014 & 0.082 & -0.003 & -1.069 & -0.201 & -0.010 & -0.010 \\
			& -0.2  & 0.015 & 0.012 & 0.043 & 0.021 & -0.778 & -0.207 & -0.042 & -0.034 & 0.010 & 0.012 & 0.050 & 0.016 & -1.038 & -0.198 & -0.014 & 0.001 \\
			\multicolumn{1}{|l|}{BIAS} & 0     & -0.005 & -0.008 & -0.012 & -0.001 & -0.782 & -0.191 & -0.023 & -0.027 & -0.025 & -0.025 & -0.034 & -0.011 & -1.107 & -0.243 & -0.082 & -0.074 \\
			& 0.2   & 0.004 & 0.004 & -0.021 & 0.005 & -0.652 & -0.047 & 0.062 & 0.063 & 0.017 & 0.015 & -0.034 & 0.020 & -1.059 & -0.158 & -0.058 & -0.055 \\
			& 0.4   & 0.002 & 0.003 & -0.036 & 0.002 & -0.578 & 0.045 & 0.103 & 0.091 & 0.020 & 0.016 & -0.053 & 0.005 & -0.905 & 0.049 & 0.026 & 0.000 \\
			\hline
			& -0.4  & 0.148 & 0.152 & 0.166 & 0.169 & 1.525 & 1.378 & 1.328 & 1.338 & 0.125 & 0.128 & 0.147 & 0.156 & 2.109 & 1.299 & 1.224 & 1.224 \\
			& -0.2  & 0.135 & 0.137 & 0.146 & 0.176 & 1.676 & 1.494 & 1.443 & 1.441 & 0.125 & 0.126 & 0.133 & 0.145 & 2.015 & 1.321 & 1.220 & 1.246 \\
			\multicolumn{1}{|l|}{MSE} & 0     & 0.134 & 0.136 & 0.148 & 0.173 & 1.574 & 1.380 & 1.324 & 1.301 & 0.149 & 0.151 & 0.161 & 0.171 & 2.158 & 1.331 & 1.189 & 1.195 \\
			& 0.2   & 0.140 & 0.141 & 0.157 & 0.174 & 1.408 & 1.394 & 1.296 & 1.279 & 0.133 & 0.134 & 0.138 & 0.151 & 2.137 & 1.326 & 1.219 & 1.272 \\
			& 0.4   & 0.158 & 0.163 & 0.187 & 0.183 & 1.410 & 1.410 & 1.307 & 1.283 & 0.131 & 0.133 & 0.150 & 0.166 & 1.957 & 1.566 & 1.417 & 1.437 \\
			\hline
			panel2 &       & \multicolumn{16}{c|}{ $h_{1}=0.02n^{-1/9}$, $h_{4}=0.16n^{-1/4}$, $h_{2}=0.1n^{-1/4}$} \\
			\hline
			& -0.4  & 0.385 & 0.392 & 0.397 & 0.440 & 1.066 & 1.266 & 1.236 & 1.244 & 0.375 & 0.379 & 0.385 & 0.418 & 0.946 & 1.125 & 1.066 & 1.095 \\
			& -0.2  & 0.386 & 0.387 & 0.389 & 0.430 & 0.992 & 1.159 & 1.172 & 1.176 & 0.357 & 0.361 & 0.370 & 0.397 & 0.937 & 1.116 & 1.083 & 1.102 \\
			\multicolumn{1}{|l|}{SD} & 0     & 0.379 & 0.384 & 0.387 & 0.411 & 1.031 & 1.227 & 1.203 & 1.213 & 0.387 & 0.388 & 0.399 & 0.420 & 0.933 & 1.118 & 1.103 & 1.095 \\
			& 0.2   & 0.388 & 0.387 & 0.403 & 0.443 & 1.089 & 1.254 & 1.225 & 1.247 & 0.375 & 0.376 & 0.383 & 0.407 & 1.028 & 1.198 & 1.161 & 1.177 \\
			& 0.4   & 0.376 & 0.379 & 0.392 & 0.411 & 1.056 & 1.238 & 1.143 & 1.175 & 0.362 & 0.366 & 0.385 & 0.411 & 1.047 & 1.197 & 1.126 & 1.174 \\
			\hline
			& -0.4  & 0.014 & 0.011 & 0.054 & 0.029 & -0.836 & -0.207 & -0.011 & -0.014 & 0.003 & 0.003 & 0.096 & 0.004 & -1.234 & -0.182 & 0.008 & 0.002 \\
			& -0.2  & -0.010 & -0.013 & 0.017 & 0.023 & -0.879 & -0.262 & -0.073 & -0.060 & 0.005 & 0.005 & 0.044 & 0.009 & -1.200 & -0.133 & 0.049 & 0.056 \\
			\multicolumn{1}{|l|}{BIAS} & 0     & 0.038 & 0.035 & 0.030 & 0.038 & -0.860 & -0.192 & -0.080 & -0.041 & -0.003 & -0.002 & -0.012 & 0.001 & -1.251 & -0.144 & -0.024 & -0.017 \\
			& 0.2   & 0.028 & 0.028 & 0.005 & 0.041 & -0.715 & -0.046 & 0.060 & 0.090 & -0.011 & -0.010 & -0.058 & -0.001 & -1.231 & -0.133 & -0.057 & -0.042 \\
			& 0.4   & 0.009 & 0.007 & -0.034 & 0.000 & -0.746 & -0.056 & -0.030 & -0.017 & -0.019 & -0.018 & -0.089 & -0.010 & -1.125 & -0.004 & -0.031 & -0.015 \\
			\hline
			& -0.4  & 0.148 & 0.154 & 0.161 & 0.194 & 1.836 & 1.646 & 1.529 & 1.548 & 0.140 & 0.144 & 0.157 & 0.174 & 2.418 & 1.299 & 1.137 & 1.199 \\
			& -0.2  & 0.149 & 0.150 & 0.151 & 0.186 & 1.756 & 1.411 & 1.378 & 1.387 & 0.128 & 0.131 & 0.139 & 0.157 & 2.319 & 1.262 & 1.176 & 1.218 \\
			\multicolumn{1}{|l|}{MSE} & 0     & 0.145 & 0.148 & 0.151 & 0.171 & 1.803 & 1.542 & 1.454 & 1.474 & 0.150 & 0.151 & 0.159 & 0.177 & 2.436 & 1.271 & 1.217 & 1.200 \\
			& 0.2   & 0.151 & 0.151 & 0.162 & 0.198 & 1.696 & 1.575 & 1.503 & 1.562 & 0.140 & 0.141 & 0.150 & 0.165 & 2.573 & 1.453 & 1.350 & 1.386 \\
			& 0.4   & 0.142 & 0.144 & 0.155 & 0.169 & 1.673 & 1.535 & 1.308 & 1.380 & 0.132 & 0.134 & 0.156 & 0.169 & 2.362 & 1.433 & 1.269 & 1.378 \\
			\hline
	\end{tabular}}%
	\label{tab2}%
\end{table}
\begin{figure}[h!]
	\centering
	\includegraphics[width=10cm]{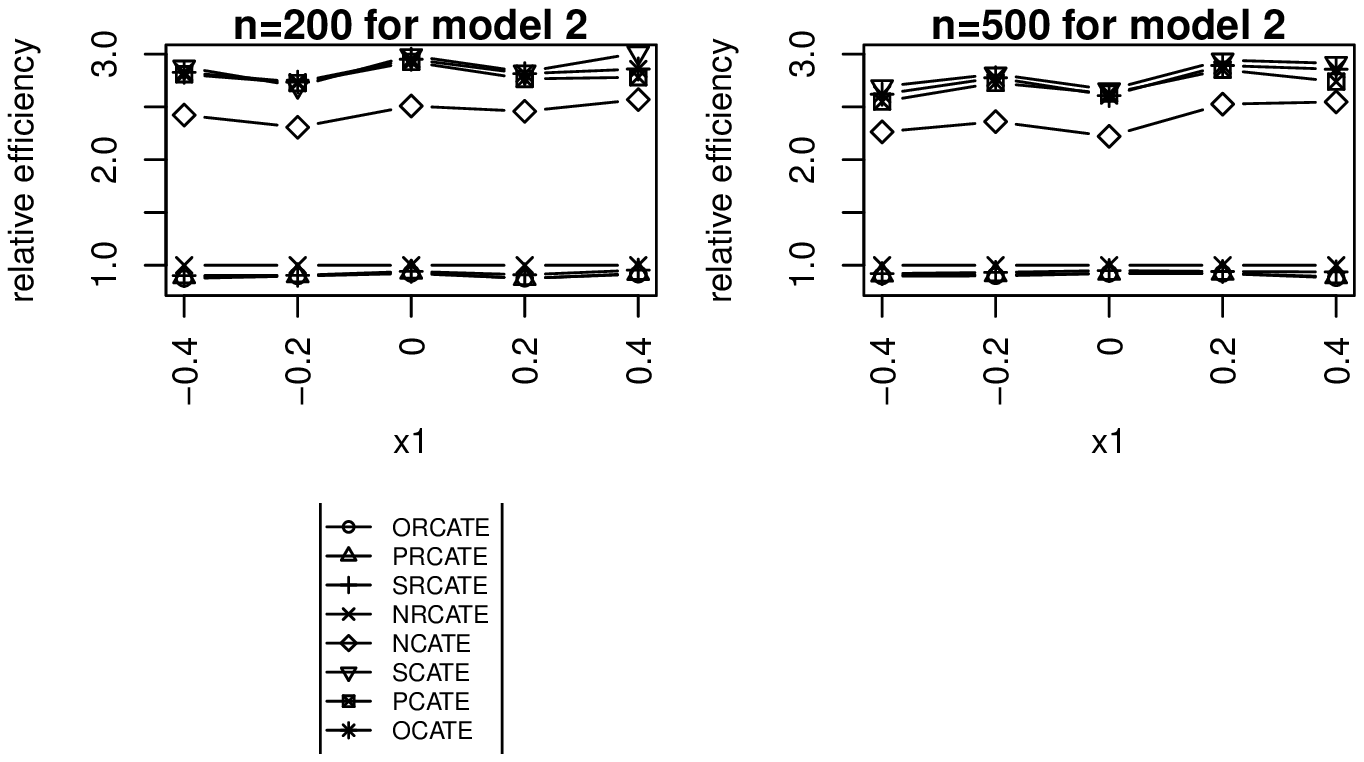}\\
	\includegraphics[width=10cm]{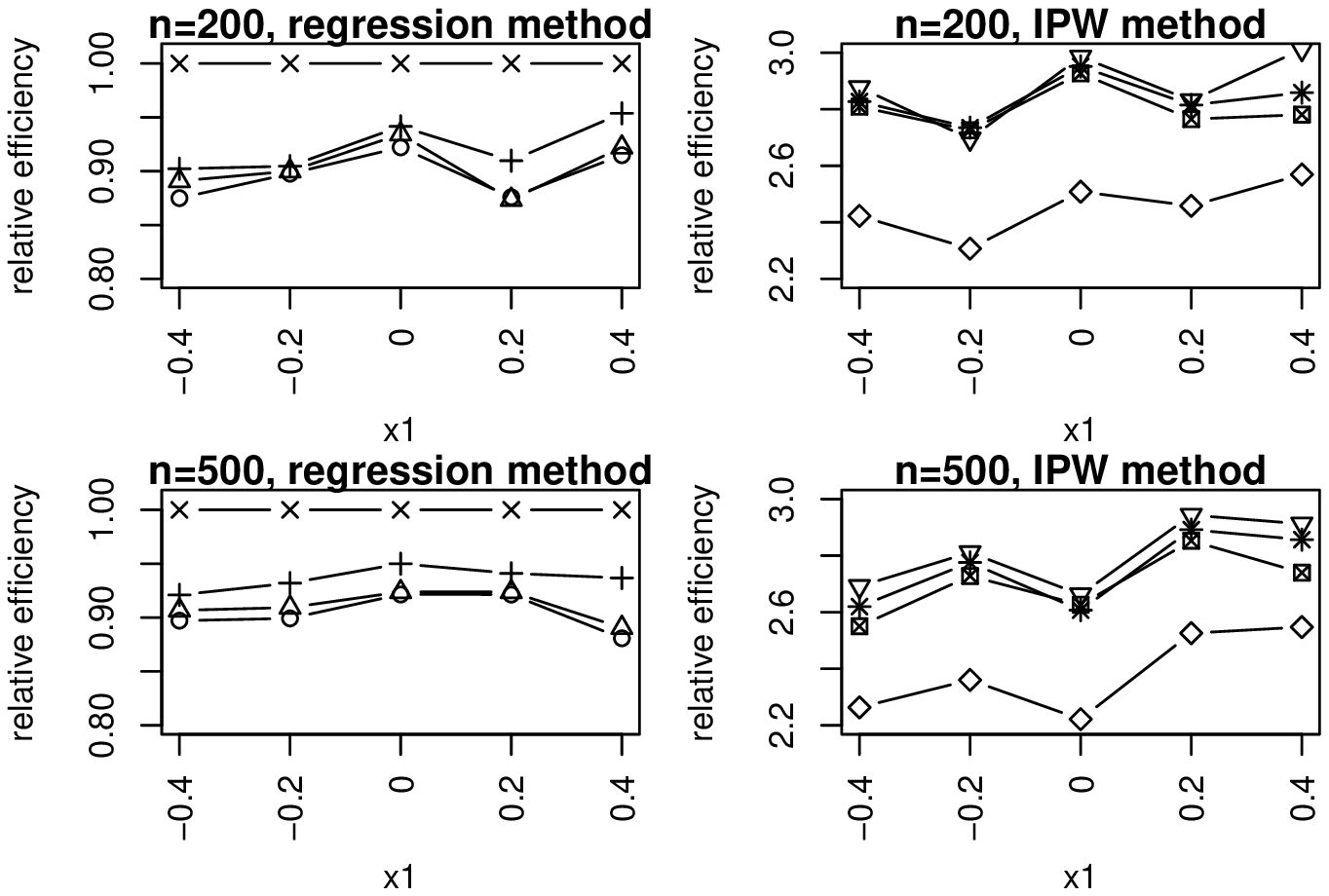}
	\caption{Relative efficiency of the CATE estimators against NRCATE for Model 2, which are based the results in panel2 of Table~\ref{tab2}.}\label{fig3}
\end{figure}

\begin{table}[h!]
	\centering
	\caption{{The distribution of $\sqrt{nh_{1}}[\widehat{\tau}(\x_{1})-\tau(\x_{1})]$ for model 3}}\resizebox{13cm}{5cm}{
		\begin{tabular}{|c|c|cccc|cccc|cccc|cccc|}
			\hline
			&       & \multicolumn{8}{c|}{n=200}                                     &        \multicolumn{8}{|c|}{n=500} \\
			\hline
			&   $\x_{1}$    & \multicolumn{1}{c}{OR} & \multicolumn{1}{c}{PR} & \multicolumn{1}{c}{SR} & \multicolumn{1}{c}{NR} & \multicolumn{1}{c}{N} & \multicolumn{1}{c}{S} & \multicolumn{1}{c}{P} & \multicolumn{1}{c}{O} &       \multicolumn{1}{c}{OR} & \multicolumn{1}{c}{PR} & \multicolumn{1}{c}{SR} & \multicolumn{1}{c}{NR} & \multicolumn{1}{c}{N} & \multicolumn{1}{c}{S} & \multicolumn{1}{c}{P} & \multicolumn{1}{c|}{O} \\
			\hline
			panel1       &       & \multicolumn{16}{c|}{ $h_{1}=0.02n^{-1/9}$, $h_{4}=0.2n^{-1/4}$, $h_{2}=0.16n^{-1/4}$} \\
			\hline
			& -0.4  & 0.327 & 0.328 & 0.330 & 0.322 & 0.498 & 0.505 & 0.538 & 0.546 & 0.282 & 0.286 & 0.287 & 0.280 & 0.481 & 0.494 & 0.492 & 0.495 \\
			& -0.2  & 0.308 & 0.310 & 0.314 & 0.310 & 0.481 & 0.480 & 0.535 & 0.530 & 0.285 & 0.286 & 0.287 & 0.282 & 0.471 & 0.474 & 0.488 & 0.486 \\
			\multicolumn{1}{|l|}{SD} & 0     & 0.301 & 0.301 & 0.306 & 0.296 & 0.452 & 0.467 & 0.514 & 0.505 & 0.287 & 0.289 & 0.294 & 0.285 & 0.479 & 0.478 & 0.506 & 0.512 \\
			& 0.2   & 0.316 & 0.319 & 0.327 & 0.317 & 0.485 & 0.500 & 0.516 & 0.516 & 0.317 & 0.317 & 0.323 & 0.314 & 0.493 & 0.492 & 0.504 & 0.500 \\
			& 0.4   & 0.290 & 0.291 & 0.301 & 0.298 & 0.485 & 0.509 & 0.514 & 0.520 & 0.297 & 0.298 & 0.299 & 0.290 & 0.476 & 0.486 & 0.493 & 0.490 \\
			\hline
			& -0.4  & -0.016 & -0.021 & -0.023 & -0.037 & -0.067 & -0.048 & -0.031 & -0.031 & -0.008 & -0.009 & -0.012 & -0.032 & -0.045 & -0.038 & -0.014 & -0.014 \\
			& -0.2  & 0.006 & 0.003 & 0.007 & 0.022 & 0.016 & 0.021 & 0.007 & 0.009 & 0.002 & 0.000 & 0.001 & 0.019 & 0.015 & 0.024 & 0.010 & 0.010 \\
			\multicolumn{1}{|l|}{BIAS} & 0     & -0.002 & -0.004 & 0.001 & 0.024 & 0.039 & 0.034 & 0.011 & 0.010 & -0.011 & -0.013 & -0.007 & 0.022 & 0.027 & 0.042 & 0.000 & -0.002 \\
			& 0.2   & 0.004 & 0.001 & 0.007 & 0.015 & -0.012 & -0.008 & -0.015 & -0.020 & 0.009 & 0.008 & 0.009 & 0.024 & 0.012 & 0.021 & 0.005 & 0.003 \\
			& 0.4   & 0.010 & 0.005 & 0.001 & -0.017 & -0.066 & -0.043 & -0.026 & -0.027 & -0.005 & -0.006 & -0.010 & -0.026 & -0.062 & -0.061 & -0.038 & -0.040 \\
			\hline
			& -0.4  & 0.107 & 0.108 & 0.109 & 0.105 & 0.252 & 0.257 & 0.290 & 0.299 & 0.080 & 0.082 & 0.083 & 0.080 & 0.234 & 0.245 & 0.243 & 0.245 \\
			& -0.2  & 0.095 & 0.096 & 0.098 & 0.097 & 0.231 & 0.231 & 0.287 & 0.281 & 0.081 & 0.082 & 0.082 & 0.080 & 0.222 & 0.225 & 0.238 & 0.236 \\
			\multicolumn{1}{|l|}{MSE} & 0     & 0.090 & 0.091 & 0.093 & 0.088 & 0.206 & 0.219 & 0.265 & 0.255 & 0.082 & 0.084 & 0.087 & 0.082 & 0.230 & 0.231 & 0.256 & 0.262 \\
			& 0.2   & 0.100 & 0.102 & 0.107 & 0.101 & 0.236 & 0.250 & 0.267 & 0.267 & 0.100 & 0.101 & 0.104 & 0.099 & 0.243 & 0.242 & 0.254 & 0.250 \\
			& 0.4   & 0.084 & 0.085 & 0.091 & 0.089 & 0.240 & 0.261 & 0.265 & 0.271 & 0.088 & 0.089 & 0.090 & 0.085 & 0.230 & 0.240 & 0.244 & 0.242 \\
			\hline
			panel2     &       & \multicolumn{16}{c|}{ $h_{1}=0.02n^{-1/9}$, $h_{4}=0.18n^{-1/4}$, $h_{2}=0.17n^{-1/4}$} \\
			\hline
			& -0.4  & 0.329 & 0.334 & 0.337 & 0.324 & 0.498 & 0.497 & 0.515 & 0.522 & 0.284 & 0.288 & 0.291 & 0.283 & 0.490 & 0.497 & 0.505 & 0.510 \\
			& -0.2  & 0.304 & 0.308 & 0.314 & 0.301 & 0.432 & 0.441 & 0.464 & 0.453 & 0.297 & 0.301 & 0.307 & 0.295 & 0.479 & 0.473 & 0.499 & 0.498 \\
			\multicolumn{1}{|l|}{SD} & 0     & 0.314 & 0.319 & 0.325 & 0.303 & 0.486 & 0.485 & 0.545 & 0.540 & 0.317 & 0.317 & 0.321 & 0.309 & 0.484 & 0.474 & 0.510 & 0.512 \\
			& 0.2   & 0.301 & 0.308 & 0.314 & 0.298 & 0.462 & 0.467 & 0.499 & 0.500 & 0.292 & 0.293 & 0.292 & 0.284 & 0.464 & 0.463 & 0.482 & 0.483 \\
			& 0.4   & 0.302 & 0.306 & 0.313 & 0.296 & 0.503 & 0.510 & 0.525 & 0.525 & 0.293 & 0.298 & 0.301 & 0.289 & 0.472 & 0.485 & 0.477 & 0.479 \\
			\hline
			& -0.4  & -0.021 & -0.016 & -0.019 & -0.027 & -0.042 & -0.036 & -0.009 & -0.010 & 0.000 & 0.002 & 0.000 & -0.019 & 0.007 & -0.002 & 0.032 & 0.034 \\
			& -0.2  & 0.004 & 0.007 & 0.015 & 0.029 & 0.019 & 0.029 & 0.012 & 0.008 & -0.015 & -0.011 & -0.011 & 0.007 & 0.008 & 0.021 & -0.002 & -0.001 \\
			\multicolumn{1}{|l|}{BIAS} & 0     & -0.012 & -0.011 & -0.009 & 0.014 & -0.007 & 0.000 & -0.042 & -0.045 & 0.006 & 0.010 & 0.012 & 0.046 & 0.021 & 0.039 & -0.006 & -0.002 \\
			& 0.2   & 0.022 & 0.023 & 0.032 & 0.039 & 0.037 & 0.048 & 0.036 & 0.035 & 0.004 & 0.008 & 0.009 & 0.028 & 0.014 & 0.025 & 0.005 & 0.008 \\
			& 0.4   & -0.023 & -0.020 & -0.022 & -0.034 & -0.040 & -0.027 & -0.009 & -0.012 & 0.003 & 0.006 & 0.007 & -0.019 & -0.004 & -0.006 & 0.016 & 0.020 \\
			\hline
			& -0.4  & 0.109 & 0.112 & 0.114 & 0.106 & 0.250 & 0.248 & 0.265 & 0.272 & 0.080 & 0.083 & 0.085 & 0.081 & 0.240 & 0.247 & 0.256 & 0.261 \\
			& -0.2  & 0.093 & 0.095 & 0.099 & 0.091 & 0.187 & 0.196 & 0.216 & 0.206 & 0.088 & 0.091 & 0.094 & 0.087 & 0.229 & 0.224 & 0.249 & 0.248 \\
			\multicolumn{1}{|l|}{MSE} & 0     & 0.099 & 0.102 & 0.106 & 0.092 & 0.236 & 0.235 & 0.299 & 0.293 & 0.101 & 0.101 & 0.103 & 0.097 & 0.234 & 0.226 & 0.260 & 0.262 \\
			& 0.2   & 0.091 & 0.095 & 0.100 & 0.090 & 0.215 & 0.221 & 0.250 & 0.251 & 0.085 & 0.086 & 0.085 & 0.082 & 0.215 & 0.215 & 0.232 & 0.234 \\
			& 0.4   & 0.092 & 0.094 & 0.098 & 0.089 & 0.255 & 0.261 & 0.276 & 0.276 & 0.086 & 0.089 & 0.091 & 0.084 & 0.223 & 0.235 & 0.228 & 0.230 \\
			\hline
	\end{tabular}}%
	\label{tab3}%
\end{table}%
\begin{figure}[h!]
	\centering
	\includegraphics[width=10cm]{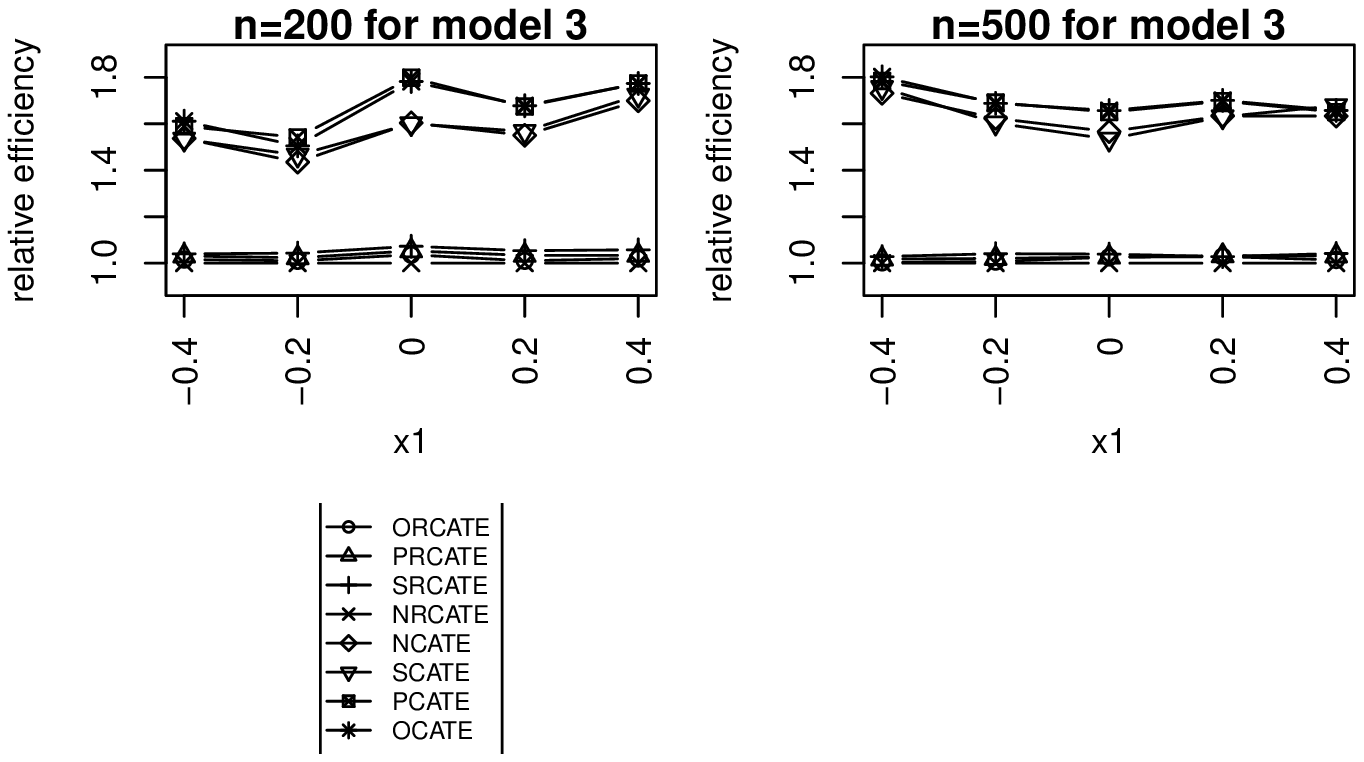}\\
	\includegraphics[width=10cm]{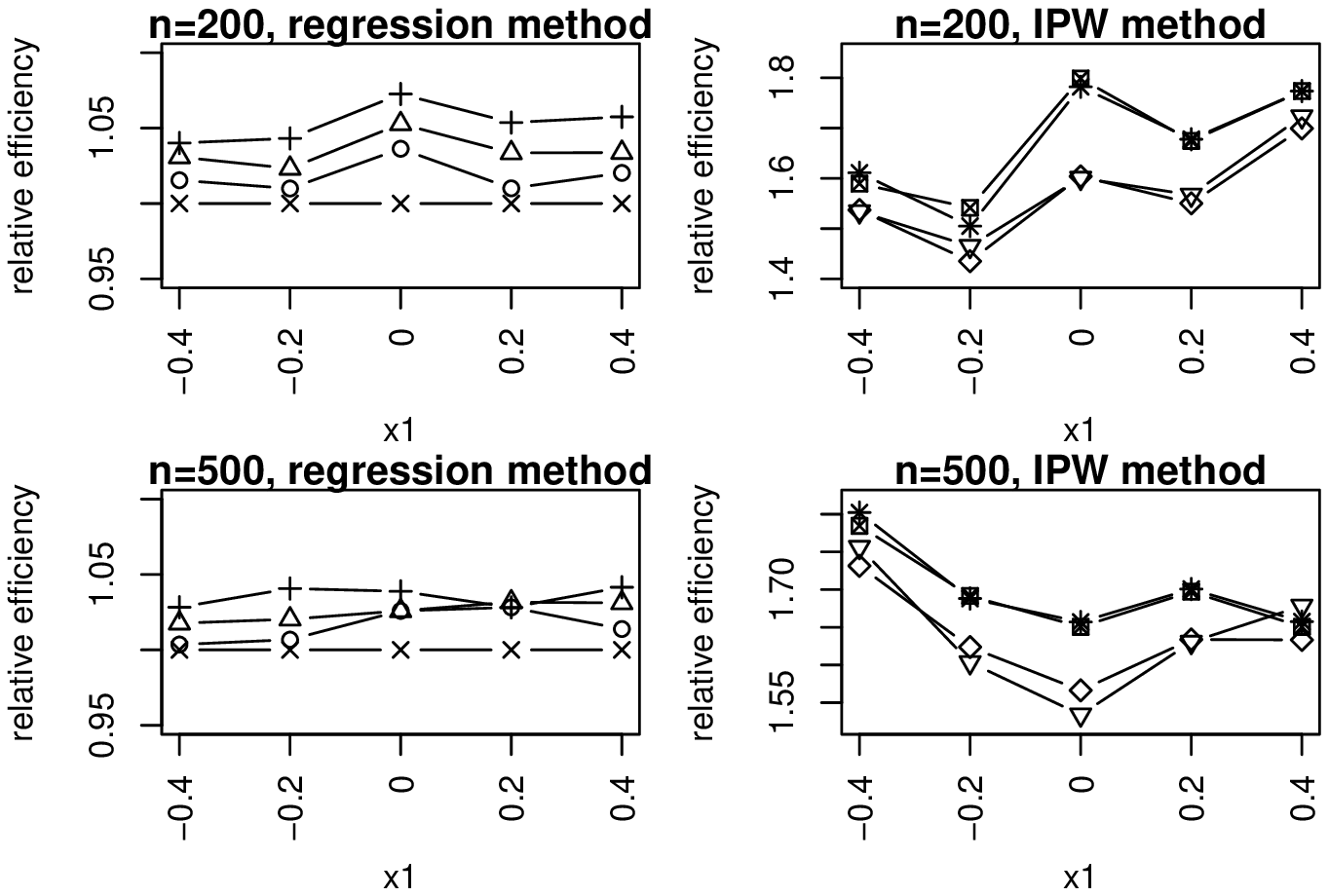}
	\caption{Relative efficiency of the CATE estimators against NRCATE for Model 3, which are based the results in panel2 of Table~\ref{tab3}.}\label{fig4}
\end{figure}

\section{Empirical applications}\label{application}
In this section, we apply SRCATE, as the dimensionality ($p=15$) of $\X$ is high, to analyse the {ACTG 175} data set that can be obtained from the R package {{speff2trial}}. This data set was collected from  a randomized clinical trial that evaluated treatment effect when either one or two therapies were used for HIV-infected adults; see \cite{Hammer:1996,Song:2008} for more details. As discussed before, our goal is to explore the heterogeneity of this treatment effect across subpopulations. Take {\it age} as $\X_{1}$ to check  how the expected pesticide effect changes with {\it age}.

A very brief description about the data set is as follows. The outcome here is CD4 T cell count at baseline and the treatment indicator variable $\D$ is a binary variable. $\D=0$ means  receiving zidovudine only and $\D=1$ means receiving two therapies simultaneously. As documented by a number of authors, we  take $\Y=\log_{10}(\text{CD4})$ and delete some infinite value after logarithmic transformation, then the number of observations is $n=2136$. Further, to guarantee the unconfoundedness assumption,  $\X$ consists of the following 15 covariates: the pidnum (patient's ID number); age (age in years at baseline); wtkg (weight in kg at baseline); hemo (hemophilia); homo (homosexual activity); drugs (history of intravenous drug use); karnof (Karnofsky score); oprior (non-zidovudine antiretroviral therapy prior to initiation of study treatment); zprior (zidovudine use prior to treatment initiation); preanti (number of days of previously received antiretroviral therapy); race; gender; str2 (antiretroviral history); offtrt (indicator of off-treatment before 96pm5 weeks); days (number of days until the first occurrence of: (i) a decline in CD4 T cell count of at least 50 (ii) an event indicating progression to AIDS, or (iii) death). 

We  now estimate CATE in the interval between 20 and 57 to avoid the boundary effect when nonparametric estimation method is involved. This range is about from $0.025$ quantile to $0.975$ quantile of the data. To apply SRCATE, we use the sufficient dimension reduction developed by \citet{Xia:2002}, 
which is now known to be MAVE to estimate the projection matrices $\be_1$ and $\be_0$, and the associated dimensions.  The results are $r(1)=2$ and $r(0)=3$. From these, we then have $s_{4}=\max\{r(1),r(0)\}+1=4$ and  $h_{4}=\widehat{\sigma}_{r} n^{-1/7}$ and $h=\widehat{\sigma}_{1} n^{-1/31}$, where $\widehat{\sigma}_{r}=\sqrt{\var(\be_{0}\ttop\X)}$, $\widehat{\be}_{0}$ is the estimated projection and $\widehat{\sigma}_{1}=2\sqrt{\var(\X_{1})}$. Similar to the simulation studies,  Gaussian kernel is used.

Figure~\ref{fig1} shows,  as a function of {\it age}, the curve of estimated CATE. Note that the curve is much above zero. In other words, receiving two therapies simultaneously has a much better treatment effect than receiving only one (zidovudine).  \cite{Song:2008} also obtained this conclusion. But the investigation on the heterogeneity shows that  the treatment effect is influenced by {\it age}. As shown in Figure~\ref{fig1}, before the age of 30, receiving two therapies leads to the immunity rise. After that, the advantage of this treatment is gradually weakened. Thus, such a treatment seems more useful for  patients whose ages are around 30.
\begin{figure}[h!]
	\centering
	\includegraphics[width=8cm]{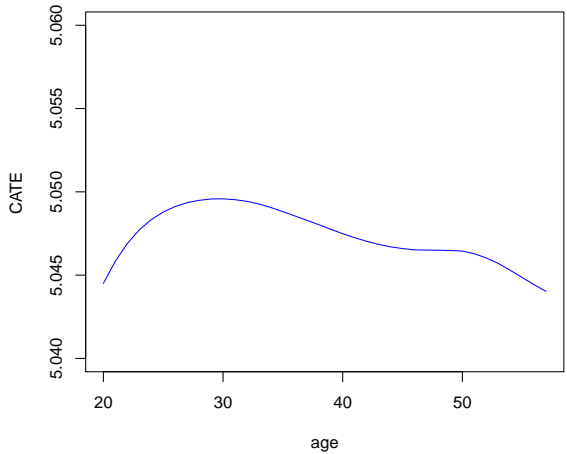}
	\caption{Conditional average treatment effect curves over age}\label{fig1}
\end{figure}

\section{Conclusion}\label{conclusion}
In this paper, we propose four regression-based estimators of CATE, aimed to capture the heterogeneity of a treatment effect across subpopulations. The systematic investigation shows the important factors that affect the asymptotic behaviours of the estimators: the convergence rates of the outcome regression functions and the affiliation of the given covariates to the set of arguments of the outcome regression functions. Further, any regression-based estimation can be asymptotically more efficient than any propensity score-based estimation, and  can at most achieve the asymptotic efficiency of nonparametric regression-based estimation in some cases. These results can give a relatively complete profile of propensity score-based and regression-based estimation for CATE. From the research, semiparametric regression-based estimation (SRCATE) is worth of recommendation as it can avoid model misspecification as well as the curse of dimensionality when some dimension reduction and feature selection approaches are combined. see \cite{Luo:2017} and \cite{Ma:2019}. In this paper, we only discuss the cases with correctly specified models. When the model is misspecified globally, further topics are about the asymptotic bias. Here global misspecification means that the assumed model is not convergent to the underlying model. If it is convergent, we call it local misspecification. Thus, we will check at which rate of convergence, the asymptotic bias vanishes and then also study its asymptotic efficiency. Another topic is about double robust estimation as it can greatly avoid model misspecification. The research is ongoing.

\section{Appendix}

\renewcommand{\theequation}{A.\arabic{equation}}
\setcounter{equation}{0}
Give some notations first.
\begin{description}
	\item (1) $C$ and $M$ stand for two generic bounded constants, $\Xi$ is the $\sigma$-field generated by $\X_{11},\ldots,\X_{1n}$.
	\item (2) $\epsilon_{ti}= Y_{i}-E(\Y_{(t)}|\X_{i})$, $\tau_{t}(\x_{1})=E[E\{\Y|\D=t,\X\}|\X_{1}=\x_{1}]$,
	$\Z^{t}=\be_{t}\ttop\X$ for $t=0,1$ and $i=1,\ldots,n$.
	\item (3) Write $\K_{1}\left(\frac{\X_{1i}-\X_{1}}{h_{1}}\right)$ as $\K_{1h}(\X_{1i})$; $\K_{2}\left(\frac{\X_{i}-\X_{j}}{h_{2}}\right)$ as $\K_{2h}(\X_{i}-\X_{j})$, and $\K_{4h}(\Z_{i}-\Z_{j})$ as $\K_{4}\left(\frac{\Z_{i}-\Z_{j}}{h_{4}}\right)$.
\end{description}
In the two-step estimation procedure for CATE, the second step involves,  
for $i=1,\ldots, n$, the quantities: 
{\small\begin{align*}
	\widehat{\K}_{1h}(\X_{1i})=\sum_{j:j\neq i}w_{ij}\K_{1h}(\X_{1j}).
	\end{align*}}We call it the estimator of $\K_{1h}(\X_{1j})$. In different circumstances, 
$w_{ij}$ can be  different. 
Take NRCATE as an example, and write $w_{ij}$ as $w_{ij}^{N}$: 
{\small\begin{align*}
	w_{ij}^{N}=\frac{\frac{1}{nh_{2}^{p}}\K_{2h}(\X_{i}-\X_{j})}{\frac{1}{nh_{2}^{p}}\sum\limits_{i=1}^{n}\K_{2h}(\X_{i}-\X_{j})\1(\D_{i}=1)}
	\end{align*}} that depends on $\X_{1},\ldots,\X_{n}$ only. 
{\lemm\label{Estimator of K1}Given assumptions (C1) - (C4) in Subsection~2.1 and (A1) - (A4) in Subsections~2.2-2.3,  
	\begin{align}\label{w1}
	|w_{ij}^{N}-w_{ji}^{N}|=\frac{\OO_{p}(h_{2})}{nh_{2}^{p}}\left|\K_{2h}\left(\X_{i}-\X_{j}\right)\right|,
	\end{align}}

\noindent {\sc Proof of Lemma \ref{Estimator of K1}.} By assumption (A2), $w_{ij}^{N}=w_{ji}^{N}=0$ for $||\X_{j}-\X_{i}||_{\infty}>h_{2}$ \citep{Abrevaya:2015}. Suppose that $||\X_{j}-\X_{i}||_{\infty}\leq h_{2}$. For all $j$, we define
{\small\begin{align*}
	\widehat{f}(\X_{j})=\frac{1}{nh_{2}^{p}}\sum\limits_{i:i\neq j}^{n}\K_{2h}\left(\X_{i}-\X_{j}\right).
	\end{align*}} It is clear that
{\small\begin{equation}\label{expansion of w}
	\begin{aligned}
	w_{ij}^{N}&=\frac{\frac{1}{nh_{2}^{p}}\K_{2h}\left(\X_{i}-\X_{j}\right)}{\frac{1}{nh_{2}^{p}}\sum\limits_{i=1}^{n}\K_{2h}\left(\X_{i}-\X_{j}\right)\1(\D_{i}=1)}=\frac{\frac{1}{nh_{2}^{p}}\K_{2h}\left(\X_{i}-\X_{j}\right)}{\frac{1}{nh_{2}^{p}}\sum\limits_{i=1}^{n}\K_{2h}\left(\X_{i}-\X_{j}\right)}\times \frac{\frac{1}{nh_{2}^{p}}\sum\limits_{i=1}^{n}\K_{2h}\left(\X_{i}-\X_{j}\right)}{\frac{1}{nh_{2}^{p}}\sum\limits_{i=1}^{n}\K_{2h}
		\left(\X_{i}-\X_{j}\right)\1(\D_{i}=1)}\\&=\frac{\frac{1}{nh_{2}^{p}}\K_{2h}\left(\X_{i}-\X_{j}\right)}{\widehat{f}(\X_{j})\widehat{p}(\X_{j})}.
	\end{aligned}\end{equation}}Then we have
{\small\begin{equation}\label{wij}\begin{aligned}
	&\quad|w_{ij}^{N}-w_{ji}^{N}|\\
	&=\frac{1}{nh_{2}^{p}}\left|\frac{\K_{2h}\left(\X_{i}-\X_{j}\right)}{\widehat{p}(\X_{j})\widehat{f}(\X_{j})}-\frac{\K_{2h}\left(\X_{j}-\X_{i}\right)}{\widehat{p}(\X_{i})\widehat{f}(\X_{i})}\right|=\frac{1}{nh_{2}^{p}}\left|\K_{2h}\left(\X_{i}-\X_{j}\right)\right|\left|\frac{1}{\widehat{p}(\X_{j})\widehat{f}(\X_{j})}-\frac{1}{\widehat{p}(\X_{i})\widehat{f}(\X_{i})}\right|\\
	&\leq\frac{1}{nh_{2}^{p}}\left|\K_{2h}\left(\X_{i}-\X_{j}\right)\right|\Bigg\{\left|\frac{1}{\widehat{p}(\X_{j})\widehat{f}(\X_{j})}-\frac{1}{p(\X_{j})f(\X_{j})}\right|+\left|\frac{1}{\widehat{p}(\X_{i})\widehat{f}(\X_{i})}-\frac{1}{p(\X_{i})f(\X_{i})}\right|\\
	&\quad+\left|\frac{1}{p(\X_{j})f(\X_{j})}-\frac{1}{p(\X_{i})f(\X_{i})}\right|\Bigg\}\\
	&=\frac{1}{nh_{2}^{p}}\left|\K_{2h}\left(\X_{i}-\X_{j}\right)\right|\Bigg\{\left|\frac{\widehat{p}(\X_{j})\widehat{f}(\X_{j})-p(\X_{j})f(\X_{j})}{\widehat{p}(\X_{j})p(\X_{j})\widehat{f}(\X_{j})f(\X_{j})}\right|+\left|\frac{\widehat{p}(\X_{i})\widehat{f}(\X_{i})-p(\X_{i})f(\X_{i})}{\widehat{p}(\X_{i})p(\X_{i})\widehat{f}(\X_{i})f(\X_{i})}\right|\\
	&\quad+\left|\frac{p(\X_{i})f(\X_{i})-p(\X_{j})f(\X_{j})}{p(\X_{i})p(\X_{j})f(\X_{i})f(\X_{j})}\right|\Bigg\}.
	\end{aligned}\end{equation}}Under conditions (C1)-(C4) and (A1)-(A4) for nonparametric estimation,
{\small\begin{align*}
	\underset{i}\sup|\widehat{f}(\X_{i})-f(\X_{i})|&=\OO_{p}\left(h_{2}^{s_{2}}+\sqrt{\frac{\log n}{nh_{2}^{p}}}\right),\\
	\underset{i}\sup|\widehat{p}(\X_{i})-p(\X_{i})|&=\OO_{p}\left(h_{2}^{s_{2}}+\sqrt{\frac{\log n}{nh_{2}^{p}}}\right).
	\end{align*}}Since $s_{2}\geq p\geq 2$, assumption (A3) implies that $\underset{i}\sup|\widehat{f}(\X_{i})-f(\X_{i})|=\oo_{p}(h_{2})$ and $\underset{i}\sup|\widehat{p}(\X_{i})-p(\X_{i})|=\oo_{p}(h_{2})$. By the mean value theorem,
{\small\begin{align*}
	\underset{j}\sup\left|\frac{\widehat{p}(\X_{j})\widehat{f}(\X_{j})-p(\X_{j})f(\X_{j})}{\widehat{p}(\X_{j})p(\X_{j})\widehat{f}(\X_{j})f(\X_{j})}\right|\leq \underset{j}\sup\frac{1}{\widetilde{p}^{2}(\X_{j})\widetilde{f}^{2}(\X_{j})}\underset{j}\sup\left|\widehat{p}(\X_{j})\widehat{f}(\X_{j})-p(\X_{j})f(\X_{j})\right|,
	\end{align*}}where $\widetilde{p}(\X_{j})$ is a quantity between $\widehat{p}(\X_{j})$ and $p(\X_{j})$, similarly, $\widetilde{f}(\X_{j})$ is also a quantity between $\widehat{f}(\X_{j})$ and $f(\X_{j})$. Owing to that $f$ and $p$ are bounded away from zero, $\underset{j}\sup\widetilde{p}^{-2}_{j}f^{-2}_{j}=\OO_{p}(1)$. After a simple calculation, we have
{\small\begin{align*}
	\underset{j}\sup\left|\widehat{p}(\X_{j})\widehat{f}(\X_{j})-p(\X_{j})f(\X_{j})\right|
	=\OO_{p}\left(h_{2}^{s_{2}}+\sqrt{\frac{\log n}{nh_{2}^{p}}}\right)=\oo_{p}(h_{2}).
	\end{align*}}Therefore,
\begin{align*}
\underset{j}\sup\left|\frac{\widehat{p}(\X_{j})\widehat{f}(\X_{j})-p(\X_{j})f(\X_{j})}{\widehat{p}(\X_{j})p(\X_{j})\widehat{f}(\X_{j})f(\X_{j})}\right|=\oo_{p}(h_{2}),\quad \underset{i}\sup\left|\frac{\widehat{p}(\X_{i})\widehat{f}(\X_{i})-p(\X_{i})f(\X_{i})}{\widehat{p}(\X_{i})p(\X_{i})\widehat{f}(\X_{i})f(\X_{i})}\right|=\oo_{p}(h_{2}).
\end{align*}As for the last term in \eqref{wij}, noticing that $f$ and $p$ are continuously differentiable on its compact support and bounded away from zero, we have  $\left|\frac{1}{f(\x_{1})p(\x_{1})}-\frac{1}{f(\x_{2})p(\x_{2})}\right|\leq M||\x_{1}-\x_{2}||_{\infty}$ for all $\x_{1},\x_{2}\in\CX$ and a constant $M>0$. $||\X_{j}-\X_{i}||_{\infty}\leq h_{2}$ leads to $\left|\frac{1}{f(\X_{i})p(\X_{i})}-\frac{1}{f(\X_{j})p(\X_{j})}\right|=\OO(h_{2})$. Combining all results yields \eqref{w1}.
\eop

\noindent {\sc Proof of Theorem \ref{Parametric}.} We can rewrite $\widehat{\m}_{1}(\X_{i})-\widehat{\m}_{0}(\X_{i})-\tau (x_1)$ as $\{\widehat{\m}_{1}(\X_{i})-\tau_{1}(\x_{1})\}-\{\widehat{\m}_{0}(\X_{i})-\tau_{0}(\x_{1})\}$. Then based on \eqref{parametric estimation},
{\small\begin{equation}\label{P1}\begin{aligned}
	&\sqrt{nh_{1}^{k}}(\widehat{\tau}(\x_{1})-\tau(\x_{1}))\\
	=&\frac{\frac{1}{\sqrt{nh_{1}^{k}}}\sum\limits_{i=1}^{n}\K_{1h}\left(\X_{1i}\right)\{[\widehat{\m}_{1}(\X_{i})-\tau_{1}(\x_{1})]-[\widehat{\m}_{0}(\X_{i})-\tau_{0}(\x_{1})]\}}{\frac{1}{nh_{1}^{k}}\sum\limits_{i=1}^{n}\K_{1h}(\X_{1i})}\\
	=&\frac{\frac{1}{\sqrt{nh_{1}^{k}}}\sum\limits_{i=1}^{n}\K_{1h}\left(\X_{1i}\right)\{[\widehat{\m}_{1}(\X_{i})-\tau_{1}(\x_{1})]-[\widehat{\m}_{0}(\X_{i})-\tau_{0}(\x_{1})]\}}{f(x_1)}(1+\oo_p(1)),
	\end{aligned}\end{equation}} as
{\small\begin{align*}
	\sup_{x_1}|\frac{1}{nh_{1}^{k}}\sum\limits_{i=1}^{n}\K_{1h}\left(\X_{1i}\right)-f(\x_{1})|=\oo_p(1).
	\end{align*}} First, deal with $\{\widehat{\m}_{1}(\X_{i})-\tau_{1}(\x_{1})\}$ in  \eqref{P1}. It is clear that
{\small\begin{equation}\label{I1}\begin{aligned}
	&\quad\frac{1}{\sqrt{nh_{1}^{k}}}\sum\limits_{i=1}^{n}\K_{1h}\left(\X_{1i}\right)[\widehat{\m}_{1}(\X_{i})-\tau_{1}(\x_{1})]\\
	&=\frac{1}{\sqrt{nh_{1}^{k}}}\left\{\sum\limits_{i=1}^{n}\K_{1h}\left(\X_{1i}\right)[\widehat{\m}_{1}(\X_{i})-\m_{1}(\X_{i})]+\sum\limits_{i=1}^{n}\K_{1h}\left(\X_{1i}\right)[\m_{1}(\X_{i})-\tau_{1}(\x_{1})]\right\}\\
	&=:\frac{1}{\sqrt{nh_{1}^{k}}} (\I_{n,1}+\I_{n,2}).
	\end{aligned}\end{equation}}
A simple calculation yields that
{\small\begin{align*}
	|\frac{1}{\sqrt{nh_{1}^{k}}}\I_{n,1}|\leq \underset{\x}\sup\left|\widehat{\m}_{1}(\X_{i})-\m_{1}(\X_{i})\right|\frac{1}{nh_{1}^{k}}\sum\limits_{i=1}^{n}\left|\K_{1h}\left(\X_{1i}\right)\right|.
	\end{align*}} As $h_{1}\rightarrow 0$ and $\frac{1}{nh_{1}^{k}}\sum\limits_{i=1}^{n}|\K_{1h}(\X_{1i})|=\OO_{p}(1)$, we then have $\frac{1}{\sqrt{nh_{1}^{k}}}\I_{n,1}=\OO_{p}(\sqrt{h_{1}^{k}})=\oo_{p}(1)$. Thus, equation \eqref{I1} becomes
{\small\begin{align*}
	\frac{1}{\sqrt{nh_{1}^{k}}}\sum\limits_{i=1}^{n}\K_{1h}\left(\X_{1i}\right)[\widehat{\m}_{1}(\X_{i})-\tau_{1}(\x_{1})]=\frac{1}{\sqrt{nh_{1}^{k}}}\sum\limits_{i=1}^{n}\K_{1h}\left(\X_{1i}\right)[\m_{1}(\X_{i})-\tau_{1}(\x_{1})]+\oo_{p}(1).
	\end{align*}}Similarly,
{\small\begin{align*}
	\frac{1}{\sqrt{nh_{1}^{k}}}\sum\limits_{i=1}^{n}\K_{1h}\left(\X_{1i}\right)[\widehat{\m}_{0}(\X_{i})-\tau_{0}(\x_{1})]=\frac{1}{\sqrt{nh_{1}^{k}}}\sum\limits_{i=1}^{n}\K_{1h}\left(\X_{1i}\right)[\m_{0}(\X_{i})-\tau_{0}(\x_{1})]+\oo_{p}(1).
	\end{align*}} Altogether, the asymptotically linear representation of $\widehat{\tau}(\x_{1})$ is
{\small\begin{align*}
	\sqrt{nh_{1}^{k}}\{\widehat{\tau}(\x_{1})-\tau(\x_{1})\}&=\frac{\frac{1}{\sqrt{nh_{1}^{k}}}\sum\limits_{i=1}^{n}\K_{1h}\left(\X_{1i}\right)\{\m_{1}(\X_{i})-\m_{0}(\X_{i})-\tau(\x_{1})\}}{f(x_1)}(1+\oo_{p}(1))\\
	&=\frac{\frac{1}{\sqrt{nh_{1}^{k}}}\sum\limits_{i=1}^{n}\K_{1h}\left(\X_{1i}\right)\{\m_{1}(\X_{i})-\m_{0}(\X_{i})-\tau(\x_{1})\}}{f(x_1)}+\oo_{p}(1).
	\end{align*}}The second equation is due to the asymptotic finiteness of the leading term that is asymptotically normal shown below. As it is the sum of independent variables, the asymptotic normality is easy to derive. Specifically, noticing that the random variables {\small\begin{align*}
	\left\{\K_{1h}\left(\X_{1i}\right)[\m_{1}(\X_{i})-\m_{0}(\X_{i})-\tau(\X_{1i})]\right\}_{i=1}^{n}
	\end{align*}}are i.i.d., then we can apply Lyapunov's central limit theorem to obtain the asymptotic distribution shown in Theorem \ref{Parametric}. 
Under the assumptions (C1)- (C4) and (A1), we derive that
\begin{align*}
\sqrt{nh_{1}^{k}}\left\{\widehat{\tau}(\x_{1})-\tau(\x_{1})\right\}\xrightarrow{d}\N\left(0,\frac{||\K_{1}||_{2}^{2}\sigma_{P}^{2}(\x_{1})}{f(\x_{1})}\right),
\end{align*}
we now give the formula of  $\sigma_{P}^{2}(\x_{1})$. It is easy to see that when $n\to \infty$, the variance of
{\small\begin{align*}
	\quad\frac{\frac{1}{\sqrt{nh_{1}^{k}}}\sum\limits_{i=1}^{n}\K_{1h}\left(\X_{1i}\right)\{\m_{1}(\X_{i})-\m_{0}(\X_{i})-\tau(\x_{1})\}}{f(x_1)}
	\end{align*}}
converges to
\begin{align*}
\sigma_{P}^{2}(\x_{1}):= E[\{\m_{1}(\X)-\m_{0}(\X)-\tau(\x_{1})\}^{2}|\X_{1}=\x_{1}].
\end{align*}The proof of Theorem \ref{Parametric} is finished.
\eop

\noindent {\sc Proof of Theorem \ref{NonParametric 1}.} 
First, we have
{\small\begin{equation}\label{nonparametric estimator}
	\begin{aligned}
	&\quad\sqrt{nh_{1}^{k}}(\widehat{\tau}(\x_{1})-\tau(\x_{1}))\\
	&=\frac{\frac{1}{\sqrt{nh_{1}^{k}}}\sum\limits_{i=1}^{n}\K_{1h}\left(\X_{1i}\right)\left[\widehat{\m}_{1}(\X_{i})-\tau_{1}(\x_{1})\right]}{\frac{1}{nh_{1}^{k}}\sum\limits_{i=1}^{n}\K_{1h}(\X_{1i})}-\frac{\frac{1}{\sqrt{nh_{1}^{k}}}\sum\limits_{i=1}^{n}\K_{1h}\left(\X_{1i}\right)\left[\widehat{\m}_{0}(\X_{i})-\tau_{0}(\x_{1})\right]}{\frac{1}{nh_{1}^{k}}\sum\limits_{i=1}^{n}\K_{1h}(\X_{1i})},
	\end{aligned}\end{equation}}where \begin{align*}
\widehat{\m}_{1}(\X_{i})=\frac{\frac{1}{nh_{2}^{p}}\sum\limits_{j=1}^{n}\K_{2h}\left(\X_{j}-\X_{i}\right)\Y_{1j}\1(\D_{j}=1)}{\frac{1}{nh_{2}^{p}}\sum\limits_{j=1}^{n}\K_{2h}\left(\X_{j}-\X_{i}\right)\1(\D_{j}=1)},
\end{align*}\begin{align*}
\widehat{\m}_{0}(\X_{i})=\frac{\frac{1}{nh_{2}^{p}}\sum\limits_{j=1}^{n}\K_{2h}\left(\X_{j}-\X_{i}\right)\Y_{0j}\1(\D_{j}=0)}{\frac{1}{nh_{2}^{p}}\sum\limits_{j=1}^{n}\K_{2h}\left(\X_{j}-\X_{i}\right)\1(\D_{j}=0)}.
\end{align*}
Similarly as the proof for Theorem~\ref{Parametric}, we have the following decomposition:
{\small\begin{equation}\label{para}\begin{aligned}
	&\quad\frac{1}{\sqrt{nh_{1}^{k}}}\sum\limits_{i=1}^{n}\K_{1h}\left(\X_{1i}\right)\left[\m_{1}(\X_{i})-\tau_{1}(\x_{1})\right]\\
	&=\frac{1}{\sqrt{nh_{1}^{k}}}\sum\limits_{i=1}^{n}\epsilon_{1i}\1(\D_{i}=1)\frac{\K_{1h}\left(\X_{1i}\right)}{p(\X_{i})}+\frac{1}{\sqrt{nh_{1}^{k}}}\sum\limits_{i=1}^{n}\K_{1h}\left(\X_{1i}\right)[E(\Y_{(1)}|\X_{i})-\tau_{1}(\x_{1})]\\
	&\quad+\frac{1}{\sqrt{nh_{1}^{k}}}\sum\limits_{i=1}^{n}\epsilon_{1i}\1(\D_{i}=1)\sum\limits_{j=1}^{n}\K_{1h}\left(\X_{1j}\right)(w_{ij}^{N}-w_{ji}^{N})\\
	&\quad+\frac{1}{\sqrt{nh_{1}^{k}}}\sum\limits_{i=1}^{n}\epsilon_{1i}\1(\D_{i}=1)\left[\sum\limits_{j=1}^{n}\K_{1h}\left(\X_{1j}\right)w_{ji}^{N}-\frac{\K_{1h}\left(\X_{1i}\right)}{p(\X_{i})}\right]\\
	&\quad+\frac{1}{\sqrt{nh_{1}^{k}}}\sum\limits_{i=1}^{n}\K_{1h}\left(\X_{1i}\right)\left[\frac{\frac{1}{nh_{2}^{p}}\sum\limits_{j=1}^{n}\K_{2h}\left(\X_{j}-\X_{i}\right)\1(\D_{j}=1)E\Y_{(1)}|\X_{j}}{\frac{1}{nh_{2}^{p}}\sum\limits_{j=1}^{n}\K_{2h}\left(\X_{j}-\X_{i}\right)\1(\D_{j}=1)}-E\Y_{(1)}|\X_{i}\right]\\
	&=:\I_{n,3}+\I_{n,4}+\I_{n,5}+\I_{n,6}+\I_{n,7},
	\end{aligned}\end{equation}}where
\begin{align*}
w_{ij}^{N}=\frac{\frac{1}{nh_{2}^{p}}\K_{2h}\left(\X_{i}-\X_{j}\right)}{\frac{1}{nh_{2}^{p}}\sum\limits_{i=1}^{n}\K_{2h}\left(\X_{i}-\X_{j}\right)\1(\D_{i}=1)}, \quad \epsilon_{1i}=\Y_{i}-E(\Y_{(1)}|\X_{i}).
\end{align*}Note that $\I_{n,3}$ and $\I_{n,4}$ in equation \eqref{para} yield the final expression in Theorem \ref{NonParametric 1}. Therefore, we need to show that $\I_{n,5}$, $\I_{n,6}$ and $\I_{n,7}$ in equation \eqref{para} are all $\oo_{p}(1)$.

First  show that $\I_{n,5}=\oo_{p}(1)$. From Lemma~\ref{Estimator of K1},
{\small\begin{align*}
	&\quad\frac{1}{\sqrt{h_{1}^{k}}}\underset{i}\sup\left|\sum\limits_{j=1}^{n}\K_{1h}\left(\X_{1j}\right)(w_{ij}^{N}-w_{ji}^{N})\right|
	\leq \frac{1}{\sqrt{h_{1}^{k}}}\underset{i}\sup\sum\limits_{j:j\neq i}(w_{ij}^{N}-w_{ji}^{N})\left|\K_{1h}\left(\X_{1j}\right)\right|\\
	&\leq\frac{MC}{h_{2}}\times\frac{h_{2}}{\sqrt{h_{1}^{k}}}\times\underset{i}\sup\sum\limits_{j:j\neq i}\frac{1}{nh_{2}^{p}}\left|\K_{2h}\left(\X_{i}-\X_{j}\right)\right|=\OO_{p}(1)\times\oo_{p}(1)\times\OO_{p}(1)=\oo_{p}(1),
	\end{align*}}
Further, $\frac{1}{\sqrt{n}}\sum\limits_{i=1}^{n}\epsilon_{1i}\1(\D_{i}=1)$ has finite limit and thus, is bounded by $\OO_p(1)$ and then $\I_{n,5}=\oo_{p}(1)$.

Deal with $\I_{n,6}$. As
{\small\begin{align*}
	\sum\limits_{j=1}^{n}\K_{1h}\left(\X_{1j}\right)w_{ji}^{N}=\frac{\frac{1}{nh_{2}^{p}}\sum\limits_{j=1}^{n}\K_{2h}\left(\X_{j}-\X_{i}\right)}{\frac{1}{nh_{2}^{p}}\sum\limits_{j=1}^{n}\K_{2h}\left(\X_{j}-\X_{i}\right)\1(\D_{j}=1)}\frac{\frac{1}{nh_{2}^{p}}\sum\limits_{j=1}^{n}\K_{2h}\left(\X_{j}-\X_{i}\right)\K_{1h}\left(\X_{1j}\right)}{\frac{1}{nh_{2}^{p}}\sum\limits_{j=1}^{n}\K_{2h}\left(\X_{j}-\X_{i}\right)},
	\end{align*}} we can then regard $\sum\limits_{j=1}^{n}\K_{1h}\left(\X_{1j}\right)w_{ji}^{N}$ as an estimator of $\frac{\K_{1h}\left(\X_{1i}\right)}{p(\X_{i})}$. Consider {\small\begin{align*}
	\left[\sum\limits_{j=1}^{n}\K_{1h}\left(\X_{1j}\right)w_{ji}^{N}-\frac{\K_{1h}\left(\X_{1i}\right)}{p(\X_{i})}\right],
	\end{align*}} which is the bias of $\frac{\K_{1h}\left(\X_{1i}\right)}{p(\X_{i})}$ to $\frac{\K_{1h}\left(\X_{1i}\right)}{p(\X_{i})}$. Write $\X=(\X_{1},\X_{(2)})$ and {\small\begin{align*}
	\K_{2h}\left(\X-\X_{j}\right)=\K_{21}\left(\frac{\X_{1}-\X_{1j}}{h_{2}}\right)\K_{22}\left(\frac{\X_{(2)}-\X_{2j}}{h_{2}}\right).
	\end{align*}}Since $\widehat{f}-f=\oo_{p}(1)$, and the kernel function is $s^{*}$ ($\geq s_{2}$) times continuously differentiable, we have
{\small\begin{equation}\label{conditional expectaion}\begin{aligned}
	&\quad E\left\{\sum\limits_{j=1}^{n}\K_{1h}\left(\X_{1j}\right)w_{ji}^{N}\Bigg|\X_{i}\right\}\\
	&=\frac{1+\oo_{p}(1)}{h_{2}^{p}f(\X_{i})p(\X_{i})}\int \K_{21}\left(\frac{u_{1j}-\X_{1i}}{h_{2}}\right)\K_{22}\left(\frac{u_{2j}-\X_{2i}}{h_{2}}\right)\K_{1h}\left(u_{1j}\right)f(u_{i})du\\
	&=\frac{1+\oo_{p}(1)}{f(\X_{i})p(\X_{i})}\int \K_{21}(v_{1})\K_{22}(v_{2})\K_{1}\left(\frac{\X_{1i}-\X_{1}}{h_{1}}+v_{1}\frac{h_{2}}{h_{1}}\right)f(\X_{i}+h_{2}v)dv\\
	&=\frac{\K_{1h}\left(\X_{1i}\right)}{p(\X_{i})}+\OO_{p}\left(\frac{h_{2}^{s_{2}}}{h_{1}^{s_{2}}}\right).
	\end{aligned}\end{equation}}Note that
{\small\begin{align*}
	&\quad\frac{\widehat{\K}_{1h}\left(\X_{1i}\right)}{\widehat{p}(\X_{i})}-\frac{\K_{1h}\left(\X_{1i}\right)}{p(\X_{i})}\\
	&=\left\{\frac{1}{\widehat{p}(\X_{i})}-\frac{1}{p(\X_{i})}+\frac{1}{p(\X_{i})}\right\}\left\{\widehat{\K}_{1h}\left(\X_{1i}\right)-\K_{1h}\left(\X_{1i}\right)+\K_{1h}\left(\X_{1i}\right)\right\}\\
	&\quad-\frac{\K_{1h}\left(\X_{1i}\right)}{p(\X_{i})}=\left\{\frac{1}{\widehat{p}(\X_{i})}-\frac{1}{p(\X_{i})}\right\}\left\{\widehat{\K}_{1h}\left(\X_{1i}\right)-\K_{1h}\left(\X_{1i}\right)\right\}\\
	&\quad+\frac{1}{p(\X_{i})}\left\{\widehat{\K}_{1h}\left(\X_{1i}\right)-\K_{1h}\left(\X_{1i}\right)\right\}+\left\{\frac{1}{\widehat{p}(\X_{i})}-\frac{1}{p(\X_{i})}\right\}\K_{1h}\left(\X_{1i}\right)\\
	&=\OO_{p}\left(\frac{h_{2}^{s_{2}}}{h_{1}^{s_{2}}}+h_{2}^{s_{2}}+\sqrt{\frac{\log n}{nh_{2}^{p}}}\right)=\OO_{p}\left(\frac{h_{2}^{s_{2}}}{h_{1}^{s_{2}}}\right).
	\end{align*}}Thus, $\underset{i}\sup\left|\sum\limits_{j=1}^{n}\K_{1h}\left(\X_{1j}\right)w_{ji}^{N}-\frac{\K_{1h}\left(\X_{1i}\right)}{p(\X_{i})}\right|=\OO_{p}\left(\frac{h_{2}^{s_{2}}}{h_{1}^{s_{2}}}\right)$. Owing to assumption (A4) that $\frac{h_{2}^{2s_{2}}}{h_{1}^{2s_{2}+k}}\rightarrow 0$, we have
{\small\begin{align*}
	\underset{i}\sup\left|\frac{1}{\sqrt{h_{1}^{k}}}\left[\sum\limits_{j=1}^{n}\K_{1h}\left(\X_{1j}\right)w_{ji}^{N}-\frac{\K_{1h}\left(\X_{1i}\right)}{p(\X_{i})}\right]\right|=\OO_{p}\left(\frac{h_{2}^{2s_{2}}}{h_{1}^{2s_{2}+k}}\right)=\oo_{p}(1).
	\end{align*}}Since $\epsilon_{1i}=\Y_{i}-E(\Y_{(1)}|\X_{i})$ are mutually independent,  we have $\I_{n,6}=\oo_{p}(1)$ in equation \eqref{para}.
Finally, to show that $\I_{n,7}=\oo_{p}(1)$ of equation \eqref{para}. Note that
\begin{align*}
&\quad\frac{\frac{1}{nh_{2}^{p}}\sum\limits_{j=1}^{n}\K_{2h}\left(\X_{j}-\X_{i}\right)\1(\D_{j}=1)E\Y_{(1)}|\X_{j}}{\frac{1}{nh_{2}^{p}}\sum\limits_{j=1}^{n}\K_{2h}\left(\X_{j}-\X_{i}\right)\1(\D_{j}=1)}\\
&=\frac{\frac{1}{nh_{2}^{p}}\sum\limits_{j=1}^{n}\K_{2h}\left(\X_{j}-\X_{i}\right)\1(\D_{j}=1)E\Y_{(1)}|\X_{j}}{\frac{1}{nh_{2}^{p}}\sum\limits_{j=1}^{n}\K_{2h}\left(\X_{j}-\X_{i}\right)}\cdot\frac{\frac{1}{nh_{2}^{p}}\sum\limits_{j=1}^{n}\K_{2h}\left(\X_{j}-\X_{i}\right)}{\frac{1}{nh_{2}^{p}}\sum\limits_{j=1}^{n}\K_{2h}\left(\X_{j}-\X_{i}\right)\1(\D_{j}=1)},
\end{align*} which can be viewed as an estimator of $\frac{E\{\1(D=1)\Y_{(1)}|\X_{i}\}}{p(\X_{i})}$. Denote $\A(\X_{i})=E\{\1(D=1)\Y_{(1)}|\X_{i}\}$. We can derive easily that
{\small\begin{align*}
	&\quad\frac{\widehat{\A}(\X_{i})}{\widehat{p}(\X_{i})}-\frac{\A(\X_{i})}{p(\X_{i})}\\
	&=\left\{\widehat{\A}(\X_{i})-\A(\X_{i})+\A(\X_{i})\right\}\left\{\frac{1}{\widehat{p}(\X_{i})}-\frac{1}{p(\X_{i})}+\frac{1}{p(\X_{i})}\right\}-\frac{\A(\X_{i})}{p(\X_{i})}\\
	&=\left\{\widehat{\A}(\X_{i})-\A(\X_{i})\right\}\left\{\frac{1}{\widehat{p}(\X_{i})}-\frac{1}{p(\X_{i})}\right\}+\A(\X_{i})\left\{\frac{1}{\widehat{p}(\X_{i})}-\frac{1}{p(\X_{i})}\right\}\\
	&\quad+\left\{\widehat{\A}(\X_{i})-\A(\X_{i})\right\}\frac{1}{p(\X_{i})}=\OO_{p}\left(h_{2}^{s_{2}}+\sqrt{\frac{\log n}{nh_{2}^{p}}}\right).
	\end{align*}}Thus
\begin{align*}
\underset{i}\sup\left|\frac{\frac{1}{nh_{2}^{p}}\sum\limits_{j=1}^{n}\K_{2h}\left(\X_{j}-\X_{i}\right)\1(\D_{j}=1)E\Y_{(1)}|\X_{j}}{\frac{1}{nh_{2}^{p}}\sum\limits_{j=1}^{n}\K_{2h}\left(\X_{j}-\X_{i}\right)\1(\D_{j}=1)}-E\Y_{(1)}|\X_{i}\right|=\OO_{p}\left(h_{2}^{s_{2}}+\sqrt{\frac{\log n}{nh_{2}^{p}}}\right).
\end{align*}Then, we can bound $\I_{n,7}$ as follows:
{\small\begin{align*}
	|\I_{n,7}|&=\left|\frac{1}{\sqrt{nh_{1}^{k}}}\sum\limits_{i=1}^{n}\K_{1h}\left(\X_{1i}\right)\left[\frac{\frac{1}{nh_{2}^{p}}\sum\limits_{j=1}^{n}\K_{2h}\left(\X_{j}-\X_{i}\right)\1(\D_{j}=1)E\Y_{(1)}|\X_{j}}{\frac{1}{nh_{2}^{p}}\sum\limits_{j=1}^{n}\K_{2h}\left(\X_{j}-\X_{i}\right)\1(\D_{j}=1)}-E\Y_{(1)}|\X_{i}\right]\right|\\
	&\leq\sqrt{nh_{1}^{k}}\underset{i}\sup\left|\frac{\frac{1}{nh_{2}^{p}}\sum\limits_{j=1}^{n}\K_{2h}\left(\X_{j}-\X_{i}\right)\1(\D_{j}=1)E\Y_{(1)}|\X_{j}}{\frac{1}{nh_{2}^{p}}\sum\limits_{j=1}^{n}\K_{2h}\left(\X_{j}-\X_{i}\right)\1(\D_{j}=1)}-E\Y_{(1)}|\X_{i}\right|\frac{1}{nh_{1}^{k}}\sum\limits_{i=1}^{n}\left|\K_{1h}\left(\X_{1i}\right)\right|\\
	&=\sqrt{nh_{1}^{k}}\OO_{p}\left(h_{2}^{s_{2}}+\sqrt{\frac{\log n}{nh_{2}^{p}}}\right)\cdot\OO_{p}(1)=\oo_{p}(1)\cdot\OO_{p}(1)=\oo_{p}(1),
	\end{align*}}where  assumption (A4) is used for the second equation. Thus, together with $\I_{n,5}=\oo_{p}(1)$, $\I_{n,6}=\oo_{p}(1)$ and $\I_{n,7}=\oo_{p}(1)$, equation \eqref{para} becomes
{\small\begin{align*}
	&\quad\frac{1}{\sqrt{nh_{1}^{k}}}\sum\limits_{i=1}^{n}\K_{1h}\left(\X_{1i}\right)\left[\frac{\frac{1}{nh_{2}^{p}}\sum\limits_{j=1}^{n}\K_{2h}\left(\X_{j}-\X_{i}\right)\Y_{1j}\1(\D_{j}=1)}{\frac{1}{nh_{2}^{p}}\sum\limits_{j=1}^{n}\K_{2h}\left(\X_{j}-\X_{i}\right)\1(\D_{j}=1)}-\tau_{1}(\x_{1})\right]\\
	&=\I_{n,3}+\I_{n,4}+\oo_{p}(1).
	\end{align*}}
Similarly, we can also deal with $\widehat{\m}_{0}(\X_{i})-\tau_{0}(\x_{1})$ of \eqref{nonparametric estimator} to have
{\small\begin{align*}
	&\quad\frac{1}{\sqrt{nh_{1}^{k}}}\sum\limits_{i=1}^{n}\K_{1h}\left(\X_{1i}\right)\left[\frac{\frac{1}{nh_{2}^{p}}\sum\limits_{j=1}^{n}\K_{2h}\left(\X_{j}-\X_{i}\right)\Y_{0j}\1(\D_{j}=0)}{\frac{1}{nh_{2}^{p}}\sum\limits_{j=1}^{n}\K_{2h}\left(\X_{j}-\X_{i}\right)\1(\D_{j}=0)}-\tau_{1}(\x_{1})\right]\\
	&:=\I_{n,8}+\I_{n,9}+\oo_{p}(1),
	\end{align*}}where
{\small\begin{align*}
	\I_{n,8}&=\frac{1}{\sqrt{nh_{1}^{k}}}\sum\limits_{i=1}^{n}\epsilon_{0i}\1(\D_{i}=0)\frac{\K_{1h}\left(\X_{1i}\right)}{1-p(\X_{i})},\quad \I_{n,9}=\frac{1}{\sqrt{nh_{1}^{k}}}\sum\limits_{i=1}^{n}\K_{1h}\left(\X_{1i}\right)E\Y_{(0)}|\X_{i},\\
	\epsilon_{0i}&=\Y_{i}-E\Y_{(0)}|\X_{i}.
	\end{align*}}Hence, we get the asymptotic linear representation of $\widehat{\tau}(\x_{1})$  as
\begin{align*}
\sqrt{nh_{1}^{k}}\{\widehat{\tau}(\x_{1})-\tau(\x_{1})\}=\frac{1}{\sqrt{nh_{1}^{k}}}\frac{1}{f(\x_{1})}\sum\limits_{i=1}^{n}\{\Psi_{1}(\X_{i},\Y_{i},\D_{i})-\tau(\x_{1})\}\K_{1h}\left(\X_{1i}\right)+\oo_{p}(1),
\end{align*}
which can be asymptotically normal.
Again, we compute its asymptotic variance. Similarly as the proof for Theorem \ref{Parametric}, we have
\begin{align*}
\text{Var}\{\widehat{\tau}(\x_{1})\}=\frac{1}{nh_{1}^{k}}\frac{||\K_{1}||_{2}^{2}\sigma_{N}^{2}(\x_{1})}{f(\x_{1})}+\oo\left(\frac{1}{nh_{1}^{k}}\right).
\end{align*} Then by assumptions (C1)-- (C4) and (A1) -- (A4)  for some $s^{*}\geq s_{2}\geq p$, we can derive that
\begin{align*}
\sqrt{nh_{1}^{k}}\left\{\widehat{\tau}(\x_{1})-\tau(\x_{1})\right\}\xrightarrow{d}\N\left(0,\frac{||\K_{1}||_{2}^{2}\sigma_{N}^{2}(\x_{1})}{f(\x_{1})}\right),
\end{align*}where \begin{align*}
\sigma_{N}^{2}(\x_{1})\equiv E[\{\Psi_{1}(\X,\Y,\D)-\tau(\x_{1})\}^{2}|\X_{1}=\x_{1}].
\end{align*} The proof  is concluded.
\eop

\noindent {\sc Proof of Theorem \ref{SemiParametric 1}.}  Inspired by the proof of Theorem 2 of \cite{Luo:2017}, we have
{\small\begin{equation}\label{non}\begin{aligned}
	&\quad\sqrt{nh_{1}^{k}}(\widehat{\tau}(\x_{1})-\tau(\x_{1}))\\
	&=\frac{\frac{1}{\sqrt{nh_{1}^{k}}}\sum\limits_{i=1}^{n}\K_{1h}\left(\X_{1i}\right)\left[\widehat{\m}_{1}(\widehat{\be}_{1}\ttop\X)-\tau_{1}(\x_{1})\right]}{\frac{1}{nh_{1}^{k}}\sum\limits_{i=1}^{n}\K_{1h}(\X_{1i})}-\frac{\frac{1}{\sqrt{nh_{1}^{k}}}\sum\limits_{i=1}^{n}\K_{1h}\left(\X_{1i}\right)\left[\widehat{\m}_{0}(\widehat{\be}_{0}\ttop\X)-\tau_{0}(\x_{1})\right]}{\frac{1}{nh_{1}^{k}}\sum\limits_{i=1}^{n}\K_{1h}(\X_{1i})}\\
	&=\frac{\frac{1}{\sqrt{nh_{1}^{k}}}\sum\limits_{i=1}^{n}\K_{1h}\left(\X_{1i}\right)\left[\widehat{\m}_{1}(\be_{1}\ttop\X)-\tau_{1}(\x_{1})\right]}{\frac{1}{nh_{1}^{k}}\sum\limits_{i=1}^{n}\K_{1h}(\X_{1i})}-\frac{\frac{1}{\sqrt{nh_{1}^{k}}}\sum\limits_{i=1}^{n}\K_{1h}\left(\X_{1i}\right)\left[\widehat{\m}_{0}(\be_{0}\ttop\X)-\tau_{0}(\x_{1})\right]}{\frac{1}{nh_{1}^{k}}\sum\limits_{i=1}^{n}\K_{1h}(\X_{1i})}\\
	&\quad +\OO_{p}(\sqrt{nh_{1}^{k}}||\widehat{\be}_{1}-\be_{1}||+\sqrt{nh_{1}^{k}}||\widehat{\be}_{0}-\be_{0}||),
	\end{aligned}\end{equation}}where{\small\begin{align*}
	\widehat{\m}_{1}(\widehat{\be}_{1}\ttop\X)&=\frac{\frac{1}{nh_{4}^{r(1)}}\sum\limits_{j=1}^{n}\K_{4h}\left(\widehat{\Z}_{j}^{1}-\widehat{\Z}_{i}^{1}\right)\Y_{1j}\1(\D_{j}=1)}{\frac{1}{nh_{4}^{r(1)}}\sum\limits_{j=1}^{n}\K_{4h}\left(\widehat{\Z}_{j}^{1}-\widehat{\Z}_{i}^{1}\right)\1(\D_{j}=1)},\quad \widehat{\Z}^{1}=\widehat{\be}_{1}\ttop\X,\\
	\widehat{\m}_{0}(\widehat{\be}_{0}\ttop\X)&=\frac{\frac{1}{nh_{4}^{r(1)}}\sum\limits_{j=1}^{n}\K_{4h}\left(\widehat{\Z}_{j}^{0}-\widehat{\Z}_{i}^{0}\right)\Y_{0j}\1(\D_{j}=0)}{\frac{1}{nh_{4}^{r(1)}}\sum\limits_{j=1}^{n}\K_{4h}\left(\widehat{\Z}_{j}^{0}-\widehat{\Z}_{i}^{0}\right)\1(\D_{j}=0)},\quad \widehat{\Z}^{0}=\widehat{\be}_{0}\ttop\X.
	\end{align*}}Under  assumptions (A8), $\OO_{p}(\sqrt{nh_{1}^{k}}||\widehat{\be}_{1}-\be_{1}||+\sqrt{nh_{1}^{k}}||\widehat{\be}_{0}-\be_{0}||)=\OO_{p}(\sqrt{h_{1}^{k}})=\oo_{p}(1)$ as $h_{1}\rightarrow 0$.
Therefore, equation \eqref{non} becomes
{\small\begin{equation}\label{semiparametric estimator}
	\begin{aligned}
	&\quad\sqrt{nh_{1}^{k}}(\widehat{\tau}(\x_{1})-\tau(\x_{1}))\\
	&=\frac{\frac{1}{\sqrt{nh_{1}^{k}}}\sum\limits_{i=1}^{n}\K_{1h}\left(\X_{1i}\right)\left[\widehat{\m}_{1}(\be_{1}\ttop\X)-\tau_{1}(\x_{1})\right]}{\frac{1}{nh_{1}^{k}}\sum\limits_{i=1}^{n}\K_{1h}(\X_{1i})}-\frac{\frac{1}{\sqrt{nh_{1}^{k}}}\sum\limits_{i=1}^{n}\K_{1h}\left(\X_{1i}\right)\left[\widehat{\m}_{0}(\be_{0}\ttop\X)-\tau_{0}(\x_{1})\right]}{\frac{1}{nh_{1}^{k}}\sum\limits_{i=1}^{n}\K_{1h}(\X_{1i})}\\
	&+\oo_{p}(1).
	\end{aligned}\end{equation}}Similarly as the proof for Theorem \ref{NonParametric 1}, we have
{\small\begin{align*}
	&\quad\frac{1}{\sqrt{nh_{1}^{k}}}\sum\limits_{i=1}^{n}\K_{1h}\left(\X_{1i}\right)\left[\widehat{\m}_{1}(\be_{1}\ttop\X)-\tau_{1}(\x_{1})\right]\\
	&=\frac{1}{\sqrt{nh_{1}^{k}}}\sum\limits_{i=1}^{n}\K_{1h}\left(\X_{1i}\right)[E\Y_{(1)}|\X_{i}-\tau_{1}(\x_{1})]+\frac{1}{\sqrt{nh_{1}^{k}}}\sum\limits_{i=1}^{n}\epsilon_{1i}\1(\D_{i}=1)\sum\limits_{j=1}^{n}\K_{1h}\left(\X_{1j}\right)(w_{ij}^{S_{1}}-w_{ji}^{S_{1}})\\
	&\quad+\frac{1}{\sqrt{nh_{1}^{k}}}\sum\limits_{i=1}^{n}\epsilon_{1i}\1(\D_{i}=1)\sum\limits_{j=1}^{n}\K_{1h}\left(\X_{1j}\right)w_{ji}^{S_{1}}\\
	&\quad+\frac{1}{\sqrt{nh_{1}^{k}}}\sum\limits_{i=1}^{n}\K_{1h}\left(\X_{1i}\right)\left[\frac{\frac{1}{nh_{4}^{r(1)}}\sum\limits_{j=1}^{n}\K_{4h}\left(\Z_{j}^{1}-\Z_{i}^{1}\right)\1(\D_{j}=1)E\Y_{(1)}|\X_{j}}{\frac{1}{nh_{4}^{r(1)}}\sum\limits_{j=1}^{n}\K_{4h}\left(\Z_{j}^{1}-\Z_{i}^{1}\right)\1(\D_{j}=1)}-E\Y_{(1)}|\X_{i}\right]\\
	&=:\I_{n,10}+\I_{n,11}+\I_{n,12}+\I_{n,13},
	\end{align*}}where
{\small\begin{align*}
	w_{ij}^{S_{1}}=\frac{\frac{1}{nh_{4}^{r(1)}}\K_{4h}\left(\Z_{i}^{1}-\Z_{j}^{1}\right)}{\frac{1}{nh_{4}^{r(1)}}\sum\limits_{i=1}^{n}\K_{4h}\left(\Z_{i}^{1}-\Z_{j}^{1}\right)\1(\D_{i}=1)}, \quad \epsilon_{1i}=\Y_{i}-E\Y_{(1)}|\X_{i}.
	\end{align*}}Similarly, we can decompose $\widehat{\m}_{0}(\be_{0}\ttop\X)-\tau_{0}(\x_{1})$ as  {\small\begin{align*}
	&\quad\frac{1}{\sqrt{nh_{1}^{k}}}\sum\limits_{i=1}^{n}\K_{1h}\left(\X_{1i}\right)\left[\widehat{\m}_{0}(\be_{0}\ttop\X)-\tau_{0}(\x_{1})\right]\\
	&=\frac{1}{\sqrt{nh_{1}^{k}}}\sum\limits_{i=1}^{n}\K_{1h}\left(\X_{1i}\right)[E\Y_{(0)}|\X_{i}-\tau_{0}(\x_{1})]+\frac{1}{\sqrt{nh_{1}^{k}}}\sum\limits_{i=1}^{n}\epsilon_{0i}\1(\D_{i}=0)\sum\limits_{j=1}^{n}\K_{1h}\left(\X_{1j}\right)(w_{ij}^{S_{2}}-w_{ji}^{S_{2}})\\
	&\quad+\frac{1}{\sqrt{nh_{1}^{k}}}\sum\limits_{i=1}^{n}\epsilon_{0i}\1(\D_{i}=0)\sum\limits_{j=1}^{n}\K_{1h}\left(\X_{1j}\right)w_{ji}^{S_{2}}\\
	&\quad+\frac{1}{\sqrt{nh_{1}^{k}}}\sum\limits_{i=1}^{n}\K_{1h}\left(\X_{1i}\right)\left[\frac{\frac{1}{nh_{4}^{r(0)}}\sum\limits_{j=1}^{n}\K_{4h}\left(\Z_{j}^{0}-\Z_{i}^{0}\right)\1(\D_{j}=0)E\Y_{(0)}|\X_{j}}{\frac{1}{nh_{4}^{r(0)}}\sum\limits_{j=1}^{n}\K_{4h}\left(\Z_{j}^{0}-\Z_{i}^{0}\right)\1(\D_{j}=0)}-E\Y_{(0)}|\X_{i}\right]\\
	&=:\I_{n,10}^{'}+\I_{n,11}^{'}+\I_{n,12}^{'}+\I_{n,13}^{'},
	\end{align*}}where
{\small\begin{align*}
	w_{ij}^{S_{2}}=\frac{\frac{1}{nh_{4}^{r(0)}}\K_{4h}\left(\Z_{i}^{0}-\Z_{j}^{0}\right)}{\frac{1}{nh_{4}^{r(0)}}\sum\limits_{i=1}^{n}\K_{4h}\left(\Z_{i}^{0}-\Z_{j}^{0}\right)\1(\D_{i}=0)}, \quad \epsilon_{0i}=\Y_{i}-E\Y_{(0)}|\X_{i}.
	\end{align*}It is easy  to show that $\I_{n,11}$, $\I_{n,11}^{'}$, $\I_{n,13}$ and $\I_{n,13}^{'}$ are $\oo_{p}(1)$ following the same arguments for proving that  $\I_{n,5}=\oo_{p}(1)$ and $\I_{n,7}=\oo_{p}(1)$ for Theorem \ref{NonParametric 1}. The details are omitted here. We now deal with $\I_{n,12}$ and $\I_{n,12}^{'}$. 
	{\lemm\label{small term}Suppose assumptions (C1) -- (C4), (A1) and (A5) -- (A7) are satisfied. Then, for each point $\x_{1}$ in the support of $\X_{1}$, \\
		(1) If $\X_{1}\subset^{k-q}\be_{1}\ttop\X$ and $\X_{1}\subset^{k-q}\be_{0}\ttop\X$ with $s_{4}(2-k/q)+k>0$ and $0<q\leq k$, we have
		\begin{align}\label{I16,1}
		\I_{n,12}=\oo_{p}(1),\quad \I_{n,12}^{'}=\oo_{p}(1).
		\end{align}The corresponding asymptotically linear representation is then
		\begin{align*}
		\sqrt{nh_{1}^{k}}\{\widehat{\tau}(\x_{1})-\tau(\x_{1})\}=\frac{1}{\sqrt{nh_{1}^{k}}}\frac{1}{f(\x_{1})}\sum\limits_{i=1}^{n}\{\m_{1}(\X_{i})-\m_{0}(\X_{i})-\tau(\x_{1})\}\K_{1h}\left(\X_{1i}\right)+\oo_{p}(1).
		\end{align*}
		(2) If $\X_{1}\subset\be_{1}\ttop\X$ and $\X_{1}\subset^{k-q}\be_{0}\ttop\X$ with $s_{4}(2-k/q)+k>0$ and $0<q\leq k$, we have
		\begin{align}\label{I16,2}
		\I_{n,12}=\frac{1}{\sqrt{nh_{1}^{k}}}\sum\limits_{i=1}^{n}\epsilon_{1i}\1(\D_{i}=1)\frac{\K_{1h}\left(\X_{1i}\right)}{p(\X_{i})}+\oo_{p}(1),\quad \I_{n,12}^{'}=\oo_{p}(1).
		\end{align}Then we have
		\begin{align*}
		\sqrt{nh_{1}^{k}}\{\widehat{\tau}(\x_{1})-\tau(\x_{1})\}=\frac{1}{\sqrt{nh_{1}^{k}}}\frac{1}{f(\x_{1})}\sum\limits_{i=1}^{n}\{\Psi_{2}(\X_{i},\Y_{i},\D_{i})-\tau(\x_{1})\}\K_{1h}\left(\X_{1i}\right)+\oo_{p}(1).
		\end{align*}
		(3) If $\X_{1}\subset^{k-q}\be_{1}\ttop\X$ and $\X_{1}\subset\be_{0}\ttop\X$ with $s_{4}(2-k/q)+k>0$ and $0<q\leq k$, we have
		\begin{align}\label{I16,2}
		\I_{n,12}=\oo_{p}(1),\quad \I_{n,12}^{'}=\frac{1}{\sqrt{nh_{1}^{k}}}\sum\limits_{i=1}^{n}\epsilon_{0i}\1(\D_{i}=0)\frac{\K_{1h}\left(\X_{1i}\right)}{p(\X_{i})}+\oo_{p}(1).
		\end{align}The corresponding asymptotically linear representation is
		\begin{align*}
		\sqrt{nh_{1}^{k}}\{\widehat{\tau}(\x_{1})-\tau(\x_{1})\}=\frac{1}{\sqrt{nh_{1}^{k}}}\frac{1}{f(\x_{1})}\sum\limits_{i=1}^{n}\{\Psi_{3}(\X_{i},\Y_{i},\D_{i})-\tau(\x_{1})\}\K_{1h}\left(\X_{1i}\right)+\oo_{p}(1).
		\end{align*}
		(4) If $\X_{1}\subset\be_{1}\ttop\X$ and $\X_{1}\subset\be_{0}\ttop\X$, we have
		\begin{align}\label{I16,2}
		\I_{n,12}=\frac{1}{\sqrt{nh_{1}^{k}}}\sum\limits_{i=1}^{n}\epsilon_{1i}\1(\D_{i}=1)\frac{\K_{1h}\left(\X_{1i}\right)}{p(\X_{i})}+\oo_{p}(1),\quad \I_{n,12}^{'}=\frac{1}{\sqrt{nh_{1}^{k}}}\sum\limits_{i=1}^{n}\epsilon_{0i}\1(\D_{i}=0)\frac{\K_{1h}\left(\X_{1i}\right)}{p(\X_{i})}+\oo_{p}(1).
		\end{align} We have
		\begin{align*}
		\sqrt{nh_{1}^{k}}\{\widehat{\tau}(\x_{1})-\tau(\x_{1})\}=\frac{1}{\sqrt{nh_{1}^{k}}}\frac{1}{f(\x_{1})}\sum\limits_{i=1}^{n}\{\Psi_{4}(\X_{i},\Y_{i},\D_{i})-\tau(\x_{1})\}\K_{1h}\left(\X_{1i}\right)+\oo_{p}(1).
		\end{align*}}
	\noindent {\sc Proof of Lemma \ref{small term}.} We need to show that $\I_{n,12}=\oo_{p}(1)$ if $\X_{1}\subset^{k-q}\be_{1}\ttop\X$ with $s_{4}(2-k/q)+k>0$ and $0<q\leq k$. Let $\X_{1}=v_{1}$, $\be_{1}\ttop\X=v_{2}$, and denote $\left(\frac{v_{1}-v_{1i}}{h_{4}},\frac{v_{2}-v_{2i}}{h_{4}}\right)$ as $(t_{1},t_{2})$. We have
	{\small\begin{align*}
		&\quad E\left\{\sum\limits_{j=1}^{n}\K_{1h}\left(\X_{1j}\right)w_{ji}^{S_{1}}\Bigg|\X_{i}\right\}\\
		&=\frac{1+\oo_{p}(1)}{h_{4}^{r(1)}f(v_{2i})p(v_{2i})}\int \K_{4}\left(\frac{v_{2j}-\be_{1}\ttop\X_{i}}{h_{4}}\right)\K_{1}\left(\frac{v_{1j}-\X_{1}}{h_{1}}\right)f(v_{i})dv\\
		&=h_{4}\frac{1+\oo_{p}(1)}{f(v_{2i})p(v_{2i})}\int \K_{4}(t_{2})\K_{1}\left(\frac{v_{1i}-\X_{1}}{h_{1}}+t_{1}\frac{h_{4}}{h_{1}}\right)f_{12}(v_{1i}+h_{4}t_{1},v_{2i}+h_{4}t_{2})dt_{1}dt_{2}\\
		&=h_{4}^{q}\K_{1}\left(\frac{v_{1i}-\X_{1}}{h_{1}}\right)\frac{f_{12}(v_{1i},v_{2i})}{f(v_{2i})p(v_{2i})}\int \K_{4}(t_{2})dt_{1}dt_{2}\\
		&\quad+\frac{h_{4}^{q+1}}{h_{1}}\K_{1}^{'}\left(\frac{v_{1i}-\X_{1}}{h_{1}}\right)\frac{f_{12}(v_{1i},v_{2i})}{f(v_{2i})p(v_{2i})}\int t_{1}\K_{4}(t_{2})dt_{1}dt_{2}+\oo_{p}\left(\frac{h_{4}^{2}}{h_{1}}\right),
		\end{align*}}where $f_{12}(v_{1i},v_{2i})$ is the joint density function of $(\X_{1},\be_{1}\ttop\X)$. Under assumptions (A5) -- (A7), we have
	{\small\begin{align*} E\left\{\sum\limits_{j=1}^{n}\K_{1h}\left(\X_{1j}\right)w_{ji}^{S_{1}}\Bigg|\X_{i}\right\}=C_{2}h_{4}^{q}\K_{1h}\left(\X_{1i}\right)\frac{f_{12}(\X_{1i},\be_{1}\ttop\X_{i})}{f(\X_{i})p(\be_{1}\ttop\X_{i})}+\OO_{p}\left(\frac{h_{4}^{q+1}}{h_{1}}\right)=\OO_{p}\left(h_{4}^{q}+\frac{h_{4}^{q+1}}{h_{1}}\right).
		\end{align*}}Hence, under assumptions (A6), (A7), $s_{4}(2-k/q)+k>0$ and $0<q\leq k$, $$\frac{1}{\sqrt{nh_{1}^{k}}}\sum\limits_{i=1}^{n}\epsilon_{1i}\1(\D_{i}=1)\sum\limits_{j=1}^{n}\K_{1h}\left(\X_{1j}\right)w_{ji}^{S_{1}}=\frac{1}{\sqrt{n}}\sum\limits_{i=1}^{n}\epsilon_{1i}\1(\D_{i}=1)\OO_{p}\left(\frac{h_{4}^{q}}{h_{1}^{k/2}}+\frac{h_{4}^{q+1}}{h_{1}^{k/2+1}}\right)=\oo_{p}(1).$$ Analogously, we get $\I_{n,12}^{'}=\oo_{p}(1)$ if $\X_{1}\subset^{k-q}\be_{0}\ttop\X$. Next, we prove that $$\I_{n,12}=\frac{1}{\sqrt{nh_{1}^{k}}}\sum\limits_{i=1}^{n}\epsilon_{1i}\1(\D_{i}=1)\frac{\K_{1h}\left(\X_{1i}\right)}{p(\X_{i})}+\oo_{p}(1),$$ if $\X_{1}\subset\be_{1}\ttop\X$. As that case that  $\X_{1}\subset\be_{1}\ttop\X$ is similar to  that  $\X_{1}\subset\X$ in nonparametric case, then parallelling  to derive equation \eqref{conditional expectaion}, we get the desired result. Similarly, we have $\I_{n,12}^{'}=\frac{1}{\sqrt{nh_{1}^{k}}}\sum\limits_{i=1}^{n}\epsilon_{0i}\1(\D_{i}=0)\frac{\K_{1h}\left(\X_{1i}\right)}{p(\X_{i})}+\oo_{p}(1)$ if $\X_{1}\subset\be_{0}\ttop\X$. The proof for Lemma \ref{small term} is concluded.
	
	\noindent {\sc Proof of Corollary \ref{NonParametric 2}.} Consider the case where 
	$\X_{1}\not\subset\widetilde{\X}\in\mR^{q}$. Similarly as before, we derive that
	{\small\begin{equation}\label{nonparametric estimator 2}
		\begin{aligned}
		&\quad\sqrt{nh_{1}^{k}}(\widehat{\tau}(\x_{1})-\tau(\x_{1}))\\
		&=\frac{\frac{1}{\sqrt{nh_{1}^{k}}}\sum\limits_{i=1}^{n}\K_{1h}\left(\X_{1i}\right)\left[\widehat{\m}_{1}(\widetilde{\X}_{i})-\tau_{1}(\x_{1})\right]}{\frac{1}{nh_{1}^{k}}\sum\limits_{i=1}^{n}\K_{1h}(\X_{1i})}-\frac{\frac{1}{\sqrt{nh_{1}^{k}}}\sum\limits_{i=1}^{n}\K_{1h}\left(\X_{1i}\right)\left[\widehat{\m}_{0}(\widetilde{\X}_{i})-\tau_{0}(\x_{1})\right]}{\frac{1}{nh_{1}^{k}}\sum\limits_{i=1}^{n}\K_{1h}(\X_{1i})},
		\end{aligned}\end{equation}}where
	\begin{align*}
	\widehat{\m}_{1}(\widetilde{\X}_{i})=\frac{\frac{1}{nh_{2}^{q}}\sum\limits_{j=1}^{n}\K_{2h}\left(\widetilde{\X}_{j}-\widetilde{\X}_{i}\right)\Y_{1j}\1(\D_{j}=1)}{\frac{1}{nh_{2}^{q}}\sum\limits_{j=1}^{n}\K_{2h}\left(\widetilde{\X}_{j}-\widetilde{\X}_{i}\right)\1(\D_{j}=1)},\quad
	\widehat{\m}_{0}(\widetilde{\X}_{i})=\frac{\frac{1}{nh_{2}^{q}}\sum\limits_{j=1}^{n}\K_{2h}\left(\widetilde{\X}_{j}-\widetilde{\X}_{i}\right)\Y_{0j}\1(\D_{j}=0)}{\frac{1}{nh_{2}^{q}}\sum\limits_{j=1}^{n}\K_{2h}\left(\widetilde{\X}_{j}-\widetilde{\X}_{i}\right)\1(\D_{j}=0)}.
	\end{align*} Some similar calculations lead to  $\widehat{\m}_{1}(\widetilde{\X}_{i})-\tau_{1}(\x_{1})$.
	{\small\begin{align*}
		&\quad\frac{1}{\sqrt{nh_{1}^{k}}}\sum\limits_{i=1}^{n}\K_{1h}\left(\X_{1i}\right)\left[\widehat{\m}_{1}(\widetilde{\X}_{i})-\tau_{1}(\x_{1})\right]\\
		&=\frac{1}{\sqrt{nh_{1}^{k}}}\sum\limits_{i=1}^{n}\K_{1h}\left(\X_{1i}\right)\left[E\Y_{(1)}|\X_{i}-\tau_{1}(\x_{1})\right]+\frac{1}{\sqrt{nh_{1}^{k}}}\sum\limits_{i=1}^{n}\epsilon_{1i}\1(\D_{i}=1)\sum\limits_{j=1}^{n}\K_{1h}\left(\X_{1j}\right)(w_{ij}^{N_{1}}-w_{ji}^{N_{1}})\\
		&\quad+\frac{1}{\sqrt{nh_{1}^{k}}}\sum\limits_{i=1}^{n}\epsilon_{1i}\1(\D_{i}=1)\sum\limits_{j=1}^{n}\K_{1h}\left(\X_{1j}\right)w_{ji}^{N_{1}}\\
		&\quad+\frac{1}{\sqrt{nh_{1}^{k}}}\sum\limits_{i=1}^{n}\K_{1h}\left(\X_{1i}\right)\left[\frac{\frac{1}{nh_{2}^{q}}\sum\limits_{j=1}^{n}\K_{2h}\left(\widetilde{\X}_{j}-\widetilde{\X}_{i}\right)\1(\D_{j}=1)E\Y_{(1)}|\X_{j}}{\frac{1}{nh_{2}^{q}}\sum\limits_{j=1}^{n}\K_{2h}\left(\widetilde{\X}_{j}-\widetilde{\X}_{i}\right)\1(\D_{j}=1)}-E\Y_{(1)}|\X_{i}\right]\\
		&=:\I_{n,14}+\I_{n,15}+\I_{n,16}+\I_{n,17},
		\end{align*}}where
	\begin{align*}
	w_{ij}^{N_{1}}=\frac{\frac{1}{nh_{2}^{q}}\K_{2h}\left(\widetilde{\X}_{i}-\widetilde{\X}_{j}\right)}{\frac{1}{nh_{2}^{q}}\sum\limits_{i=1}^{n}\K_{2h}\left(\widetilde{\X}_{i}-\widetilde{\X}_{j}\right)\1(\D_{i}=1)}.
	\end{align*}Then we can prove that $\I_{n,15}$ and $\I_{n,17}$ are $\oo_{p}(1)$ by the same arguments as those used to handle $\I_{n,5}$ and $\I_{n,7}$ for proving  Theorem \ref{NonParametric 1}. Owing to $\X_{1}\not\subset\widetilde{\X}$, similar arguments for proving Lemma~\ref{small term} implies that $\I_{n,16}=\oo_{p}(1)$. The proof for Corollary \ref{NonParametric 2} is concluded.
	\eop
	
	\noindent {\sc Proof of Corollary \ref{Comparallel to general case}.} From the proof for Theorem \ref{NonParametric 1}, we can see that
	{\small\begin{align*}
		E\left\{\sum\limits_{j=1}^{n}\K_{1h}\left(\X_{1j}\right)w_{ji}^{N}\Bigg|\X_{i}\right\}=\OO_{p}\left(h_{2}+\frac{h_{2}^{s_{2}}}{h_{1}^{s_{2}}}\right),
		\end{align*}}by the condition $\sqrt{nh_{1}^{k}}\left(h_{2}^{s}+\sqrt{\log(n)/nh_{2}^{p}}\right)=\oo(1)$. Then NRCATE shares the same asymptotic distribution as PRCATE. For SRCATE, we can use similar arguments to show the same result. The proof is finished.
	\eop


\bibliographystyle{elsarticle-harv} 
\bibliography{cate_ref}
\begin{table}[h]
	\centering
	\caption{{The distribution of $\sqrt{nh_{1}}[\widehat{\tau}(\x_{1})-\tau(\x_{1})]$ for model 1}}\resizebox{13cm}{8cm}{
		\begin{tabular}{|c|c|cccc|cccc|cccc|cccc|}
			\hline
			&       & \multicolumn{8}{c|}{n=200}                                     &        \multicolumn{8}{|c|}{n=500} \\
			\hline
			&       & \multicolumn{1}{c}{OR} & \multicolumn{1}{c}{PR} & \multicolumn{1}{c}{SR} & \multicolumn{1}{c}{NR} & \multicolumn{1}{c}{N} & \multicolumn{1}{c}{S} & \multicolumn{1}{c}{P} & \multicolumn{1}{c}{O} &       \multicolumn{1}{c}{OR} & \multicolumn{1}{c}{PR} & \multicolumn{1}{c}{SR} & \multicolumn{1}{c}{NR} & \multicolumn{1}{c}{N} & \multicolumn{1}{c}{S} & \multicolumn{1}{c}{P} & \multicolumn{1}{c|}{O} \\
			\hline
			&       & \multicolumn{16}{c|}{$h_{1}=0.05n^{-1/9}$, $h_{4}=0.6n^{-1/4}$, $h_{2}=0.4n^{-1/4}$} \\
			\hline
			& -0.4  & 0.178 & 0.224 & 0.213 & 0.223 & 0.352 & 0.361 & 0.396 & 0.407 & 0.188 & 0.224 & 0.208 & 0.225 & 0.371 & 0.388 & 0.404 & 0.409 \\
			& -0.2  & 0.182 & 0.192 & 0.181 & 0.196 & 0.351 & 0.365 & 0.377 & 0.380 & 0.186 & 0.193 & 0.191 & 0.213 & 0.368 & 0.383 & 0.389 & 0.395 \\
			\multicolumn{1}{|l|}{SD} & 0     & 0.199 & 0.210 & 0.231 & 0.232 & 0.420 & 0.440 & 0.460 & 0.491 & 0.198 & 0.205 & 0.206 & 0.217 & 0.415 & 0.430 & 0.466 & 0.476 \\
			& 0.2   & 0.208 & 0.216 & 0.248 & 0.243 & 0.466 & 0.476 & 0.503 & 0.525 & 0.195 & 0.203 & 0.231 & 0.226 & 0.423 & 0.438 & 0.484 & 0.509 \\
			& 0.4   & 0.195 & 0.215 & 0.239 & 0.236 & 0.377 & 0.395 & 0.415 & 0.426 & 0.202 & 0.222 & 0.250 & 0.247 & 0.364 & 0.372 & 0.415 & 0.432 \\
			\hline
			& -0.4  & 0.005 & 0.011 & -0.032 & 0.001 & -0.024 & -0.001 & -0.004 & -0.006 & 0.021 & 0.026 & -0.097 & 0.011 & 0.022 & 0.055 & 0.043 & 0.037 \\
			& -0.2  & -0.002 & 0.006 & 0.094 & 0.043 & -0.011 & 0.011 & 0.014 & 0.017 & -0.005 & -0.003 & 0.119 & 0.034 & -0.036 & -0.008 & -0.013 & -0.013 \\
			\multicolumn{1}{|l|}{BIAS} & 0     & 0.005 & 0.013 & 0.040 & 0.032 & -0.033 & 0.004 & -0.007 & 0.012 & 0.007 & 0.007 & 0.057 & 0.030 & -0.026 & 0.006 & -0.004 & 0.007 \\
			& 0.2   & 0.005 & 0.009 & 0.008 & 0.013 & -0.006 & 0.035 & 0.001 & 0.014 & 0.003 & 0.001 & 0.005 & 0.007 & -0.030 & -0.004 & -0.005 & 0.006 \\
			& 0.4   & 0.006 & 0.004 & -0.014 & -0.008 & 0.033 & 0.066 & 0.041 & 0.027 & 0.015 & 0.013 & 0.000 & 0.008 & 0.032 & 0.038 & 0.017 & 0.012 \\
			\hline
			& -0.4  & 0.032 & 0.050 & 0.047 & 0.050 & 0.125 & 0.130 & 0.157 & 0.165 & 0.036 & 0.051 & 0.052 & 0.051 & 0.138 & 0.153 & 0.165 & 0.169 \\
			& -0.2  & 0.033 & 0.037 & 0.042 & 0.040 & 0.124 & 0.133 & 0.142 & 0.145 & 0.035 & 0.037 & 0.051 & 0.047 & 0.137 & 0.147 & 0.152 & 0.156 \\
			\multicolumn{1}{|l|}{MSE} & 0     & 0.040 & 0.044 & 0.055 & 0.055 & 0.177 & 0.194 & 0.212 & 0.241 & 0.039 & 0.042 & 0.046 & 0.048 & 0.173 & 0.185 & 0.217 & 0.226 \\
			& 0.2   & 0.043 & 0.047 & 0.061 & 0.059 & 0.217 & 0.228 & 0.253 & 0.276 & 0.038 & 0.041 & 0.054 & 0.051 & 0.180 & 0.192 & 0.234 & 0.259 \\
			& 0.4   & 0.038 & 0.046 & 0.057 & 0.056 & 0.143 & 0.160 & 0.174 & 0.182 & 0.041 & 0.049 & 0.062 & 0.061 & 0.133 & 0.140 & 0.173 & 0.187 \\
			\hline
			&       & \multicolumn{16}{c|}{ $h_{1}=0.05n^{-1/9}$, $h_{4}=0.5n^{-1/4}$, $h_{2}=0.4n^{-1/4}$} \\
			\hline
			& -0.4  & 0.195 & 0.235 & 0.219 & 0.227 & 0.358 & 0.371 & 0.403 & 0.392 & 0.181 & 0.220 & 0.212 & 0.225 & 0.365 & 0.381 & 0.406 & 0.405 \\
			& -0.2  & 0.191 & 0.198 & 0.196 & 0.211 & 0.386 & 0.398 & 0.413 & 0.410 & 0.192 & 0.201 & 0.194 & 0.214 & 0.373 & 0.386 & 0.406 & 0.408 \\
			\multicolumn{1}{|l|}{SD} & 0     & 0.199 & 0.206 & 0.214 & 0.216 & 0.391 & 0.415 & 0.429 & 0.435 & 0.196 & 0.209 & 0.219 & 0.230 & 0.418 & 0.436 & 0.466 & 0.484 \\
			& 0.2   & 0.202 & 0.207 & 0.235 & 0.231 & 0.440 & 0.455 & 0.495 & 0.525 & 0.203 & 0.209 & 0.231 & 0.227 & 0.419 & 0.431 & 0.468 & 0.493 \\
			& 0.4   & 0.207 & 0.222 & 0.248 & 0.245 & 0.375 & 0.380 & 0.429 & 0.441 & 0.196 & 0.212 & 0.231 & 0.229 & 0.361 & 0.370 & 0.416 & 0.426 \\
			\hline
			& -0.4  & 0.011 & 0.019 & -0.043 & 0.012 & 0.023 & 0.046 & 0.038 & 0.034 & 0.015 & 0.003 & -0.126 & -0.011 & -0.008 & 0.024 & 0.005 & 0.000 \\
			& -0.2  & 0.000 & 0.001 & 0.081 & 0.035 & -0.033 & -0.011 & -0.002 & -0.006 & 0.011 & 0.009 & 0.126 & 0.045 & -0.021 & 0.002 & -0.002 & -0.001 \\
			\multicolumn{1}{|l|}{BIAS} & 0     & -0.012 & -0.016 & 0.013 & 0.006 & -0.033 & 0.000 & -0.009 & -0.003 & 0.009 & 0.013 & 0.064 & 0.038 & -0.012 & 0.010 & 0.014 & 0.027 \\
			& 0.2   & -0.003 & -0.008 & -0.008 & -0.004 & -0.041 & -0.014 & -0.035 & -0.019 & -0.009 & -0.004 & -0.002 & -0.001 & -0.019 & -0.008 & -0.009 & 0.007 \\
			& 0.4   & -0.007 & -0.010 & -0.025 & -0.022 & 0.017 & 0.037 & 0.030 & 0.026 & 0.017 & 0.019 & 0.010 & 0.015 & 0.055 & 0.055 & 0.046 & 0.047 \\
			\hline
			& -0.4  & 0.038 & 0.056 & 0.050 & 0.051 & 0.129 & 0.140 & 0.164 & 0.155 & 0.033 & 0.048 & 0.061 & 0.051 & 0.133 & 0.145 & 0.165 & 0.164 \\
			& -0.2  & 0.037 & 0.039 & 0.045 & 0.046 & 0.150 & 0.159 & 0.171 & 0.168 & 0.037 & 0.040 & 0.053 & 0.048 & 0.139 & 0.149 & 0.165 & 0.167 \\
			\multicolumn{1}{|l|}{MSE} & 0     & 0.040 & 0.043 & 0.046 & 0.047 & 0.154 & 0.172 & 0.184 & 0.189 & 0.039 & 0.044 & 0.052 & 0.054 & 0.175 & 0.190 & 0.217 & 0.235 \\
			& 0.2   & 0.041 & 0.043 & 0.055 & 0.053 & 0.195 & 0.207 & 0.246 & 0.276 & 0.041 & 0.044 & 0.053 & 0.051 & 0.176 & 0.186 & 0.219 & 0.243 \\
			& 0.4   & 0.043 & 0.049 & 0.062 & 0.061 & 0.141 & 0.146 & 0.185 & 0.195 & 0.039 & 0.045 & 0.053 & 0.053 & 0.133 & 0.140 & 0.175 & 0.184 \\
			\hline
			&       & \multicolumn{16}{c|}{ $h_{1}=0.05n^{-1/9}$, $h_{4}=0.5n^{-1/4}$, $h_{2}=0.45n^{-1/4}$} \\
			\hline
			& -0.4  & 0.183 & 0.222 & 0.218 & 0.219 & 0.357 & 0.375 & 0.405 & 0.398 & 0.191 & 0.214 & 0.207 & 0.215 & 0.376 & 0.390 & 0.400 & 0.408 \\
			& -0.2  & 0.195 & 0.203 & 0.186 & 0.197 & 0.360 & 0.366 & 0.380 & 0.384 & 0.186 & 0.196 & 0.183 & 0.198 & 0.364 & 0.372 & 0.386 & 0.391 \\
			\multicolumn{1}{|l|}{SD} & 0     & 0.193 & 0.206 & 0.214 & 0.217 & 0.441 & 0.453 & 0.474 & 0.486 & 0.193 & 0.201 & 0.207 & 0.211 & 0.432 & 0.442 & 0.478 & 0.491 \\
			& 0.2   & 0.200 & 0.213 & 0.237 & 0.232 & 0.460 & 0.476 & 0.516 & 0.525 & 0.194 & 0.202 & 0.230 & 0.227 & 0.479 & 0.489 & 0.526 & 0.541 \\
			& 0.4   & 0.198 & 0.220 & 0.241 & 0.239 & 0.407 & 0.414 & 0.453 & 0.460 & 0.211 & 0.231 & 0.257 & 0.255 & 0.406 & 0.408 & 0.455 & 0.474 \\
			\hline
			& -0.4  & 0.005 & 0.000 & -0.064 & -0.009 & -0.013 & 0.010 & 0.010 & 0.007 & -0.004 & -0.004 & -0.130 & -0.016 & -0.006 & 0.019 & 0.009 & 0.020 \\
			& -0.2  & 0.001 & -0.002 & 0.079 & 0.049 & -0.044 & -0.029 & -0.024 & -0.022 & -0.002 & -0.004 & 0.118 & 0.041 & -0.025 & -0.007 & -0.007 & -0.004 \\
			\multicolumn{1}{|l|}{BIAS} & 0     & 0.003 & 0.000 & 0.022 & 0.017 & -0.029 & -0.010 & -0.018 & -0.005 & 0.008 & 0.006 & 0.056 & 0.034 & -0.034 & -0.014 & -0.003 & -0.008 \\
			& 0.2   & 0.016 & 0.013 & 0.013 & 0.018 & -0.034 & -0.001 & -0.026 & -0.016 & 0.000 & -0.001 & 0.002 & 0.006 & -0.030 & -0.026 & -0.022 & -0.028 \\
			& 0.4   & 0.004 & -0.001 & -0.014 & -0.015 & 0.014 & 0.035 & 0.013 & 0.000 & 0.021 & 0.021 & 0.009 & 0.010 & 0.030 & 0.036 & 0.008 & -0.003 \\
			\hline
			& -0.4  & 0.034 & 0.049 & 0.052 & 0.048 & 0.128 & 0.141 & 0.164 & 0.159 & 0.037 & 0.046 & 0.060 & 0.046 & 0.142 & 0.152 & 0.160 & 0.167 \\
			& -0.2  & 0.038 & 0.041 & 0.041 & 0.041 & 0.131 & 0.134 & 0.145 & 0.148 & 0.035 & 0.038 & 0.047 & 0.041 & 0.133 & 0.139 & 0.149 & 0.153 \\
			\multicolumn{1}{|l|}{MSE} & 0     & 0.037 & 0.042 & 0.046 & 0.047 & 0.195 & 0.205 & 0.225 & 0.236 & 0.037 & 0.040 & 0.046 & 0.046 & 0.188 & 0.195 & 0.228 & 0.241 \\
			& 0.2   & 0.040 & 0.045 & 0.056 & 0.054 & 0.213 & 0.226 & 0.267 & 0.276 & 0.038 & 0.041 & 0.053 & 0.052 & 0.230 & 0.240 & 0.278 & 0.293 \\
			& 0.4   & 0.039 & 0.048 & 0.058 & 0.057 & 0.166 & 0.172 & 0.205 & 0.211 & 0.045 & 0.054 & 0.066 & 0.065 & 0.165 & 0.168 & 0.207 & 0.225 \\
			\hline
	\end{tabular}}%
	\label{tab12}%
\end{table}%

\begin{table}[h]
	\centering
	\caption{{The distribution of $\sqrt{nh_{1}}[\widehat{\tau}(\x_{1})-\tau(\x_{1})]$ for model 2}}\resizebox{13cm}{8cm}{
		\begin{tabular}{|c|c|cccc|cccc|cccc|cccc|}
			\hline
			&       & \multicolumn{8}{c|}{n=200}                                     &        \multicolumn{8}{|c|}{n=500} \\
			\hline
			&       & \multicolumn{1}{c}{OR} & \multicolumn{1}{c}{PR} & \multicolumn{1}{c}{SR} & \multicolumn{1}{c}{NR} & \multicolumn{1}{c}{N} & \multicolumn{1}{c}{S} & \multicolumn{1}{c}{P} & \multicolumn{1}{c}{O} &       \multicolumn{1}{c}{OR} & \multicolumn{1}{c}{PR} & \multicolumn{1}{c}{SR} & \multicolumn{1}{c}{NR} & \multicolumn{1}{c}{N} & \multicolumn{1}{c}{S} & \multicolumn{1}{c}{P} & \multicolumn{1}{c|}{O} \\
			\hline
			&       & \multicolumn{16}{c|}{ $h_{1}=0.02n^{-1/9}$, $h_{4}=0.2n^{-1/4}$, $h_{2}=0.15n^{-1/4}$} \\
			\hline
			& -0.4  & 0.384 & 0.386 & 0.391 & 0.409 & 0.950 & 1.103 & 1.101 & 1.098 & 0.341 & 0.348 & 0.356 & 0.376 & 0.959 & 1.066 & 1.061 & 1.079 \\
			& -0.2  & 0.389 & 0.395 & 0.399 & 0.430 & 0.968 & 1.100 & 1.106 & 1.114 & 0.360 & 0.362 & 0.365 & 0.396 & 0.962 & 1.091 & 1.088 & 1.099 \\
			\multicolumn{1}{|l|}{SD} & 0     & 0.386 & 0.388 & 0.393 & 0.419 & 1.000 & 1.165 & 1.141 & 1.136 & 0.373 & 0.376 & 0.380 & 0.397 & 0.940 & 1.077 & 1.053 & 1.064 \\
			& 0.2   & 0.379 & 0.378 & 0.381 & 0.406 & 1.011 & 1.175 & 1.122 & 1.116 & 0.357 & 0.361 & 0.368 & 0.398 & 0.998 & 1.140 & 1.120 & 1.121 \\
			& 0.4   & 0.384 & 0.390 & 0.412 & 0.435 & 1.011 & 1.150 & 1.105 & 1.103 & 0.390 & 0.394 & 0.413 & 0.438 & 1.045 & 1.182 & 1.129 & 1.157 \\
			\hline
			& -0.4  & 0.010 & 0.010 & 0.059 & 0.031 & -0.667 & -0.063 & 0.081 & 0.089 & 0.020 & 0.018 & 0.108 & 0.028 & -1.033 & -0.118 & 0.032 & 0.015 \\
			& -0.2  & -0.017 & -0.017 & 0.012 & 0.003 & -0.740 & -0.158 & -0.008 & 0.011 & -0.005 & -0.007 & 0.033 & -0.004 & -1.078 & -0.151 & -0.023 & -0.036 \\
			\multicolumn{1}{|l|}{BIAS} & 0     & -0.013 & -0.015 & -0.022 & -0.017 & -0.751 & -0.143 & -0.025 & -0.009 & 0.004 & 0.002 & -0.003 & 0.028 & -0.996 & -0.082 & 0.050 & 0.052 \\
			& 0.2   & -0.002 & -0.002 & -0.023 & 0.013 & -0.650 & 0.004 & 0.058 & 0.078 & -0.011 & -0.013 & -0.068 & -0.031 & -0.968 & -0.008 & 0.028 & 0.035 \\
			& 0.4   & 0.060 & 0.058 & 0.013 & 0.041 & -0.566 & 0.095 & 0.103 & 0.104 & 0.005 & 0.004 & -0.067 & -0.007 & -0.892 & 0.120 & 0.020 & 0.019 \\
			\hline
			& -0.4  & 0.148 & 0.149 & 0.156 & 0.168 & 1.348 & 1.220 & 1.218 & 1.213 & 0.117 & 0.121 & 0.139 & 0.142 & 1.987 & 1.150 & 1.127 & 1.165 \\
			& -0.2  & 0.152 & 0.157 & 0.160 & 0.185 & 1.483 & 1.234 & 1.222 & 1.240 & 0.129 & 0.131 & 0.134 & 0.157 & 2.087 & 1.213 & 1.184 & 1.209 \\
			\multicolumn{1}{|l|}{MSE} & 0     & 0.149 & 0.151 & 0.155 & 0.176 & 1.564 & 1.377 & 1.303 & 1.291 & 0.139 & 0.142 & 0.145 & 0.158 & 1.876 & 1.166 & 1.111 & 1.136 \\
			& 0.2   & 0.143 & 0.143 & 0.146 & 0.165 & 1.445 & 1.380 & 1.262 & 1.252 & 0.128 & 0.130 & 0.140 & 0.160 & 1.932 & 1.300 & 1.256 & 1.258 \\
			& 0.4   & 0.151 & 0.156 & 0.170 & 0.191 & 1.342 & 1.332 & 1.232 & 1.226 & 0.152 & 0.155 & 0.175 & 0.192 & 1.888 & 1.411 & 1.275 & 1.339 \\
			\hline
			&       & \multicolumn{16}{c|}{ $h_{1}=0.02n^{-1/9}$, $h_{4}=0.18n^{-1/4}$, $h_{2}=0.13n^{-1/4}$} \\
			\hline
			& -0.4  & 0.368 & 0.374 & 0.376 & 0.397 & 1.003 & 1.204 & 1.140 & 1.161 & 0.346 & 0.348 & 0.358 & 0.372 & 0.945 & 1.093 & 1.077 & 1.102 \\
			& -0.2  & 0.399 & 0.397 & 0.407 & 0.443 & 1.011 & 1.192 & 1.183 & 1.177 & 0.368 & 0.369 & 0.371 & 0.394 & 0.889 & 1.030 & 1.028 & 1.039 \\
			\multicolumn{1}{|l|}{SD} & 0     & 0.389 & 0.390 & 0.392 & 0.413 & 1.029 & 1.192 & 1.188 & 1.197 & 0.362 & 0.364 & 0.373 & 0.408 & 0.966 & 1.101 & 1.077 & 1.099 \\
			& 0.2   & 0.387 & 0.387 & 0.395 & 0.411 & 1.048 & 1.254 & 1.207 & 1.198 & 0.328 & 0.330 & 0.333 & 0.376 & 0.966 & 1.097 & 1.078 & 1.104 \\
			& 0.4   & 0.391 & 0.398 & 0.420 & 0.432 & 1.041 & 1.202 & 1.131 & 1.147 & 0.370 & 0.377 & 0.390 & 0.391 & 1.019 & 1.172 & 1.089 & 1.114 \\
			\hline
			& -0.4  & 0.023 & 0.027 & 0.079 & 0.019 & -0.811 & -0.173 & -0.012 & -0.030 & -0.023 & -0.020 & 0.070 & -0.015 & -1.169 & -0.194 & -0.027 & -0.018 \\
			& -0.2  & 0.003 & 0.005 & 0.033 & 0.015 & -0.754 & -0.101 & 0.052 & 0.049 & 0.005 & 0.007 & 0.046 & 0.002 & -1.101 & -0.141 & 0.039 & 0.050 \\
			\multicolumn{1}{|l|}{BIAS} & 0     & 0.008 & 0.009 & 0.011 & 0.019 & -0.781 & -0.109 & -0.014 & -0.005 & 0.000 & 0.001 & -0.007 & 0.001 & -1.103 & -0.121 & 0.013 & 0.021 \\
			& 0.2   & 0.027 & 0.025 & -0.006 & 0.031 & -0.653 & 0.054 & 0.122 & 0.119 & 0.003 & 0.003 & -0.046 & 0.016 & -1.060 & -0.011 & 0.054 & 0.052 \\
			& 0.4   & 0.023 & 0.020 & -0.025 & 0.018 & -0.588 & 0.124 & 0.157 & 0.128 & 0.014 & 0.013 & -0.058 & 0.008 & -0.986 & 0.057 & 0.027 & 0.021 \\
			\hline
			& -0.4  & 0.136 & 0.141 & 0.147 & 0.158 & 1.664 & 1.479 & 1.300 & 1.348 & 0.120 & 0.121 & 0.133 & 0.139 & 2.258 & 1.232 & 1.161 & 1.215 \\
			& -0.2  & 0.159 & 0.158 & 0.166 & 0.196 & 1.592 & 1.431 & 1.402 & 1.388 & 0.135 & 0.136 & 0.139 & 0.155 & 2.002 & 1.081 & 1.058 & 1.083 \\
			\multicolumn{1}{|l|}{MSE} & 0     & 0.151 & 0.152 & 0.154 & 0.171 & 1.668 & 1.433 & 1.412 & 1.433 & 0.131 & 0.133 & 0.139 & 0.166 & 2.150 & 1.227 & 1.159 & 1.207 \\
			& 0.2   & 0.150 & 0.151 & 0.156 & 0.170 & 1.525 & 1.577 & 1.473 & 1.450 & 0.108 & 0.109 & 0.113 & 0.141 & 2.057 & 1.203 & 1.165 & 1.222 \\
			& 0.4   & 0.154 & 0.159 & 0.177 & 0.187 & 1.429 & 1.460 & 1.304 & 1.331 & 0.137 & 0.142 & 0.155 & 0.153 & 2.011 & 1.376 & 1.187 & 1.241 \\
			\hline
			&       & \multicolumn{16}{c|}{$h_{1}=0.02n^{-1/9}$, $h_{4}=0.15n^{-1/4}$, $h_{2}=0.1n^{-1/4}$} \\
			\hline
			& -0.4  & 0.364 & 0.370 & 0.379 & 0.409 & 0.997 & 1.172 & 1.140 & 1.164 & 0.358 & 0.363 & 0.370 & 0.389 & 0.970 & 1.134 & 1.109 & 1.140 \\
			& -0.2  & 0.388 & 0.392 & 0.406 & 0.436 & 1.042 & 1.225 & 1.227 & 1.230 & 0.364 & 0.362 & 0.365 & 0.408 & 0.901 & 1.086 & 1.054 & 1.054 \\
			\multicolumn{1}{|l|}{SD} & 0     & 0.396 & 0.398 & 0.413 & 0.446 & 0.992 & 1.180 & 1.166 & 1.161 & 0.371 & 0.374 & 0.382 & 0.417 & 0.919 & 1.113 & 1.077 & 1.084 \\
			& 0.2   & 0.388 & 0.389 & 0.397 & 0.436 & 1.029 & 1.254 & 1.161 & 1.182 & 0.369 & 0.370 & 0.374 & 0.389 & 1.021 & 1.199 & 1.148 & 1.168 \\
			& 0.4   & 0.375 & 0.379 & 0.403 & 0.430 & 1.151 & 1.360 & 1.261 & 1.280 & 0.364 & 0.370 & 0.386 & 0.409 & 1.049 & 1.243 & 1.132 & 1.160 \\
			\hline
			& -0.4  & -0.001 & 0.002 & 0.047 & 0.008 & -0.838 & -0.202 & -0.010 & -0.019 & -0.013 & -0.010 & 0.087 & 0.000 & -1.255 & -0.212 & -0.034 & -0.021 \\
			& -0.2  & 0.008 & 0.012 & 0.038 & 0.020 & -0.872 & -0.245 & -0.067 & -0.058 & 0.021 & 0.022 & 0.061 & 0.018 & -1.145 & -0.086 & 0.120 & 0.121 \\
			\multicolumn{1}{|l|}{BIAS} & 0     & 0.022 & 0.023 & 0.023 & 0.032 & -0.850 & -0.196 & -0.036 & -0.030 & -0.001 & -0.003 & -0.014 & -0.005 & -1.265 & -0.190 & -0.023 & -0.024 \\
			& 0.2   & -0.007 & -0.007 & -0.036 & -0.014 & -0.839 & -0.140 & -0.053 & -0.042 & 0.007 & 0.005 & -0.047 & -0.001 & -1.213 & -0.103 & 0.006 & -0.007 \\
			& 0.4   & 0.011 & 0.007 & -0.036 & 0.000 & -0.759 & -0.075 & -0.013 & -0.021 & -0.005 & -0.009 & -0.082 & -0.013 & -1.191 & -0.073 & -0.075 & -0.103 \\
			\hline
			& -0.4  & 0.133 & 0.137 & 0.146 & 0.167 & 1.695 & 1.414 & 1.299 & 1.355 & 0.129 & 0.132 & 0.144 & 0.151 & 2.517 & 1.330 & 1.230 & 1.300 \\
			& -0.2  & 0.150 & 0.154 & 0.166 & 0.191 & 1.846 & 1.562 & 1.510 & 1.515 & 0.133 & 0.132 & 0.137 & 0.167 & 2.121 & 1.187 & 1.125 & 1.126 \\
			\multicolumn{1}{|l|}{MSE} & 0     & 0.158 & 0.159 & 0.171 & 0.200 & 1.706 & 1.431 & 1.360 & 1.350 & 0.138 & 0.140 & 0.146 & 0.174 & 2.445 & 1.275 & 1.161 & 1.176 \\
			& 0.2   & 0.150 & 0.151 & 0.159 & 0.190 & 1.764 & 1.592 & 1.351 & 1.400 & 0.136 & 0.137 & 0.142 & 0.151 & 2.512 & 1.447 & 1.319 & 1.363 \\
			& 0.4   & 0.141 & 0.144 & 0.164 & 0.185 & 1.901 & 1.856 & 1.590 & 1.638 & 0.132 & 0.137 & 0.156 & 0.168 & 2.519 & 1.551 & 1.286 & 1.356 \\
			\hline
	\end{tabular}}%
	\label{tab22}%
\end{table}%

\begin{table}[h]
	\centering
	\caption{{The distribution of $\sqrt{nh_{1}}[\widehat{\tau}(\x_{1})-\tau(\x_{1})]$ for model 3}}\resizebox{13cm}{8cm}{
		\begin{tabular}{|c|c|cccc|cccc|cccc|cccc|}
			\hline
			&       & \multicolumn{8}{c|}{n=200}                                     &        \multicolumn{8}{|c|}{n=500} \\
			\hline
			&       & \multicolumn{1}{c}{OR} & \multicolumn{1}{c}{PR} & \multicolumn{1}{c}{SR} & \multicolumn{1}{c}{NR} & \multicolumn{1}{c}{N} & \multicolumn{1}{c}{S} & \multicolumn{1}{c}{P} & \multicolumn{1}{c}{O} &       \multicolumn{1}{c}{OR} & \multicolumn{1}{c}{PR} & \multicolumn{1}{c}{SR} & \multicolumn{1}{c}{NR} & \multicolumn{1}{c}{N} & \multicolumn{1}{c}{S} & \multicolumn{1}{c}{P} & \multicolumn{1}{c|}{O} \\
			\hline
			&       & \multicolumn{16}{c|}{ $h_{1}=0.02n^{-1/9}$, $h_{4}=0.18n^{-1/4}$, $h_{2}=0.16n^{-1/4}$} \\
			\hline
			& -0.4  & 0.320 & 0.321 & 0.328 & 0.319 & 0.516 & 0.532 & 0.570 & 0.548 & 0.306 & 0.310 & 0.315 & 0.307 & 0.502 & 0.501 & 0.501 & 0.499 \\
			& -0.2  & 0.341 & 0.343 & 0.352 & 0.345 & 0.477 & 0.477 & 0.514 & 0.526 & 0.298 & 0.301 & 0.304 & 0.292 & 0.476 & 0.472 & 0.497 & 0.501 \\
			\multicolumn{1}{|l|}{SD} & 0     & 0.301 & 0.306 & 0.312 & 0.304 & 0.450 & 0.458 & 0.487 & 0.495 & 0.292 & 0.296 & 0.299 & 0.290 & 0.484 & 0.466 & 0.512 & 0.514 \\
			& 0.2   & 0.320 & 0.320 & 0.322 & 0.313 & 0.493 & 0.486 & 0.521 & 0.514 & 0.296 & 0.298 & 0.304 & 0.296 & 0.470 & 0.455 & 0.489 & 0.488 \\
			& 0.4   & 0.306 & 0.314 & 0.319 & 0.312 & 0.501 & 0.525 & 0.534 & 0.530 & 0.301 & 0.305 & 0.308 & 0.296 & 0.473 & 0.477 & 0.483 & 0.491 \\
			\hline
			& -0.4  & -0.023 & -0.025 & -0.028 & -0.044 & -0.038 & -0.029 & 0.001 & 0.003 & 0.026 & 0.027 & 0.025 & 0.005 & 0.027 & 0.021 & 0.050 & 0.054 \\
			& -0.2  & 0.026 & 0.022 & 0.020 & 0.035 & 0.004 & 0.009 & -0.006 & -0.009 & -0.006 & -0.006 & -0.003 & 0.010 & 0.013 & 0.022 & 0.003 & 0.004 \\
			\multicolumn{1}{|l|}{BIAS} & 0     & 0.003 & -0.001 & 0.011 & 0.035 & 0.048 & 0.051 & 0.019 & 0.022 & 0.010 & 0.011 & 0.013 & 0.042 & 0.020 & 0.043 & -0.014 & -0.015 \\
			& 0.2   & 0.003 & 0.000 & 0.001 & 0.014 & 0.015 & 0.026 & 0.008 & 0.010 & -0.012 & -0.011 & -0.010 & 0.009 & 0.033 & 0.044 & 0.023 & 0.022 \\
			& 0.4   & -0.004 & -0.006 & -0.011 & -0.018 & -0.023 & -0.012 & 0.011 & 0.013 & 0.001 & 0.001 & -0.004 & -0.023 & -0.044 & -0.049 & -0.010 & -0.008 \\
			\hline
			& -0.4  & 0.103 & 0.104 & 0.109 & 0.103 & 0.267 & 0.284 & 0.324 & 0.301 & 0.094 & 0.097 & 0.100 & 0.094 & 0.252 & 0.252 & 0.253 & 0.252 \\
			& -0.2  & 0.117 & 0.118 & 0.124 & 0.120 & 0.227 & 0.227 & 0.265 & 0.277 & 0.089 & 0.091 & 0.092 & 0.085 & 0.227 & 0.223 & 0.247 & 0.251 \\
			\multicolumn{1}{|l|}{MSE} & 0     & 0.091 & 0.094 & 0.098 & 0.094 & 0.205 & 0.212 & 0.237 & 0.246 & 0.085 & 0.088 & 0.090 & 0.086 & 0.234 & 0.219 & 0.262 & 0.265 \\
			& 0.2   & 0.102 & 0.102 & 0.104 & 0.098 & 0.244 & 0.237 & 0.272 & 0.264 & 0.088 & 0.089 & 0.093 & 0.088 & 0.222 & 0.209 & 0.240 & 0.239 \\
			& 0.4   & 0.093 & 0.098 & 0.102 & 0.097 & 0.251 & 0.276 & 0.285 & 0.281 & 0.091 & 0.093 & 0.095 & 0.088 & 0.226 & 0.230 & 0.234 & 0.241 \\
			\hline
			&       & \multicolumn{16}{c|}{$h_{1}=0.02n^{-1/9}$, $h_{4}=0.2n^{-1/4}$, $h_{2}=0.15n^{-1/4}$} \\
			\hline
			& -0.4  & 0.297 & 0.301 & 0.313 & 0.306 & 0.465 & 0.480 & 0.502 & 0.506 & 0.285 & 0.290 & 0.296 & 0.289 & 0.497 & 0.512 & 0.507 & 0.513 \\
			& -0.2  & 0.311 & 0.313 & 0.319 & 0.314 & 0.423 & 0.428 & 0.467 & 0.462 & 0.303 & 0.307 & 0.309 & 0.300 & 0.471 & 0.475 & 0.482 & 0.478 \\
			\multicolumn{1}{|l|}{SD} & 0     & 0.316 & 0.320 & 0.322 & 0.322 & 0.470 & 0.465 & 0.532 & 0.522 & 0.321 & 0.325 & 0.331 & 0.325 & 0.483 & 0.487 & 0.521 & 0.521 \\
			& 0.2   & 0.318 & 0.323 & 0.328 & 0.323 & 0.460 & 0.462 & 0.501 & 0.502 & 0.291 & 0.298 & 0.301 & 0.297 & 0.468 & 0.471 & 0.484 & 0.485 \\
			& 0.4   & 0.301 & 0.305 & 0.306 & 0.305 & 0.489 & 0.493 & 0.530 & 0.516 & 0.310 & 0.311 & 0.312 & 0.307 & 0.518 & 0.536 & 0.524 & 0.522 \\
			\hline
			& -0.4  & 0.003 & 0.002 & -0.001 & -0.013 & -0.019 & -0.004 & 0.024 & 0.025 & -0.003 & -0.004 & -0.007 & -0.025 & -0.043 & -0.025 & -0.018 & -0.015 \\
			& -0.2  & -0.025 & -0.023 & -0.023 & -0.013 & -0.023 & -0.024 & -0.043 & -0.044 & 0.011 & 0.011 & 0.015 & 0.029 & 0.007 & 0.016 & -0.004 & -0.004 \\
			\multicolumn{1}{|l|}{BIAS} & 0     & 0.001 & 0.003 & 0.009 & 0.028 & 0.019 & 0.028 & -0.009 & -0.014 & 0.014 & 0.015 & 0.023 & 0.051 & 0.050 & 0.059 & 0.012 & 0.018 \\
			& 0.2   & 0.008 & 0.009 & 0.017 & 0.025 & 0.024 & 0.029 & 0.017 & 0.011 & 0.010 & 0.011 & 0.015 & 0.026 & 0.019 & 0.032 & 0.012 & 0.011 \\
			& 0.4   & -0.010 & -0.009 & -0.014 & -0.025 & -0.055 & -0.048 & -0.010 & -0.012 & -0.004 & -0.004 & -0.010 & -0.025 & -0.034 & -0.013 & -0.009 & -0.008 \\
			\hline
			& -0.4  & 0.088 & 0.090 & 0.098 & 0.094 & 0.217 & 0.230 & 0.253 & 0.257 & 0.081 & 0.084 & 0.088 & 0.084 & 0.248 & 0.263 & 0.257 & 0.264 \\
			& -0.2  & 0.097 & 0.099 & 0.102 & 0.099 & 0.179 & 0.183 & 0.220 & 0.216 & 0.092 & 0.094 & 0.095 & 0.091 & 0.222 & 0.226 & 0.232 & 0.229 \\
			\multicolumn{1}{|l|}{MSE} & 0     & 0.100 & 0.103 & 0.104 & 0.105 & 0.221 & 0.217 & 0.283 & 0.272 & 0.103 & 0.106 & 0.110 & 0.108 & 0.236 & 0.241 & 0.271 & 0.272 \\
			& 0.2   & 0.101 & 0.104 & 0.108 & 0.105 & 0.212 & 0.215 & 0.251 & 0.252 & 0.085 & 0.089 & 0.091 & 0.089 & 0.220 & 0.223 & 0.234 & 0.235 \\
			& 0.4   & 0.091 & 0.093 & 0.094 & 0.094 & 0.243 & 0.246 & 0.281 & 0.267 & 0.096 & 0.097 & 0.097 & 0.095 & 0.270 & 0.288 & 0.275 & 0.273 \\
			\hline
			&       & \multicolumn{16}{c|}{ $h_{1}=0.02n^{-1/9}$, $h_{4}=0.18n^{-1/4}$, $h_{2}=0.15n^{-1/4}$} \\
			\hline
			& -0.4  & 0.309 & 0.313 & 0.323 & 0.316 & 0.515 & 0.541 & 0.531 & 0.519 & 0.282 & 0.286 & 0.286 & 0.283 & 0.487 & 0.493 & 0.495 & 0.491 \\
			& -0.2  & 0.315 & 0.315 & 0.326 & 0.312 & 0.495 & 0.534 & 0.528 & 0.539 & 0.299 & 0.300 & 0.304 & 0.294 & 0.480 & 0.481 & 0.488 & 0.486 \\
			\multicolumn{1}{|l|}{SD} & 0     & 0.322 & 0.328 & 0.329 & 0.324 & 0.473 & 0.485 & 0.526 & 0.531 & 0.304 & 0.304 & 0.308 & 0.299 & 0.472 & 0.476 & 0.493 & 0.487 \\
			& 0.2   & 0.302 & 0.303 & 0.311 & 0.302 & 0.460 & 0.433 & 0.480 & 0.485 & 0.305 & 0.307 & 0.308 & 0.305 & 0.465 & 0.466 & 0.485 & 0.487 \\
			& 0.4   & 0.306 & 0.311 & 0.314 & 0.309 & 0.476 & 0.500 & 0.506 & 0.514 & 0.283 & 0.285 & 0.286 & 0.284 & 0.449 & 0.471 & 0.467 & 0.464 \\
			\hline
			& -0.4  & -0.011 & -0.011 & -0.016 & -0.026 & -0.026 & 0.002 & 0.021 & 0.016 & 0.016 & 0.014 & 0.006 & -0.009 & -0.020 & -0.018 & 0.012 & 0.009 \\
			& -0.2  & -0.012 & -0.012 & -0.012 & -0.001 & -0.014 & -0.011 & -0.024 & -0.027 & 0.007 & 0.005 & 0.006 & 0.022 & 0.009 & 0.025 & 0.006 & 0.004 \\
			\multicolumn{1}{|l|}{BIAS} & 0     & -0.012 & -0.012 & -0.005 & 0.018 & 0.058 & 0.048 & 0.024 & 0.023 & -0.002 & -0.004 & -0.003 & 0.024 & 0.021 & 0.034 & -0.014 & -0.021 \\
			& 0.2   & 0.017 & 0.018 & 0.021 & 0.034 & 0.040 & 0.042 & 0.032 & 0.029 & -0.007 & -0.009 & -0.009 & 0.007 & -0.004 & 0.001 & -0.016 & -0.019 \\
			& 0.4   & 0.012 & 0.013 & 0.012 & -0.002 & -0.043 & -0.017 & 0.003 & 0.004 & -0.009 & -0.011 & -0.011 & -0.029 & -0.044 & -0.032 & -0.005 & -0.007 \\
			\hline
			& -0.4  & 0.096 & 0.098 & 0.105 & 0.101 & 0.266 & 0.293 & 0.282 & 0.270 & 0.080 & 0.082 & 0.082 & 0.080 & 0.237 & 0.244 & 0.245 & 0.241 \\
			& -0.2  & 0.100 & 0.099 & 0.106 & 0.098 & 0.245 & 0.285 & 0.279 & 0.291 & 0.090 & 0.090 & 0.092 & 0.087 & 0.230 & 0.232 & 0.238 & 0.237 \\
			\multicolumn{1}{|l|}{MSE} & 0     & 0.104 & 0.108 & 0.108 & 0.106 & 0.227 & 0.238 & 0.278 & 0.282 & 0.092 & 0.092 & 0.095 & 0.090 & 0.223 & 0.228 & 0.244 & 0.238 \\
			& 0.2   & 0.091 & 0.092 & 0.097 & 0.093 & 0.213 & 0.189 & 0.232 & 0.236 & 0.093 & 0.095 & 0.095 & 0.093 & 0.217 & 0.217 & 0.235 & 0.238 \\
			& 0.4   & 0.094 & 0.097 & 0.098 & 0.096 & 0.228 & 0.251 & 0.256 & 0.265 & 0.080 & 0.081 & 0.082 & 0.081 & 0.204 & 0.223 & 0.218 & 0.216 \\
			\hline
	\end{tabular}}
	\label{tab32}%
\end{table}%
\end{document}